\input amssym.def  
\input epsf

\font\piccolo=cmr10 at 8truept
\font \piccobo=cmbx10 at 8truept
\def \hb{\hfill \break}  
\def \pf{\noindent{$\underline{\hbox{Proof}}$.\ }}
\def \sym{{\rm{Sym}}}
\def \rk{{\rm{rk}}} 
\def \codim{{\rm{codim\,}}} 
\def \al{\alpha} 
\def \be{\beta}   
\def \xx{\overline X}
\def \aa{\overline A}
\def \bb{\overline B}
 
\def \dd{\overline D}  
\def \ss{\overline S}
\def \vv{\overline V} 
\def \uu{\overline U}
\def \rr{\overline R}

\def \pp{\overline P}
\def \qq{\overline Q}
\def \ff{\overline f}

\def \ddd{\overline {\dd}}
\def \del{\overline {\Delta}}
\def \wdd{\widetilde {D}}
\def \wss{\widetilde {S}}

\def \wll{\widetilde {L}} 
\def \wmm{\widetilde {M}}  
\def \wxx{\widetilde {X}}

\def \wff{\widetilde {f}}
\def \cvd{\hfill {$\diamond$} \smallskip}
\def \hs#1{{\Bbb F}_{#1}}
\def \Hs {Hirzebruch surface }
\def \grc {\hbox{\bf \piccobo GRC}}
\def\mapright#1,#2 {\smash{\mathop
{{\hbox to #2pt {\rightarrowfill}}}\limits^{#1}} }
\def\mapleft#1,#2 {\smash{\mathop
{{\hbox to #2pt {\leftarrowfill}}}\limits_{#1}} }
\def\mapdown#1{\vbox to 15 pt{} 
   \Big\downarrow\rlap{$\vcenter{\hbox{$\scriptstyle#1$}}$}}
\def\lmapdown#1{\vbox to 15 pt{} 
   \llap{$\vcenter{\hbox{$\scriptstyle#1$}}$}\Big\downarrow}

\def\lmapup#1{\vbox to 15 pt{}
    \llap{$\vcenter{\hbox{$\scriptstyle#1$}}$}\Big\uparrow}

\centerline{\bf Stratification of the moduli space of four--gonal curves}
\medskip
\centerline{Michela Brundu - Gianni Sacchiero}

\bigskip
\centerline{\bf Abstract}
\medskip
\noindent
{\piccolo 
Let $X$ be a smooth irreducible projective curve of genus $g$ and gonality 4.
We show that the canonical model of $X$ is contained in a uniquely
defined surface, ruled by conics, whose geometry is deeply related to  that of  $X$.
This surface allows us to
 define four invariants of
$X$ and hence to stratify the moduli space of  four--gonal curves by
means of closed irreducible subvarieties whose dimensions we compute.}


\bigskip
\noindent
{\bf Contents}

\medskip
\noindent

0. Preliminaries

1. The surface $S$ of minimum degree,  ruled by conics and containing
$X_K$

2.  Birational models of $X_K \subset S$

3. Singularities of a birational model $X_0$

4.  ``Standard'' birational models of $S$ and $X_K$

5. Bounds of the invariants $\lambda$ and $t$

6. Geometric meaning of the invariant $\lambda$ 

7. Bounds of the invariants $a$ and $b$ 

8. Existence of curves with given invariants  $\lambda, a, b$ in the case $t=0$

9. Proof of the key--lemma

10. Moduli spaces of 4-gonal curves with $t=0$

11. Moduli spaces of 4-gonal curves with $t\ge 1$

\bigskip
\noindent
{\bf Introduction}
\medskip
Let $X$ be a smooth irreducible curve of genus $g$ and gonality
$\gamma$, i.e. $\gamma$ is the minimal degree of a
base--point--free linear series on $X$. 
Let ${\cal M}_g$ denote the moduli space of  curves of genus
$g$ and  ${\cal M}_{g,\gamma} \subset {\cal M}_g$ denote the
variety parametrizing the $\gamma$-gonal curves; it is
well--known that  ${\cal M}_{g,\gamma}$ is an irreducible
variety of dimension $2g+2\gamma -5$, as far as $2 \le \gamma \le {g \over 2}+1$
(see [{\bf 13}] and [{\bf 1}]). 

The structure of ${\cal M}_{g,\gamma}$ is completely understood in the cases
$\gamma =2$ (hyperelliptic curves) and $\gamma =3$ (trigonal
curves). In this paper we are interested in the study of four--gonal
curves. Let us briefly recall the setting in the
trigonal case.

Let $K$ denote the canonical divisor on $X$ and $X_K \subset
\Bbb P^{g-1}$ be the canonical model of $X$. {}From the
Geometric Riemann--Roch Theorem, any trigonal divisor spans a
line in $\Bbb P^{g-1}$, therefore $X_K$ is contained in a
rational normal ruled surface,   $R$ say. It is clear that $R$ is of
the form $\Bbb P({\cal O}(m) \oplus {\cal O}(g-2-m))$;
assuming  $m \le g-2-m$, the integer $m$ is uniquely
determined and it is called the {\it Maroni invariant} of $X$.

Set ${\cal M}_{g,3}(m)$ the variety parametrizing the trigonal
curves of Maroni invariant not bigger than $m$. The following
fact holds:

\medskip
\noindent
{\bf Theorem. } {\sl If ${{g-4} \over 3} \le m < {{g-2} \over 2}$ (resp. $m={{g-2} \over
2}$) then 
${\cal M}_{g,3}(m)$ is a locally closed subset of  ${\cal M}_{g,3}$ of
dimension $g+2m+4$ (resp. $2g+1$).}
\medskip
\noindent
(See [{\bf 14}], Proposition 1.2).

\medskip
One can see that for each curve of genus $g\ge 5$ of Maroni invariant
$m$ there exists a unique linear series $g^1_\lambda$, where $\lambda$ is the
minimum integer bigger than $3$ and $\lambda =g-m-1$. Hence $\lambda$ is
uniquely determined by $m$ and the above filtration of ${\cal M}_{g,3}$
given by the varieties ${\cal M}_{g,3}(m)$ can be rewritten in terms
of $\lambda$.

In general, it seems interesting to find ``good invariants"  arising from the geometric properties of
$\gamma$--gonal canonical curves, in order to obtain an analogous
stratification of the moduli space ${\cal M}_{g,\gamma}$.

As in the trigonal case, one can introduce the rational normal scroll
$V$, whose fibres are the $(\gamma-2)$--planes spanned by the
$\gamma$--gonal divisor on $X$. Clearly 
$V=\Bbb P({\cal O}(a_1) \oplus \cdots 
\oplus {\cal O}(a_{\gamma -1}))$, where $a_1 + \cdots +
a_{\gamma-1} = g-\gamma+1$; in this way the integers
$a_1, \dots, a_{\gamma-2}$ play the role of the Maroni invariant
$m$ in the trigonal case.

\medskip
In this paper we focus  on $4$--gonal curves. We show that in
the volume $V=\Bbb P({\cal O}(a) \oplus {\cal O}(b) \oplus
{\cal O}(c))$ there exists an (almost always) uniquely determined ``minimal"
surface, ruled by conics, containing $X_K$. \hb
Such  a surface $S$ gives rise to other two invariants: on one hand, one
defines the number $t$ which is the uniquely determined
invariant of a suitable geometrically ruled surface birationally equivalent to $S$. On the
other hand, analyzing the embedding of $X$ in $S$, we obtain
another number $\lambda >4$ which turns out to be the minimum
degree of  a linear series on $X$ different from the gonal one. \hb
Comparing the configuration $X_K \subset S \subset V$ in
the $4$--gonal case with the analougous situation $X_K \subset R$
of the trigonal case, it is clear that the invariant $m$ has been
replaced, in some sense, by $a$, $b$ and $t$. Finally, one can prove that 
$\lambda$ is now independent of $a,b$ and $t$; so a four--gonal
curve is determined by the four invariants $a,b,\lambda,t$.

\medskip
In Section 6 we describe the geometric meaning of $\lambda$, while, in Sections 5 and
7,  we find the ranges for the above invariants $\lambda, t$ and $a,b$, respectively.

If $t=0$ the cited ranges become:
$$
{{g+3}\over 3} \le \lambda \le 
{{g+3}\over 2}
\eqno(R_1)
$$
$$
a_{\min} \le a \le {{g -3} \over 3 }
\eqno (R_2)
$$
$$
g-\lambda -1 \le a+b \le {{2(g -3)} \over 3 }
\eqno (R_3)
$$
where
$$
a_{\min} = 
\left\{
\matrix{
\displaystyle{\Big\lceil {\lambda -4 \over 2} \Big\rceil} & \hbox{if} &
\lambda
\ge {2g+6\over5}
\cr 
\cr
\displaystyle{g-2\lambda   +1} & \hbox{if} & \lambda \le {2g+6\over5} \cr 
} \right. 
$$
In Section 8 (see Theorem 8.5) we then show that, if 
$(R_1),(R_2),(R_3)$ are satisfied,  there exists a $4$--gonal
curve of genus $g$ and invariants
$a,b,\lambda$ and $t=0$.

\medskip
Finally, in Section 10 we study the moduli spaces ${\cal M}_{g,4}$ of
$4$--gonal curves with $t=0$. Set ${\cal M}_g^{\lambda} \subset 
{\cal M}_{g,4}$  be the variety parametrizing the $4$--gonal curves of
invariant $\lambda$ and ${\cal M}_g^{\lambda}(a,b) \subset 
{\cal M}_g^{\lambda}$ the subvariety parametrizing the curves of
further invariants $a$ and $b$. We prove the following:

\medskip
\proclaim
Main Theorem. Let $g,\lambda,a,b$ be positive integers satisfying
$(R_1)$, $(R_2)$, $(R_3)$ and $g \ge 10$. Then:
\item {$i)$} There exists a stratification of the moduli space
${\cal M}_{g,4}$ of $4$--gonal curves given by:
$$
{\cal M}_{g,4}=
\overline{\cal M}^{\left\lceil{{g+2} \over 2}\right\rceil}_g
\supset 
\overline{\cal M}^{\left\lceil{{g} \over 2}\right\rceil}_g
\supset \cdots \supset
\overline{\cal M}^{\lambda}_{g}
\supset \cdots \supset
\overline{\cal M}^{\left\lceil{{g+3} \over 3}\right\rceil}_g
$$
and $\overline{\cal M}^{\lambda}_g$ are irreducible locally closed
subsets of dimension $g+2\lambda+1$, if $\lambda < \left\lceil{{g+2} \over 2}\right\rceil$.
\item {$ii)$} For each admissible $\lambda$, we can write:
$$
\overline{\cal M}^{\lambda}_g =
\bigcup_{a,b}\;
\overline{\cal M}^{\lambda}_g(a,b)
$$
where 
$\overline{\cal M}^{\lambda}_g(a,b)$ is a non--empty, irreducible
subvariety whose dimension is :
$$
\dim({\cal M}^{\lambda}_g(a,b))=
\cases{
2(2a+b+\lambda)+10-g-\epsilon-\tau-\xi, & if  $a \ge {{g-\lambda-1} \over 2}$ \cr
\cr
2(a+b)+\lambda +8-\epsilon-\xi, & if  $a < {{g-\lambda-1} \over 2}$ \cr
}
$$
where
$$
\epsilon := \cases{
0, & if $b<c$ \cr
1 , & if $a<b=c$ \cr
2 , & if $a=b=c$ \cr}
\quad , \quad
\tau := \cases{
0, & if $a<b$ \cr
1 , &  if $a=b$ \cr}
\quad \hbox{and} \quad
\xi := \cases{
1, & if $\lambda={{g+3} \over 2}$ \cr
0 , & otherwise \cr}.
$$
\par

\medskip
In Section 11 we briefly describe the moduli space of four--gonal curves of invariant $t\ge 1$.

\bigskip
We would like to thank Valentina Beorchia for many helpful discussions and suggestions and Gianfranco Casnati for several interesting remarks. We are also grateful to Simon Brain and Giovanni Landi for the warm support.

\vfill
\eject

\bigskip
\noindent
{\bf 0. Preliminaries}
\medskip
We say that a curve is $4$--gonal if it has a linear series 
$g^1_4$ but  no $g^1_d$, for any $d\le 3$. We also assume that such curve is not
bi--hyperelliptic (i.e. the degree four map on $\Bbb P^1$ does not
factorize through a hyperelliptic curve), in particular  that is not bielliptic.

Let $X$ be a $4$--gonal curve of genus $g$. In order to have a unique $g^1_4$ on $X$, we assume $g \ge 10$.

Denote by $\varphi_K: X \rightarrow X_K \subset \Bbb P^{g-1}$  the
canonical map associated to $X$ and by $X_K$ the canonical model of $X$. In general, if $Y$ is
a variety and
$D$ is a divisor on $Y$, we denote by $\varphi_D\, : Y \rightarrow 
\varphi_D(Y) \subset \Bbb P(H^0(Y,{\cal O}_Y(D)))$ the morphism
associated to $D$.

If $\Phi \in g_4^1$ is a $4$--gonal divisor, by
the Geometric Riemann--Roch Theorem (see [{\bf 2}], Ch. I, Sect. 2) we
have that:
$\dim\langle \varphi_K(\Phi)\rangle = 
\deg(\Phi)-h^0({\cal O}_X(\Phi)) = 2$; therefore
$$
V:= \bigcup_{\Phi \in g_4^1} \langle \varphi_K(\Phi)\rangle
\subset \Bbb P^{g-1}
$$
is a scroll, ruled by planes on $\Bbb P^1$, containing $X_K$.  
 Denote $\pi: V \longrightarrow \Bbb P^1$ the natural projection.

\medskip
Recall that a non degenerate variety $W \subset \Bbb P^r$ is said  to be {\it  projectively normal} if it is normal and, for any $k \in \Bbb N$, the homomorphism
$$
H^0(\Bbb P^r,  {\cal O}_{\Bbb P^r}(k)) \longrightarrow H^0(W, {\cal O}_W(k))
$$
induced by the exact sequence of sheaves
$$
0 \longrightarrow  {\cal I}_W \longrightarrow {\cal O}_{\Bbb P^r}
\longrightarrow {\cal O}_W \longrightarrow 0
$$
is surjective. \hb
We say that $W$ is {\it  linearly normal} if the homomorphism above is surjective for $k=1$. In particular, if $W$ is a non degenerate curve, then  it is linearly normal  if  and only if 
$h^0(W, {\cal O}_W(1))=h^0(\Bbb P^r, {\cal O}_{\Bbb P^r}(1))=r+1$.

\medskip
It is well--known that $X_K$ is projectively normal; so 
 $V$ is a rational normal scroll (hence  projectively normal as well). We then set
  $V= \Bbb P({\cal F})$, where
$\cal F$ is a vector bundle of rank 3 on $\Bbb P^1$ i.e.
$$
{\cal F} =  {\cal O}(a) \oplus 
{\cal O}(b) \oplus {\cal O}(c) ,
$$
for suitable non--negative integers $a \le b \le c$.  It is also well--known that, for any $k$, it holds:
$$
h^0(V, {\cal O}_V(k)) =
h^0(\Bbb P^1, \pi_{*}{\cal O}_V(k)) =
h^0(\Bbb P^1, \sym^k{\cal F}) \eqno(1)
$$
and  that the Riemann -- Roch
Theorem for any vector bundle $\cal G$ on $\Bbb P^1$ with non--negative splitting type gives:
$$
h^0(\Bbb P^1, {\cal G}) = \deg({\cal G}) + \rm {rk}({\cal G}).
\eqno(RR)
$$
{}From the two above relations, since $a,b,c \ge 0$, we then have: $h^0(V, {\cal O}_V(1))=h^0(\Bbb P^1, {\cal F})=
 \deg({\cal F}) + \rm {rk}({\cal F})$. Taking into account that $h^0(V, {\cal O}_V(1))=g$, we finally obtain:
$$
a+b+c = g-3. \eqno(2)
$$

\medskip
In the  following we will need some basic notations and facts about
ruled surfaces. \hb
We denote by $\hs{t}$ (where $t \ge 0$) the {\it
\Hs} of {\it invariant} $t$, i.e. the $\Bbb P^1$--bundle over $\Bbb P^1$ associated to the sheaf 
${\cal O}(-t) \oplus {\cal O}$ (here $\cal O$ means ${\cal O}_{\Bbb P^1}$).

If $1 \le a \le b$, a {\it rational ruled surface} $R_{a,b}$ is $\Bbb P({\cal O}(a) \oplus {\cal O}(b))$, naturally embedded
in
$\Bbb P^{a+b+1}$. Clearly, setting  $t:=b-a$, we have $R_{a,b} \cong \Bbb F_t$, so $t$ is 
the  {\it invariant} of 
$R_{a,b}$. 
\medskip

Let us recall the following well--known  facts  
(see [{\bf 11}], Ch. V, 2.9, 2.17 and 2.3):

\goodbreak

\medskip
\proclaim
Lemma 0.1. Let $\hs{t}$ be as before, $f$ its  generic  fibre and
 $C_0 = \Bbb P({\cal O}(-t)) \subset \hs{t}$. Then:
\item {$i)$} $C_0^2=-t$;
\item {$ii)$} if $U$ is any directrix (i.e. an irreducible unisecant curve) of $\, \hs{t}$, different from $C_0$, then
$U^2 \ge t$;
\item {$iii)$} if there exists a directrix $U$ of $R$ such
that $U^2=0$ then $t=0$, i.e. $\hs{0} \cong \Bbb P^1 \times \Bbb P^1$. \hb
Moreover, $t>0$ if and only if  $\, \hs{t}$ has exactly one unisecant curve (namely
$C_0$) having negative self--intersection.
\item {$iv)$}
$Num( \hs{t}) = \Bbb Z\langle C_0 \rangle \times \Bbb Z\langle f \rangle$.
\par

\medskip
Finally let us recall three classical formulas concerning ruled
surfaces and scrolls, due to C. Segre.

\medskip
\proclaim
Unisecants Formula. Let $R\subset \Bbb P^{r+1}$ be a ruled
surface $R$ of degree $r$ and invariant $t$ and let
$Un^d(R)$ be the variety of the unisecant
curves on $R$ having degree $d$ and self--intersection bigger
than $t$. Then the general element of $Un^d(R)$ is
irreducible and
$$
\dim(Un^d(R))= 2d+1-r.
\eqno (UF)
$$
\par
\pf Recall that, if $U \sim  C_0+nf$ is a unisecant curve on $R$,
where $U^2 > t$, then
$$
h^0(R, {\cal O}_R(U)) = 2n -t+2
\eqno(3)
$$
(see [{\bf 11}], Ch. V, 2.19). By appliying the equality $(3)$
to the hyperplane section $H$ of $R$, we
get 
$H \sim  C_0+ {{r+t} \over 2} f$. Take $D \in Un^d(R)$; since $D
\cdot H = d$, then $D \sim  C_0+ (d-{{r-t} \over 2})f$.
Therefore, since $D^2 >t$ by assumption, we can apply $(3)$ and
obtain the required formula. 
\cvd
\medskip
\noindent
The following Genus Formula $(GF)$ is a consequence of the
Adjuction Formula.

\medskip
\proclaim
Genus Formula. If $Y$ is a $q$-secant curve on a ruled surface 
$R \subset \Bbb P^r$, then 
$$
p_a(Y) = {{q-1} \over {2}} \; \big[2(\deg(Y)-1) - q
\deg(R)\big]. \eqno (GF)
$$
\par

The following relation $(IF)$, generalizing the
analogous property  for ruled surfaces, comes from the
Intersection Law on a scroll ([{\bf 8}], 8.3.14):
\medskip
\proclaim
Intersection Formula.
Let $W$ be a rational scroll ruled by $n$--planes and let
$C_1$ and $C_2$ be two subschemes of
$W$ meeting properly  and such that
$C_i$ is $m_i$--secant, for $i=1,2$ (i.e. $C_i$ meets the general
fibre of $W$ in a variety of degree $m_i$). Then the following
equality holds:
$$
\deg(C_1 \cdot C_2)=
m_1 \deg(C_2) +m_2 \deg(C_1) - m_1m_2 \deg(W).
\eqno (IF)
$$
\par

\medskip
Let us also recall the following notions:

\medskip
\noindent
{\bf Definition.}
Let  
$D$ be a very ample bisecant divisor on a \Hs  $\hs{}$; then the surface 
$S_0:= \varphi_D(\hs{})$ is said {\it geometrically  ruled by conics} (over $\Bbb P^1$).
Equivalently,  a projective surface
$S_0 \subset \Bbb P^N$ is  geometrically ruled  by conics if there exists
a surjective morphism $\pi: \; S_0 \longrightarrow \Bbb P^1$
such that the fibre $\pi^{-1}(y)$ is a smooth rational curve of degree $2$ for every point $y \in \Bbb P^1$ and $\pi$
admits a section. \hb
We say that  a projective surface
$S \subset \Bbb P^N$ is   {\it ruled by conics} (over  $\Bbb P^1$) if it is birational to a surface geometrically 
ruled by conics. Equivalently, if there exists a surjective morphism $\pi: \; S \longrightarrow
\Bbb P^1$ and  an open subset
$U\subseteq \Bbb P^1$ such that: 
\item \item {-} the fibre $\pi^{-1}(y)$ is a curve of degree $2$ and arithmetic genus $0$ for every point $y
\in \Bbb P^1$; 
\item \item {-} the fibre $\pi^{-1}(y)$ is  smooth  for every point $y \in U$; 
\item \item {-}  $\pi$ admits a section.
\par

\medskip
 The following classification of the degenerate fibres of a surface ruled
by conics is Thm. 2.4 (see also 1.13),  [{\bf 6}].

\goodbreak

\medskip
\proclaim 
Theorem 0.2. Let $S \subset \Bbb P^N$ be a projective surface ruled by conics over a smooth irreducible curve.
Then the degenerate fibres of $S$ are of one of the following types (where $n$ is an integer $\ge 3$  in the last two statements):
\item{-} $F_1$ is the union of two distinct lines and $S$ is smooth along $F_1$; 
\item{-}  $F_2(A)$ is the union of two distinct lines, whose common point is an
ordinary double point of $S$; 
\item{-}
 $F_2(D)$ is the union of two coincident lines, containing exactly two ordinary double
points of $S$;
\item{-}  $F_n(A)$ is the union of two distinct lines, whose common point is a rational
double point of type $(A_{n-1})$;
\item{-} $F_n(D)$ is  the union of two coincident lines, containing exactly
one rational double points of $S$; in particular, this point is of type $(A_3)$, if $n=3$, and of type
$(D_n)$, if $n \ge 4$.
\par

\medskip
Since any surface $S$ ruled by conics is birational to a surface $S_0$, geometrically ruled by conics, then $S$ can
be obtained from a suitable $S_0$ by a finite number of monoidal transformations. In particular, each singular fibre
of $S$ (as described in 0.2) arises in this way. Again in [{\bf 6}] we have studied
this situation, as summarized below.

Let $\hs{}$ and $D$ be as before and
$S_0= \varphi_D(\hs{})$ be a  surface geometrically ruled by conics  via the
morphism $\pi: S_0 \longrightarrow \Bbb P^1$. Consider a point  
$P_1 \in S_0$ and let $f_0:= \pi^{-1}(y)$ be the fibre of $S_0$ containing $P_1$.  Consider the
blow--up $\sigma_{P_1}$of $S_0$ at $P_1$ and the corresponding projection on $\Bbb P^1$,  
$\pi_1$ say: 
$$
\matrix{
Bl_{P_1} (S_0) :=   & S_1  & \mapright \sigma_{P_1},30 & S_0 \cr
&{\mapdown {\pi_1}} & & {\mapdown {\pi}}  \cr
&\Bbb P^1 &  & \Bbb P^1}
$$
Denote also by $f_1 := \pi_1^{-1}(y)$ the total transform of $f_0$ via $\sigma_{P_1}$. \hb
Take now $P_2 \in f_1$ and consider the corresponding blow--up $\sigma_{P_2}: \; S_2 \longrightarrow S_1$.
With obvious notations, we can iterate this construction and  obtain a sequence of blow--ups:
$$
\matrix{
\wss_0:= S_n & \mapright \sigma_{P_n},30 &  \cdots & \longrightarrow & S_2 & \mapright \sigma_{P_2},30
& \hfill  S_1  & \mapright \sigma_{P_1},30 & \hfill S_0 \cr
\hfill \cup \; & & & & \cup && \hfill \cup \; && \hfill \cup \;\cr
\wff_0:=f_n & & & & f_2 && P_2 \in f_1 &  &  P_1 \in f_0 \cr}
$$
where, for  any $i=1, \dots, n$, we define $P_i \in f_{i-1}$, $f_i:= \pi_i^{-1}(y)$ and 
$\pi_i: S_i:= Bl_{P_i}(S_{i-1}) \longrightarrow \Bbb P^1$ is the natural projection.

\medskip
\noindent
{\bf Definition.} With the above notation, we say that $f_n = \wff_0 \subset \wss_0$ is a fibre of {\it level $n$ over
$f_0$}.

\medskip
\noindent
Denoting by
 $\sigma$ the sequence of blowing--ups of $S_0$ defined above,  setting  $\wdd$  to be 
 the strict
transform of $D$ (very ample bisecant divisor on $S_0$) via $\sigma$ and 
$B$  the base locus of $\wdd$, then $S$ can be obtained in this way:
$$
\matrix{ 
\wss_0 &  \mapright \sigma,30 & S_0 \cr
\lmapdown {\varphi_{\wdd -B}} & \nearrow_{\rho}\cr
S
}
$$
where $\rho$ is defined as the birational map such that  the diagram is commutative.

\medskip
\noindent
{\bf Definition.} We say that the fibre $f \subset S$ is an {\it embedded
fibre of level
$n$} if
$$
n = \min_i \; \{ \hbox{there exists a blow--up $\sigma: \wss_0 \rightarrow S_0$ and a fibre 
$f_i \subset \wss_0$  of level
$i$ such that} \; f =
\varphi_{\wdd -B}(f_i)
\}.
$$

\medskip
\noindent
Again in [{\bf 6}], we noted that each fibre $f \subset S$ of type $F_n(A)$ or $F_n(D)$ is an embedded fibre of level 
$n$. 
 There  we also gave the following:

\medskip
\noindent
{\bf Definition.} Let $f^{(1)}, \dots, f^{(p)}$ be the degenerate fibres of $S$ and let $l_i$ be the
level of $f^{(i)}$, for $i=1, \dots, p$. 
If $\sum_{i=1}^p l_i =L$, we say that $S$ is of {\it level} $L$.

\medskip
 Moreover, we proved  that all the surfaces geometrically ruled by conics (briefly g.r.c.)
and giving rise -- by a minimal number of elementary transformations -- to a surface $S$ ruled by conics of level $L$,  are
exactly the elements of the following set:
$$
\eqalign
{\grc_L(S):= &\{ S_0 \;| \; S_0 \; \hbox{is a g.r.c. surface and $S$ can be obtained from it} \cr
 & \; \hbox{by a sequence of $L$
blow--ups and contractions}\}. \cr}
$$

\bigskip
\noindent
{\bf 1. The surface $S$ of minimum degree,  ruled by conics and containing
$X_K$}

\medskip
\noindent
Starting from the situation  $X_K \subset V \subset \Bbb P^{g-1}$, described at the
beginning of the previous section,  we will try to ``canonically'' define a surface
(ruled by conics) containing  $X_K$ and contained in $V$.

\medskip
\noindent
{\bf Notation.} As usual, if $n$ is a rational number, $[n]$
denotes the  greatest  integer smaller or equal than $n$, while 
$\lceil n \rceil$ denotes the  smallest integer bigger or equal than
$n$.

\medskip
\proclaim
Theorem 1.1. There exists a surface $S$ ruled by conics such that  
$X_K \subset S \subset V$ and  
$\displaystyle{
\deg(S) \le \left\lceil{{3g-8} \over 2} \right\rceil}$. 
Moreover, $S$ is unique unless 
$\displaystyle{
\deg(S) = {{3g-7} \over 2}}$; in this case,
$S$ varies in a pencil.
\par
\medskip
\pf 
Let us consider the vector space 
${\cal H}:=H^0(\Bbb P^{g-1}, {\cal I}_{X_K}(2))/H^0(\Bbb P^{g-1}, {\cal I}_V(2))$ and set $N:=\dim({\cal H})$; clearly, $\Sigma:=\Bbb P({\cal
H})$ parametrizes the hyperquadrics of
$\Bbb P^{g-1}$ containing $X_K$ but not containing $V$. \hb  
Let us recall that, if $W$ is a projectively
normal subvariety of $\Bbb P^{g-1}$, then we get the cohomology
exact sequence (see Section 0)
$$
0 \longrightarrow H^0({\cal I}_W(2)) 
  \longrightarrow H^0({\cal O}_{\Bbb P^{g-1}}(2))
  \longrightarrow H^0({\cal O}_W(2))
  \longrightarrow 0
$$
hence $h^0({\cal O}_{\Bbb P^{g-1}}(2))=h^0({\cal I}_W(2)) 
+h^0({\cal O}_W(2))$. Rewriting this equality for both $X_K$ and
$V$, we get 
$h^0({\cal I}_{X_K}(2)) +h^0({\cal O}_{X_K}(2))=
h^0({\cal O}_{\Bbb P^{g-1}}(2))= 
h^0({\cal I}_V(2)) +h^0({\cal O}_V(2))$, so 
$$
N=h^0({\cal I}_{X_K}(2))-h^0({\cal I}_V(2))=
h^0({\cal O}_V(2))-h^0({\cal O}_{X_K}(2)).
$$
In order to compute $N$, recall the relations $(1)$ and  $(RR)$ on the scroll 
$V =  \Bbb P({\cal F})$:
$$
h^0(V,{\cal O}_V(2)) = 
h^0(\Bbb P^1, \sym^2({\cal F}))= \deg(\sym^2({\cal F})) + \rk(\sym^2({\cal F})).
$$
Clearly, $\sym^2({\cal F})$  is a  free bundle of degree $4(a+b+c)$ and rank 6; therefore, 
from  $(2)$ we get:
$h^0({\cal O}_V(2))=4g-6$. \hb
On the other hand, by the Riemann--Roch Theorem
$h^0({\cal O}_{X_K}(2))=3(g-1)$.
Hence the above space $\Sigma$ of hyperquadrics is a
projective space of dimension
$$
N-1=h^0({\cal O}_V(2))-h^0({\cal O}_{X_K}(2))-1=g-4.
$$
For each $Q \in \Sigma \cong \Bbb P^{g-4}$, consider the scheme--theoretic
intersection
$$
Q \cdot V = \Big(\bigcup_{i=1, \dots, h_Q} F_i \Big) \cup S_Q
$$
where the $F_i$'s are the fibres of $V$ entirely contained in $Q$,
$h_Q \ge 0$ and $S_Q$ is a surface, which is ruled in
conics (since $Q$ intersects the general fibre $F$ of $V$ in a
conic passing through the four points of the divisor $\Phi
\subset F$) and contains $X_K$.\hb 
Note that $S_Q$ is irreducible; if not $S_Q = S_1 \cup S_2$, where the $S_i$'s were
ruled surfaces; but  $X_K \subset S_Q$ and it cannot be contained in a ruled surface since each
$4$--gonal divisor spans a plane. \hb
In order to find a quadric $\qq \in \Sigma$ such that
$\deg(S_{\qq})$ is minimum, it is enough to require that
the number $h_{\qq}$ is maximum. Note that a fibre $F$ is
contained in a quadric $Q \in \Sigma$ if $Q$ contains two points,
say $P_1$ and $P_2$, belonging to $F$ and such that the $0$-cycle
of $V$ of degree 6 given by $\Phi + P_1+P_2$ does not lie on a
conic.\hb
Since $\dim(\Sigma) = g-4$, we can impose  that the space
$\Sigma$ contains $\left[{{g-4} \over {2}}\right]$   pairs of
points. If each  such a pair  of points belongs to the same fibre (and
satisfies the above conditions), then we can find a $\qq \in
\Sigma$ containing $\left[{{g-4} \over {2}}\right]$ fibres. \hb
Clearly  $\qq$ could contain further fibres, hence
$$
\deg(S_{\qq}) \le 
\deg(\qq\cap V) -\left[{{g-4} \over {2}}\right] \le 
2(g-3) -
\left[{{g-4} \over {2}}\right] =
\left\lceil{{3g-8} \over 2} \right\rceil.
$$
This proves the existence of the required surface 
$S:=S_{\qq}$. \hb 
Concerning the uniqueness, let us assume
that there are two such surfaces, say $S_1$ and $S_2$. \hb
Since $X_K \subset (S_1 \cap S_2)$, from 
$(IF)$  we get: 
$$
2g-2 =\deg(X_K)\le \int (S_1 \cdot S_2)= 
2\deg(S_1) + 2\deg(S_2)- 4\deg(V)  . 
$$
This relation is verified if and only if 
$\deg(S_1)=\deg(S_2)= (3g-7)/2$. To complete the proof, just observe
that the linear system of the quadrics $\qq \in \Sigma$ containing
$\left[{{g-4}\over 2} \right]$ fibres has dimension
$$
\dim \Sigma - 2\left[{{g-4}\over 2} \right]=
g-4-2\left({{g-5}\over 2} \right)=1
$$
therefore there is a pencil of distinct surfaces $S_{\qq}$.
\cvd

\bigskip
The existence of such surface $S$ has been proved, using a different method,  also by Schreyer in [{\bf 12}],  Sect.6.

\medskip
\noindent
{\bf Notation.} {}From now on,  $f$ will denote  the general fibre of
$S$, so
$f$ is a conic lying on a plane $F = \langle f \rangle$. Moreover, if $T$ is a surface
ruled by conics, we will denote by $V_T$ the scroll whose fibres are the
planes spanned by these conics. For example, if $S$ is the surface defined in 1.1, the
scroll $V_S$ is exactly $V$.

\medskip
\noindent
{\bf Remark 1.2.}
The fibres of the ruled surface $S$ defined in 1.1 cannot
be  all singular. Otherwise, from 1.2, [{\bf 5}], the surface $S$ would be ruled by
lines on a hyperelliptic curve,  $Y$ say, via $\al: \; S \rightarrow Y$ and the ruling 
$\pi: \; S \rightarrow \Bbb P^1$ would factorize through $\al$. \hb
Hence, 
taking into account that the restriction $X_K \rightarrow Y$ of $\al$ has
degree two, we obtain that $X_K$ is bi-hyperelliptic,  contrary to  the assumption
made before on $X$.

\medskip
\noindent
{\bf Remark 1.3.}
The surface  $S$ introduced in 1.1 is then ruled by conics in the sense  of 
the preliminary Section.

\bigskip
\noindent
{\bf 2.  Birational models of $X_K \subset S$}

\bigskip
In this section we  shall  study a surface $S$ (not necessarily of minimum degree as
that one defined in 1.1) such that $S$ is ruled by conics and $X_K \subset S \subset V$,
where $V$ denotes as usual the $3$--dimensional scroll spanned by the four--gonal divisors
on $X_K$. \hb 
Note that, since $X_K$ is linearly normal, then $S\subset \Bbb P^{g-1}$  is
linearly normal. Moreover the scroll
$V=V_S$ is not a cone  (see the forthcoming Corollary 7.9), then 0.2 holds, so the
classification of the degenerate fibres of the surface $S$ is the one described there.
 \hb
  In Section 0 we have also summarized the results (contained in [{\bf 6}]) which
allow us to associate to a surface $S$, ruled by conics and of a certain level $L$, the
set
$\grc_L(S)$ consisting of all the g.r.c. surfaces linked to $S$ via a sequence of $L$ 
monoidal transformations. \hb 
Here  we are looking for the inverse procedure: how to recover the surface $S$ (and the
curve $X_K$) starting from a g.r.c. surface $S_0 \in \grc_L(S)$.

\medskip
\noindent
{\bf Notation.} 
Since each surface $S_0 \in \grc_L(S)$  is geometrically ruled by conics, it admits an
invariant $\tau_0 := t(S_0)$, in the sense that $S_0 \cong \hs{\tau_0}$. We denote by 
$X_{\tau_0} \subset \hs{\tau_0} \cong S_0$ the corresponding model of
 $X_K \subset S$. \hb 
 Since $X_{\tau_0} \subset \hs{\tau_0}$ is a four--secant curve, then
$$
X_{\tau_0} \sim 4C_0 + (\lambda_0 +\tau_0) f \eqno{(4)}
$$
where $C_0$ and $f$ are the generators of ${\rm Num}(\Bbb F_{\tau_0})$ (see 0.1) and
$\lambda_0$ is a suitable integer.  Moreover, denoting by
$p_a(C)$ the arithmetic genus of a curve
$C$, we set 
$$
\delta_{\tau_0}: = p_a(X_{\tau_0} ) - g.
$$ 
Note that, if all the singularities of $X_{\tau_0} $ are ordinary double points, then 
$\delta_{\tau_0} = \deg(Sing(X_{\tau_0}))$.

\medskip
\noindent 
{\bf Remark 2.1.} Let us recall the Adjunction Formula for the dualizing
sheaf $\omega_{X_R}$ of a curve $X_R$ on a smooth surface $R$ (see [{\bf 7}], Ch.1,
(1.5))
$$
\omega_{X_R} = {\cal K}_R \otimes {\cal O}_R(X_R)_{|X_R} \eqno{(5)}
$$
where ${\cal K}_R = {\cal O}_R(K_R)$ denotes the canonical sheaf of $R$. Taking the
degrees we then obtain: 
$$
2p_a(X_R) -2 = X_R \cdot(X_R + K_R). \eqno{(6)}
$$
In our situation $R=\hs{\tau_0}$ and $X_R=X_{\tau_0}$. Then  ${\cal
K}_{\hs{\tau_0}} = {\cal O}_{\hs{\tau_0}}(-2C_0-(\tau_0+2)f)$, so using $(4)$
 we obtain
$$
{\cal K}_{\hs{\tau_0}} \otimes {\cal O}_{\hs{\tau_0}}(X_{\tau_0}) = {\cal
O}_{\hs{\tau_0}}(2C_0+(\lambda_0-2)f).
$$
Hence from $(5)$ we can obtain the dualizing sheaf of the curve $X_{\tau_0}$ as:
$$
\omega_{X_{\tau_0}}=
{\cal O}_{\hs{\tau_0}}(2C_0+(\lambda_0-2)f)_{|X_{\tau_0}}.
$$
Finally,  taking into account that   $K_{\hs{\tau_0}} \sim -2C_0 - (\tau_0+2)f$, from
$(6)$ and $(4)$  we obtain
$$
2p_a(X_{\tau_0}) - 2 = 6 \lambda_0 -6\tau_0 -8.
$$

\medskip
\proclaim
Proposition 2.2. The following properties hold:
\item {$i)$} the arithmetic genus of $X_{\tau_0}$ is
$p_a(X_{\tau_0})=3(\lambda_0-\tau_0-1)$; 
\item {$ii)$} $\lambda_0 \ge \max \;\{ 3\tau_0, \tau_0+5\}$;
\item{$iii)$}  $\delta_{\tau_0}=3(\lambda_0 -\tau_0-1) -g$.
\par
\pf $i)$ Immediate from the last relation of 2.1. \hb
$ii)$ {}From [{\bf 11}], Ch. V, 2.18, since $X_{\tau_0}$ is irreducible,
then $\lambda_0 + \tau_0 \ge 4\tau_0$. Therefore
$\lambda_0 \ge 3 \tau_0$. On the other hand, $p_a(X_{\tau_0}) \ge g \ge 10$ by
assumption. Then, using $(i)$, we obtain $\lambda_0 \ge \tau_0+5$. \hb
$iii)$ It follows from $\delta_{\tau_0}= p_a(X_{\tau_0}) -g$ and from $(i)$.
\cvd

\medskip 
We  wish to describe  how  to recover the canonical model
$X_K$ starting from the chosen birational model $X_{\tau_0} \subset \hs{\tau_0} \cong S_0
\in \grc_L(S)$.  \hb
Since $X_0$ is the embedded model of $X_{\tau_0}$ obtained via the dualizing
sheaf $\omega_{X_{\tau_0}}$ (described before), then, in order to obtain $X_0$, we have to 
 embed $\hs{\tau_0}$ by  the sheaf 
${\cal O}_{\hs{\tau_0}}(2C_0+(\lambda_0-2)f)$ (see 2.1). Finally, we will project the obtained curve
$X_0$ from its singular points.

\medskip
\noindent
{\bf Remark 2.3.} Note first that $\lambda_0 -2 > 2\tau_0$. In fact, if $\tau_0 \le 2$
then $\lambda_0 > \tau_0 + 4 \ge 2\tau_0 +2$.  
If $\tau_0 \ge 3$, then  $\lambda_0 \ge 3\tau_0> 2\tau_0 +2$ (both arguments follow from 2.2, $(ii)$).
\hb 
Therefore (using [{\bf 11}], Ch. V, 2.18) the
linear system $|2C_0+ (\lambda_0 -2)f|$ is very ample on $\hs{\tau_0}$. Moreover, 
from [{\bf 4}], Prop.1.8, and from 2.2, $(iii)$ we get that 
$$
h^0\left(\hs{\tau_0},{\cal O}_{\hs{\tau_0}}(2C_0+ (\lambda_0 -2)f)\right)= g+\delta_{\tau_0}.
$$
Hence  there is an isomorphism
$$
\varphi: \hs{\tau_0}\;\; \mapright \cong,25 \; \; S_0 \subset \Bbb
P^{g-1+\delta_{\tau_0}},
\quad \hbox{where} \quad
\varphi = \varphi_{2C_0+(\lambda_0-2)f}
\quad \hbox{and} \quad
S_0:= \varphi(\hs{\tau_0}).
$$

\medskip
Clearly $S_0$ is a projective ruled surface, whose fibers  are all
smooth conics and  $X_0=\varphi(X_{\tau_0}) \subset S_0$, so we have the  commutative
diagrams:
$$
\matrix{
\hs{\tau_0} & \hfill \mapright{\varphi \atop \cong},30  & S_0 & \subset & \Bbb
P^{g-1+\delta_{\tau_0}}
\cr &   & \lmapup {\rho}  \mapdown{\pi} && \mapdown{\pi}
\cr {} & & S & \subset & \Bbb P^{g-1}\cr}
\qquad 
\hbox{and} 
\qquad
\matrix{
X_{\tau_0} & \hfill \mapright{\varphi_{|X_{\tau_0}} \atop \cong},30  & X_0 & \subset &
\hfill S_0
\cr &   & \lmapup {\rho}  \mapdown{\pi} && \mapdown{\pi}
\cr {} & & X_K & \subset & \phantom{\Big\{}S\cr}
$$
where $\pi$ (which is the inverse  of the map $\rho$) is exactly the desingularization morphism of $X_0$ or,
equivalently,  the linear projection centered in $\langle \Sigma \rangle$ is generated by  the singular points of
$X_0$ (possibly infinitely near). 

\medskip
\noindent
{\bf Remark 2.4.} Since there are at most  two  singular points on each fibre, then $\langle \Sigma \rangle$ meets $S_0$ in a zero--dimensional variety of degree $\delta_{\tau_0} $.
It is then clear that  $\delta_{\tau_0} = L$ and $\deg(S) = \deg(S_0)
- \delta_{\tau_0}$.

\bigskip
\noindent
{\bf 3. Singularities of a birational model $X_0$} 

\medskip
The purpose of this section is to describe all the possible singularities of
$X_0$. \hb
Recall that, from 2.3, the projection $\pi: X_0 \subset S_0 \longrightarrow X_K \subset
S$ is centered in the singular points of $X_0$ and the singular fibres of $S$ correspond
to the fibres of $S_0$ containing the singular points of $X_0$. Therefore it is enough to
examine the singular fibres of
$S$ and the  four--gonal divisor on each of them. \hb 
In order to do this,
let us focus  on one singular fibre $f$ of
$S$ and the corresponding fibre $f_0 \subset S_0$. 

\medskip
\noindent
{\bf Remark 3.1.}
Note  that the curve $X_K \subset S$ intersects
each fibre of $S$ in four points (the $4$--gonal divisor $\Phi \in g^1_4$). In
particular, $X_K$ meets also each singular fibre $f$ in four
points. If $f = l \cup m$ and $l\ne m$ then two of them belong to the line $l$ and two
 are on the other line $m$ (possibly  coinciding);  where this not the case, 
$X_K$ would have a trisecant line, hence a trigonal series (from
the Geometric Riemann--Roch Theorem). On the other hand, if $l=m$, then the support of $\Phi= X_K \cap f$ consists of two
points, possibly  coinciding.

\medskip
\noindent
{\bf Example 3.2.}
Let $f \subset S$ be an embedded fibre of level 1. Then $\pi$ is
the projection  centered  at the point $P_0 \in f_0$, where $P_0 \in Sing(X_0)$. Clearly,
$f =f_0 + E$, where
$E$ is the exceptional divisor and $f_0$  still denotes   the
other component of
$f$.  
Setting $A:= f_0 \cdot E$, $P_i \in f_0$ and $Q_i \in E$ (where $P_i \ne A \ne Q_i$ and $P_i \ne Q_i$, for $i=1,2$), the
possible cases are the following:
$$
\eqalign{
(a) & \qquad \Phi = P_1 + P_2 + Q_1 + Q_2 \cr
(b) & \qquad \Phi = P_1 + P_2 + 2Q_1 \cr
(c) & \qquad \Phi = 2P_1  + Q_1 + Q_2 \cr
(d) & \qquad \Phi = 2P_1 + 2Q_1  \cr
(e) & \qquad \Phi = P_1 + 2A + Q_1 \cr
(f) & \qquad \Phi = P_1 + 3A \quad \hbox{(where $X_K \cdot f_0 =P_1 +A$ and $X_K \cdot E
= 2A$)}\cr (g) & \qquad \Phi = 3A + Q_1 \quad \hbox{(where $X_K \cdot f_0 =2A$ and $X_K
\cdot E = A + Q_1$)}.\cr }
$$
The  picture below illustrates the corresponding singularities of $X_0$.
\bigskip
 \centerline{
 \epsfxsize=14cm 
 \epsfbox{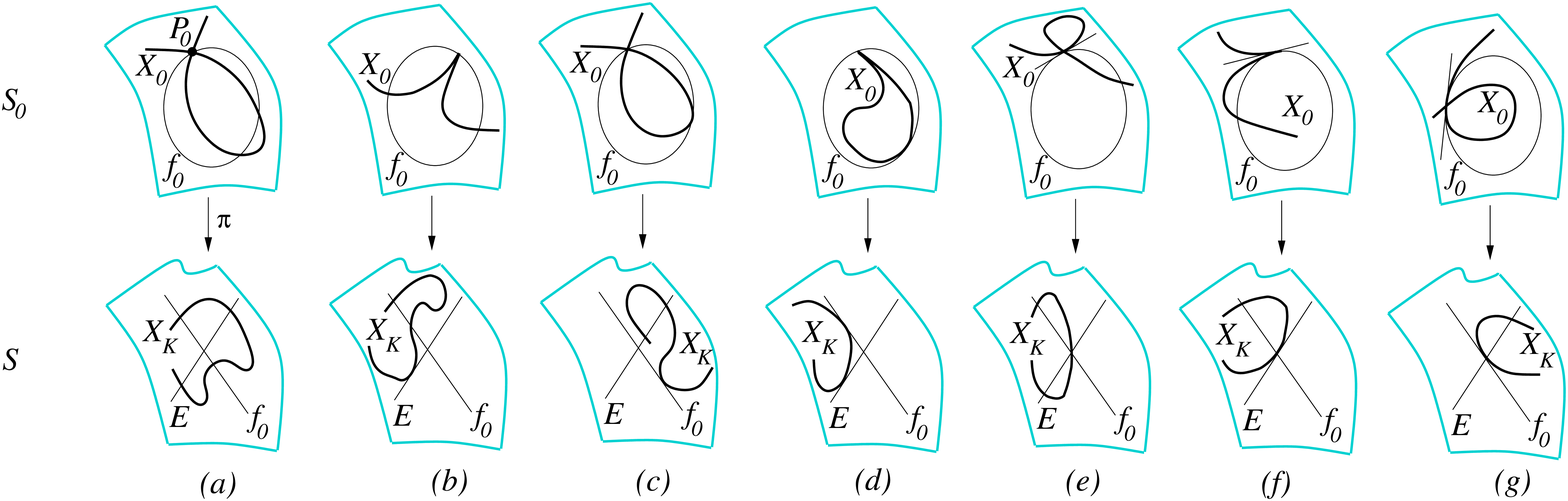}
 }
 \medskip
 \centerline {Figure 1}
  \medskip
\goodbreak

It is clear that, in all the cases above, $X_0$ has a double point: more precisely, either
a node, in cases $(a), (c), (e), (g)$, or an ordinary cusp, in cases $(b),(d),(f)$.

\medskip
\noindent
A description of the  double points of an algebraic curve  can be
found, for instance, in [{\bf 10}], Lect. 20.  \hb
Here let us just recall that a {\it node of  $n$-th kind} is a double point analitically equivalent to $y^2-x^{2n} =0$. In particular, if $n=1,2,3$, it is called {\it (ordinary) node}, {\it  tacnode}, {\it oscnode}, respectively. \hb
Moreover, a {\it cusp of  $n$-th  kind} is a double point analitically equivalent to $y^2-x^{2n+1} =0$. In particular, if $n=1,2$, it is called {\it (ordinary) cusp} or {\it  ramphoid cusp}, respectively.

\medskip
\noindent
{\bf Definition.} We say  for short  that a double point $P_0$ of $X_0$ is {\it
transversal} if the tangent line to the fibre $f_0$ at $P_0$ does not coincide with
any of the tangent lines to
$X_0$ at $P_0$; it is {\it tangent} otherwise.

\medskip
\noindent
{\bf Example 3.3.} Assume that $S$ is a surface ruled by conics having a fibre $f$ of
type $(2A)$, as defined in 0.2. Clearly (see [{\bf 6}],  Sect. 3) this fibre arises from a fibre
$f_0 \subset S_0$ by projecting it from two points. More precisely, the projection 
$\pi: S_0 \longrightarrow S$ can be factorized by  $\pi = \pi_{P_1}
\circ \pi_{P_0}$, where $P_0 \in f_0$ and $P_1 \in f_1:= f_0 + E \subset
\pi_{P_0}(S_0)$ and $P_1 \ne f_0 \cdot E$. 
There are two possibilities: either $P_1 \in f_0$ or $P_1 \in E$. \hb
In the first case, $f= E+E^{(1)}$, while in the second one, where $P_1$ is infinitely
near  to $P_0$, we have 
$f = f_0 +E^{(1)}$ (in both cases $E^{(1)}$ denotes the exceptional divisor of the
blowing--up centered at $P_1$). Moreover, in both configurations, $f$ turns out to be a
union of two lines meeting in an ordinary double point for the surface $S$. \hb 
Let us start by scketching the  situations
corresponding to the configuration $(a)$ (in both cases $f=E+E^{(1)}$ and  $f = f_0
+E^{(1)}$) and the configurations $(b)$ and $(d)$ (both in the case  $f = f_0
+E^{(1)}$).

\bigskip
 \centerline{
 \epsfxsize=9cm 
 \epsfbox{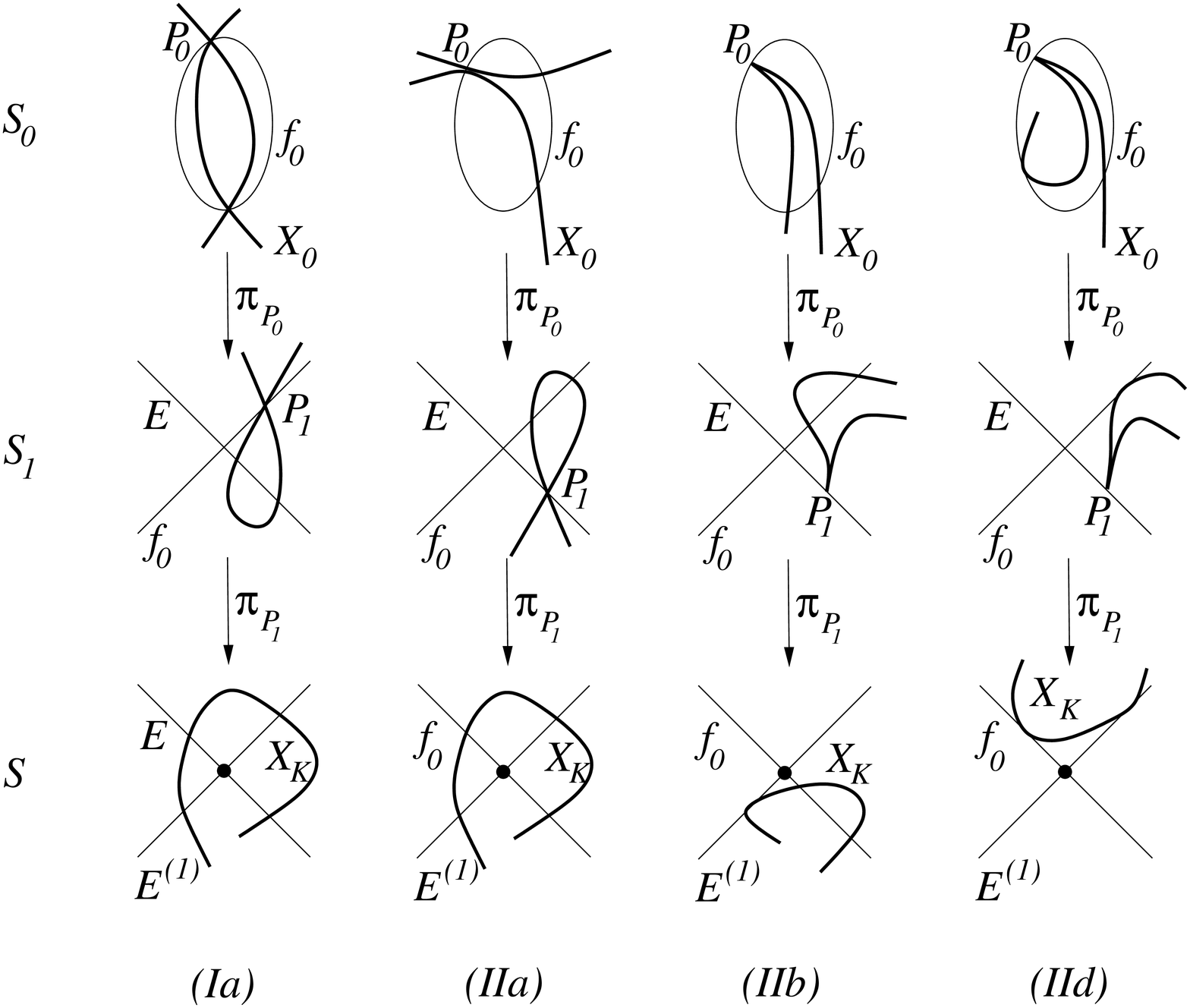}
 }
 \medskip
 \centerline {Figure 2}
  \medskip

  The construction $(Ia)$  gives  $X_0$ to have two nodes  on
the fibre $f_0$; in $(IIa)$ the curve $X_0$ has a tacnode, while in  $(IIb)$
and
$(IId)$ it has a ramphoid cusp.
Finally, one can easily see that the cases related to  $(e), (f),
(g)$ do not occur.

\medskip
\noindent
{\bf Remark 3.4.} The two examples above lead us  to a general pattern.
If $X_0$ has only one singular point  $P_0 \in f_0$ and $f$ is of type $(nA)$, then:
\item{-}  $f = f_0 + E^{(n-1)}$ and $\pi$ can be factorized by $\pi = \pi_{P_{n-1}} \circ
\cdots \circ \pi_{P_1} \circ \pi_{P_0}$, where $P_{i+1} \in E^{(i)}$ for all $i$;
\item {-} the type of the singularity of $P_0$ depends only on the intersection 
 $X_K \cdot E^{(n-1)}$ on $S$, so
 we can always assume that the two points
given by $X_K \cdot f_0$ on $S$ are distinct.

\smallskip
\noindent
 We can now  complete  3.3:  if $X_0$ has one singular point on $f_0$, then the  significant cases are $(IIa)$
and $(IIb)$, where   $X_0$ has  a
transversal tacnode or a transversal ramphoid cusp. Note that the
difference between these two cases is that $X_K$ is tangent (resp. transversal) to
$E^{(1)}$ on $S$.

\medskip
\noindent
{\bf Example 3.5.}
In the same way, we get the possible singularities in the
case $F_3(A)$,  as  this picture shows:

 \bigskip
 \centerline{
 \epsfxsize=6.8cm 
 \epsfbox{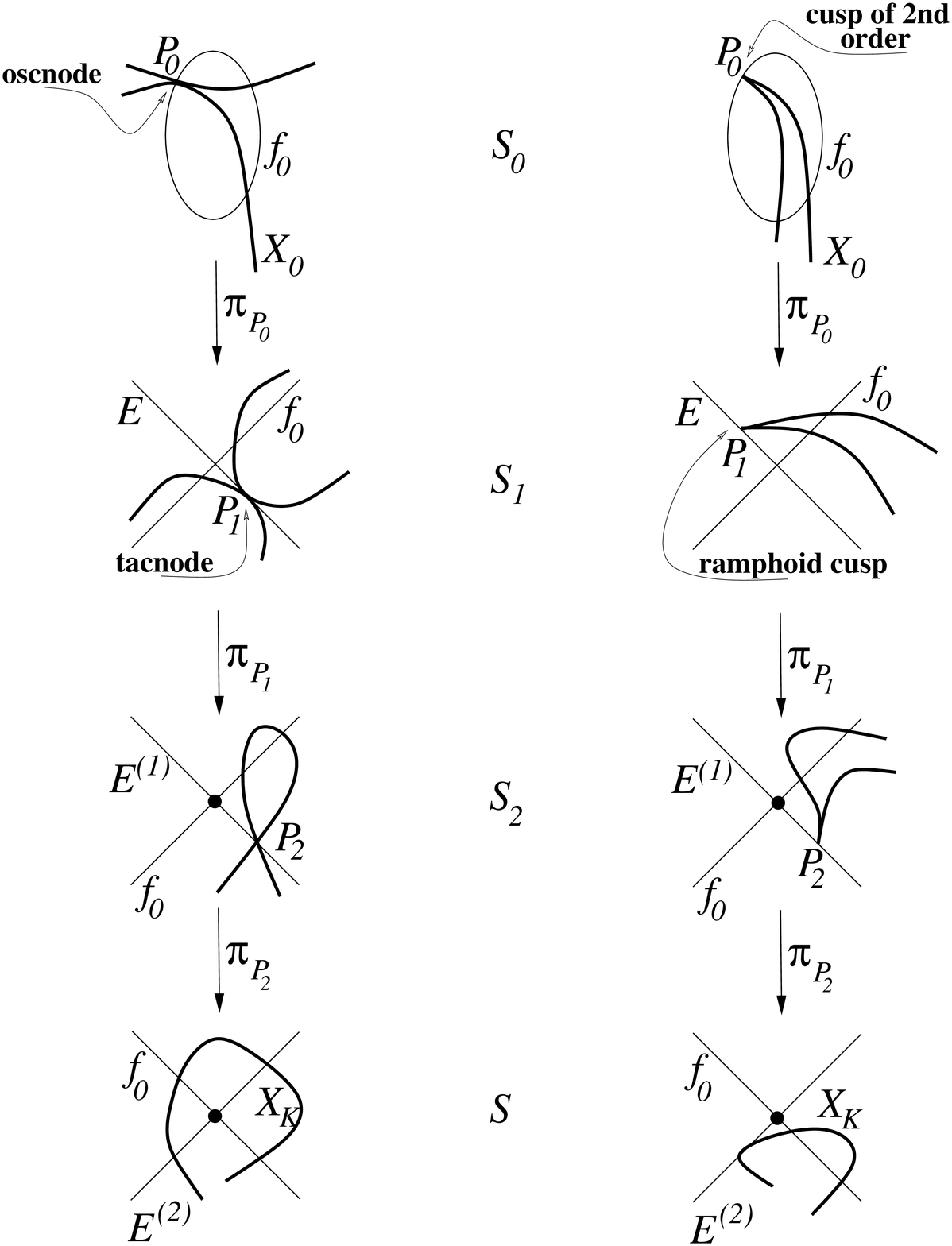}
 }
 \smallskip
 \centerline {Figure 3}
  \smallskip

The above study can be easily generalized, obtaining the following result:

\smallskip
\proclaim
Proposition 3.6. The possible singularities of $X_0 \subset S_0$ 
arising from a fibre of $S$ of type $F_n(A)$, where $n \ge 2$, are the following points
on the same fibre $f_0 \subset S_0$:
\item {$(\bullet)$} if $n=2$  there is either one double point  of second kind (either a transversal tacnode
or a transversal ramphoid cusp) or two double points of first kind (either node or cusp);
\item {$(\bullet)$} if $n \ge 3$  there is either one double point of   $n-$th  kind (transv.) or two double
points of lower kind.
 \cvd

\smallskip
\noindent
Note that in the case of two double points on $f_0$,    
these  two points are of kind $h$ and $k$, where $h+k = n$.

\smallskip
\noindent
{\bf Example 3.7.} Assume now that $S$ is a surface ruled by conics having 
 a fibre $f$ of type $(2D)$. Clearly (see [{\bf 6}], Sect. 3) this fibre arises from a fibre
$f_0 \subset S_0$ by projecting it from two infinitely near points. More precisely,
if $\pi: S_0 \longrightarrow S$ is the considered projection, then 
$\pi = \pi_{P_1} \circ \pi_{P_0}$, where $P_0 \in f_0$ and, if  $f_1:= f_0 + E
\subset \pi_{P_0}(S_0)$, then  $P_1 := f_0 \cdot E$. As noted in [{\bf 6}], the
fibre of $S$ corresponding to $f_0$ is given by
$f=2E^{(2)}$:  it is a totally degenerate conic containing two singular
points of $S$, which correspond to the lines $f_0$ and $E$. 
Since $f$
consists of a double line, the four--gonal divisor can be either  $2A +2B$
(where
$A, B \in E^{(2)}$ are distinct  points  non singular for $S$) or 
$4A$, as the following picture describes:

\smallskip
 \centerline{
\epsfxsize=4.8cm 
 \epsfbox{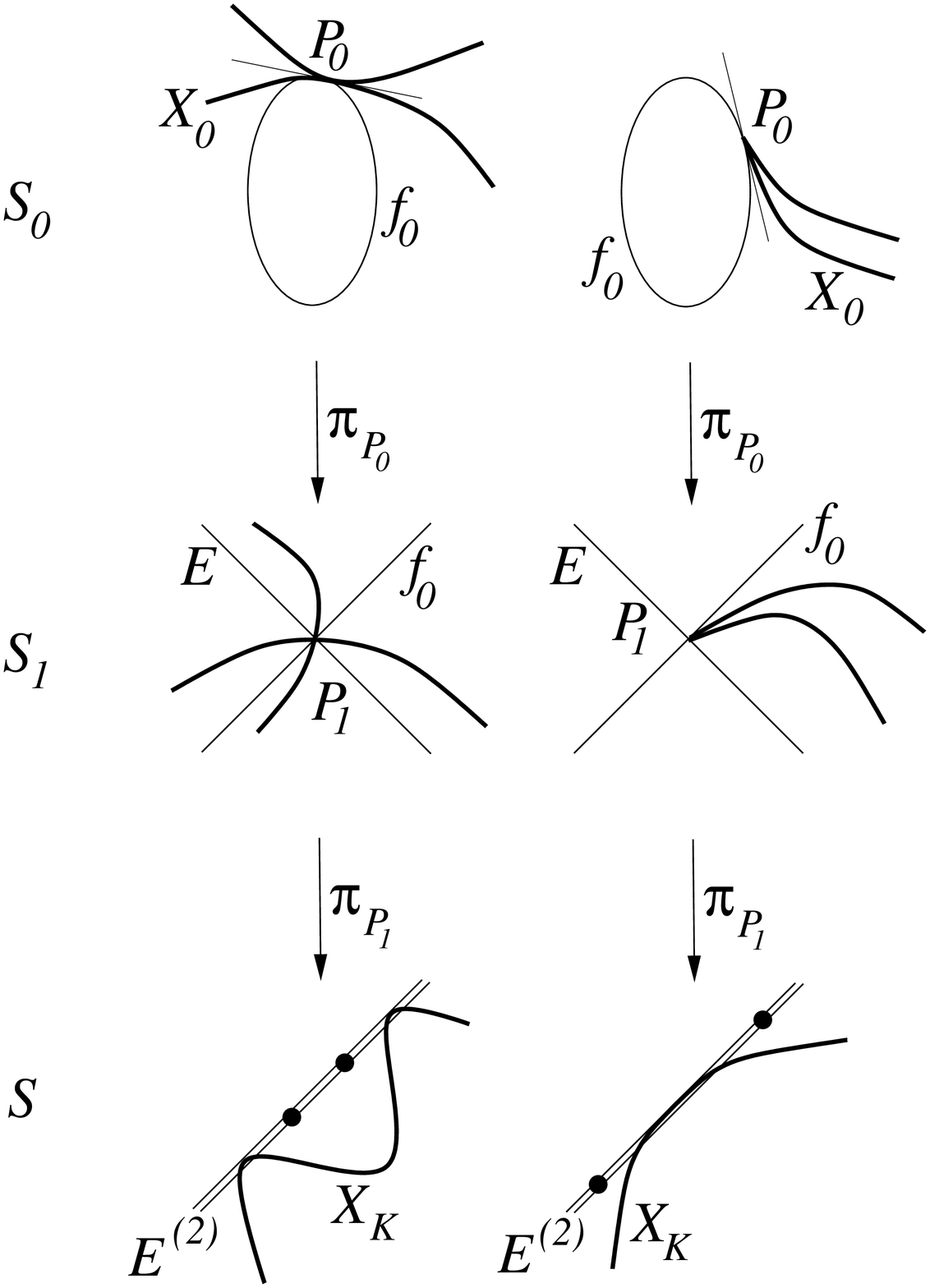}
 }
 \smallskip
 \centerline {Figure 4}
 \medskip

\goodbreak
  
It is clear that the first configuration leads to a tangential tacnode and the second one gives a tangential ramphoid
cusp of first order. 
With the same argument as before, we easily get
the following result:

\medskip
\proclaim
Proposition 3.8. The possible singularities of $X_0 \subset S_0$ arising from a
fibre of $S$ of type $F_n(D)$, where $n\ge 2$, consist of a unique singular point of the
corresponding fibre $f_0 \subset S_0$ as follows:
\item {$(\bullet)$} if $n=2$ then  there is either a tangential tacnode or a tangential ramphoid cusp;
\item {$(\bullet)$} if $n \ge 3$ then  there is a tangential double point of   $n$-th kind.

\medskip
\noindent
Collecting 3.2, 3.6, 3.8, we obtain the following complete description of the
possible singularities of $X_0$.

\medskip
\proclaim
Theorem 3.9. Let $S$ be a surface ruled by conics containing $X_K$ and let $X_0
\subset S_0$ be  birational models of $X_K$ and $S$ respectively, where $S_0$ is a g.r.c.
surface. Let
$\pi: S_0
\longrightarrow S$ be the usual projection. Assume that $f$ is the unique singular fibre
of $S$ and set $f_0$ the corresponding fibre of $S_0$. \hb
Then the singular points of $X_0$ belong to $f_0$ and are, as far as  $f$ is of type $F_1$, of one of the following types,
   $F_n(A)$,  $F_n(D)$, for $n \ge 2$: 
\item \item {$F_1$ -} one singular point: either a node or a cusp, both of them either
tangential or transversal;
\item \item {$F_n(A)$ -}  only transversal singular points and precisely:
\item\item {$\qquad (a)$} one double point of   $n$-th kind;
\item\item {$\qquad (b)$} two double points of orders $h,k < n$, where $h+k=n$;
\item \item {$F_n(D)$ -}  only one tangential double point of   $n$-th kind;
\smallskip
In particular, all the singular points of $X_0$ are double points.
\cvd
\par

\bigskip
\noindent
{\bf 4.  ``Standard'' birational models of $X_K \subset S$} 

\medskip
\noindent
In Section 2 we studied the set $\grc_L(S)$ consisting of the g.r.c. surfaces $S_0$ such
that $S$ can be obtained from $S_0$ by a sequence of $L$ monoidal transformations (here
$L$ is the level of $S$). In this section we are going to determine one of such surfaces
in a sort of ``canonical'' way: this will be called ``standard'' birational model of $S$.

\medskip
\proclaim
Proposition  4.1. Let $X_0 \subset S_0 \in \grc_L(S)$ be as usual. Then 
$$
\grc_L(S) = \{elm_\Sigma (S_0) \; | \; \Sigma \subseteq Sing(X_0)\}
$$
i.e. each $S'_0 \in \grc_L(S)$ can be obtained from $S_0$ by a sequence of elementary
transformations centered in singular points of $X_0$ (or infinitely  near  to them) and
conversely.
\par
\pf Consider a surface $S'_0 \in \grc_L(S)$ and the corresponding model of $X_K$, say
$X'_0 \subset S'_0$. As in 2.2, denote by $\pi$ and $\pi'$ the projections centered in
the singular points (possibly infinitely near) of $X_0$ and $X'_0$, respectively. We get
then the diagram
$$
\matrix{
S_0 & --- \rightarrow & S'_0 \cr
&{}_\pi\searrow \qquad  \swarrow_{\pi'} &\cr
&S 
\cr}
$$
where the horizontal arrow denotes a suitable sequence of elementary transformations
centered in (some of) the singular points of $X_0$. \hb
Conversely, note that each elementary transformation of $S_0$ can be obtained by considering
an embedded model of
$S_0$ which is ruled by lines and projecting it from a finite number of points. In this
way, we get a birational model $S'_0$ of $S$ which is a geometrically ruled surface. If
$X'_0 \subset S'_0$ is the corresponding curve, it is clear that $\delta(X'_0) =
\delta(X_0)$ if and only if the above projection is centered in singular points of $X_0$
(this is due to the fact that the singular points of $X_0$ are double points for
3.9). Therefore, if $S'_0 = elm_\Sigma(S_0)$, where $\Sigma \subseteq Sing(X_0)$,
 using 2.4, the level of $S'_0$ coincides with $\delta(X'_0) = \delta(X_0) = L$,
hence $S'_0 \in \grc_L(S)$, as requested.
\cvd

\medskip
Among the surfaces $S_0$ geometrically ruled by conics belonging to $\grc_L(S)$ (and the
corresponding curves $X_0$), we are going to establish a way for choosing one particular
model of $S$ (and hence of $X_K$). In order to do this, we give the following notion.

\medskip
\noindent
{\bf Definition.} Given a surface $S$ ruled by conics, we say that a surface $\ss_0 \in \grc_L(S)$ is a {\it
standard model} of $S$ if its invariant is
$$
t := \min \{\tau_0=t(S_0) \; | \; S_0 \in \grc_L(S) \}.
$$
Let us consider now the curve $X_K \subset S$ and the corresponding birational model,
say $\xx_0 := \rho(X_K) \subset \ss_0$, where $\ss_0$ is a standard model of $S$. We say also that $\xx_0$ is a
{\it standard model} of $X_K$. \hb
Finally, if $\ss_0$ is a standard model of $S$,  we denote  the corresponding
invariant $\lambda_0$ by $\lambda$.

\medskip
\proclaim
Theorem 4.2.
Let $S$ be as before, $L$ be its level,  $S_0 \in \grc_L(S)$ be a birational model of $S$ 
 of invariant $\tau_0$ and
$X_0$ be the  model of $X_K$ on $S_0$.  If we assume that $t>0$, then the
following facts hold:
\item {$i)$} if $S_0$ is a standard model, then the singular points of $X_0$ belong
to the minimum unisecant $C_0$ of
$S_0$;
\item {$ii)$}   there is exactly one standard model $\ss_0$ of $S$;
\item {$iii)$}  if the singular points of $X_0$ belong
to the minimum unisecant $C_0$ of
$S_0$, then $S_0 = \ss_0$.
\par
\pf 
Consider first the model $X' \subset R_{1,\tau_0 +1} \cong S_0$. 
We know that  $X' \sim 4 C_0
+ (\lambda_0 +\tau_0)f$
 and $\delta(X') = 3(\lambda_0 -\tau_0 -1) -g$  by 2.2. In particular, the level of $S$
is $L =  3(\lambda_0 -\tau_0 -1) -g$. \hb
Consider a singular point $T$ of $X'$ and the projection $\pi_T$ from $T$. {}From
4.1, $\pi_T(R_{1,\tau_0 +1})$ belongs to $\grc_L(S)$. \hb 
$(i)$ If $S_0$ is a standard model, then $\tau_0 = t$. Assume that the  point  $T$
does not belong to $C_0$. Then the invariant of  
$\pi_T(R_{1,t +1})$ is $t-1$,  while $t$ is the minimum invariant of the surfaces
belonging to $\grc_L(S)$. \hb
$(ii)$ Let $\ss_0 \cong R_{1,t+1}$ be a standard model and let $S'_0$ be another surface
in
$\grc_L(S)$. {}From 4.1, we know that $S'_0 = elm_\Sigma(\ss_0)$, where $\Sigma
\subseteq Sing(\xx_0)$. For simplicity, assume that $\Sigma = \{T\}$, where
 $T$ is a singular point of $\xx_0$. {}From $(i)$, we have that $T \in C_0$ and, from
3.9, we know that $T$ is a double point of
$\xx_0$, so $T = A_1+A_2$, where $\Phi:= A_1+A_2+A_3+A_4$ is the four--gonal divisor on
the fibre $\ff_0$ containing $T$.  \hb
Clearly,  $S'_0 = \pi_T(R_{1,t +1})$, so  the curve $X'_0$ has a double point on the
fibre $\ff'_0$ given by
$A_3+A_4$ and  such point does not belong to the unisecant curve $C'_0$ of $S'_0$.
Therefore we get from $(i)$ that $S'_0$ is not a standard model of $S$.
\hb
$(iii)$  An analogous argument.
 \cvd

\medskip
\proclaim
Proposition 4.3. With the above notation, if  $t>0$ then the singular points of
$\xx_0$ belong to distinct fibres.
\par
\pf Also in this case consider the model $X' \subset R_{1,t +1} \cong \ss_0$ and
assume that there exists a fibre containing two distinct singular points of $X'$,   
$P_1$ and
$P_2$, say . Clearly, one of them,  $P_1$ say, does not belong to $C_0$.
So, by projecting $R_{1,t +1} $ from $P_1$ we get a contraddiction with the 
argument used in 4.2.
\cvd 

\medskip
\proclaim
Theorem 4.4. With the notation above, the surface $S$ has degree
$$
\deg(S) = 4(\lambda-t-2)-\delta_t  =g+ \lambda - t -5.
$$
\par
\pf
Since $\ss_0 = \varphi_{2C_0+ (\lambda -2)f} (\hs{t})$ and $C_0^2 = -t$,
then 
$$
\deg(\ss_0 \,)= (2C_0+ (\lambda -2)f)^2=4(\lambda-t-2).
$$
Moreover, from 2.4   we have that $\deg(S) = \deg(\ss_0 \,)- \delta_t $, so
the first equality holds. The second equality   follows  immediately from
$\delta_{t}=3(\lambda -t-1) -g$ (see 2.2,
$(iii)$).
\cvd

\bigskip
\noindent
{\bf 5. Bounds on the invariants $\lambda$ and $t$} 
\medskip
Let us come back to the global description of the four--gonal curve $X$ of genus $g$ whose canonical model is 
$X_K \subset S \subset V = \Bbb P({\cal O}(a) \oplus {\cal O}(b) \oplus {\cal O}(c))
\subset \Bbb P^{g-1}$ and the surface $S$ is (as in 1.1) the surface of minimum degree. 
\hb  
We have chosen $\xx_0 \subset \ss_0 \cong \hs t$ as a   pair  of
standard models of  $X_K \subset S$ respectively.
Since  the model $X_t \subset \Bbb F_t$ is again a
four--secant curve, it is of the type $X_t
\sim 4C_0 + (\lambda + t) f$.

\medskip
So far we have defined a set of integers, $a,b,c,t, \delta, \lambda$ (here, for
simplicity,  $\delta :=
\delta_t$), that are   {\it invariants} of the curve $X$. All of them will be useful to
describe its geometry.

\medskip
Let us start  with  the dependence of the first three invariants $a,b,c$ on the others $t,
\delta, \lambda$.

\medskip
\noindent
{\bf Remark 5.1.}
Consider the isomorphism
$$
\varphi_{2C_0+(\lambda-2)f}: \;  \Bbb F_t \longrightarrow
\ss_0 \subset  \Bbb P^{g-1+\delta}
$$
and the volume $V_{\ss_0} \subset \Bbb P^{g-1+\delta}$ generated by
$\ss_0$. {}From 1.8, [{\bf 4}], we have that 
$$
V_{\ss_0} = \Bbb P({\cal O}(\lambda-2-2t) \oplus
{\cal O}(\lambda-2-t) \oplus
{\cal O}(\lambda-2) ).
$$
 If we consider the projection  
$\pi:\; \Bbb P^{g-1+\delta}
\rightarrow \Bbb P^{g-1}$ centered at the singular locus of $\xx_0$, it is clear that 
$\pi(V_{\ss_0}) =V_S$.  \hb
Using 4.2 $(i)$, if $t >0$ then the
 singular  points of $\xx_0$ are contained in the unisecant of minimum degree of $\ss_0$
and hence of
$V_{\ss_0}$. Moreover, if these points are  all distinct, then $V_S$ has the form:
$$
V_S = \Bbb P({\cal O}(\lambda-2-2t-\delta) \oplus
{\cal O}(\lambda-2-t) \oplus
{\cal O}(\lambda-2) ).
$$
On the other hand, taking into account that $c = g-3 - a-b$, the scroll above is:
$$
V_S = \Bbb P({\cal O}(a) \oplus {\cal O}(b) \oplus {\cal O}(g-3-a-b) ).
$$
Hence, comparing  the two expressions of $V_S$ and using the equality 
$\delta = 3(\lambda- t -1) -g$ (see 2.2 $(iii)$), we  obtain: 
$$
a=g+t-2\lambda+1 
\quad \hbox{and} \quad
b=\lambda-t-2.
$$
Note that, if $t >0$ but the $\delta$ double points of $\xx_0$ are not all distinct, then
$a \ge g+t-2\lambda+1$.

\medskip
\proclaim
Proposition 5.2.   With the above notation, if 
$V=\Bbb P({\cal O}(a) \oplus {\cal O}(b) \oplus {\cal O}(c))$, then
$$
a+b \ge {{g-5} \over 2} .
$$
\par
\pf 
Let us consider the curve $X_K \subset V$ and the
ruled surface $R_{a,b} = \Bbb P({\cal O}(a) \oplus {\cal O}(b)) \subset V$.  In order to apply the Intersection Formula 
$(IF)$ in Section 0, we observe first that $R_{a,b}$ and $X_K$ meet properly on $V$, i.e. 
$$
\dim(R_{a,b} \cap X_K) = \dim(R_{a,b}) + \dim(X_K) -\dim(V) =0.
$$
To see this note that  $X_K$ cannot be contained in $R_{a,b}$, otherwise the general $4$--gonal divisor on $X_K$ would span a line instead of a plane, against the Geometric Riemann--Roch Theorem. \hb
Hence $\dim(R_{a,b} \cap X_K) = 0$ and we can apply  $(IF)$, which gives the (non--negative) degree of the intersection:
$$
0 \le \deg_V(R_{a,b} \cdot X_K) = 4(a+b)+2g-2 - 4(g-3) = 2(a+b) -g +5
$$
and this proves the requested inequality.
\cvd

\medskip
The lower bound of $\lambda$ in terms of $t$ given in the previous section
can be improved. Namely, we saw that $\lambda \ge \max\{3 t, t+5\}$ (see 2.2).

\medskip
\noindent
{\bf Remark 5.3.} Assume that $t \ge 1$ and the $\delta$ singular points of $X_t$  are distinct. 
Clearly
$$
2 \delta \le \int C_0 \cdot X_t = \int C_0 \cdot (4C_0 +(\lambda +t)f) = \lambda - 3 t
$$
hence
$$ 
\lambda \ge 2 \delta + 3 t.
$$
Since $\delta  = 3(\lambda-t-1)-g$ (see 2.2 $(iii)$), we easily obtain: 
$$
\lambda \le {{2g+3t+6} \over {5}}
.\eqno(7)
$$

\goodbreak

\medskip
\proclaim
Proposition 5.4.   The following properties hold :
\item {$(i)$} for any $t$:
$$
\lambda \ge {g \over 3}+t+1 ;
$$
\item {$(ii)$} if $t = 0$ then
$$
\lambda \le  {{g+3} \over 2} ;
$$
\item {$(iii)$} if $t \ge 1$ then
$$
\lambda \le t+ {{g+3} \over 2} 
\quad \hbox{and} \quad
t \le {{g+3} \over 4} ;
$$
\item {$(iv)$} if $t \ge 1$ and the double points of $X$ are all distinct, then
$$
\lambda \le 
{{g+3} \over 2} 
\quad \hbox{and} \quad
t \le {{g+3} \over 6}.
$$
\par
\pf  $(i)$ It comes from 2.2 $(i)$, since 
$p_a(\xx_0)=3(\lambda-t-1) \ge g$. \hb
$(ii) - (iii)$
Using 1.1 and 4.4 we have
$$
g+ \lambda -t-5=\deg(S) \le 
\left\lceil{{3g-8} \over 2}\right\rceil
\quad \Rightarrow \quad
\lambda -t \le \left\lceil{{3g-8} \over 2}\right\rceil-g+5=
\left\lceil{{g+2} \over 2}\right\rceil
$$
hence,  we obtain the required bounds either if $t=0$  or  if $t \ge 1$. Moreover, from
2.2 we have $\lambda \ge 3t$; so, using the previous bound of $\lambda$ in $(iii)$, we
finally get $t \le \lambda/3 \le t/3 +  {{g+3} \over 6}$ and this concludes the proof.\hb
$(iv)$
In this case, we can apply 5.3. Using   $3(\lambda -t-1)-g=\delta \ge 0$ followed by  $(7)$, we
get:
$$
t\le \lambda - {{g+3} \over {3}} \le
 {{2g+3t+6} \over {5}}- {{g+3} \over {3}}  
\quad \Rightarrow \quad t \le {{g+3} \over 6}.
$$
Using this bound and $(7)$ we finally get
$\displaystyle {\lambda \le {{g+3} \over 2}}$.
\cvd

\bigskip
\noindent
{\bf 6. Geometric meaning of the invariant $\lambda$}

\medskip
Let us keep the notation of the previous section: $S$ is a surface  ruled by conics such that $X_K \subset S \subset V$ and $L$ denotes its level. Take  a standard model $\ss_0 \in \grc_L(S)$ and  consider its embedded 
model  $R_{1,t+1} \subset \Bbb P^{t+3}$. \hb  
 Let us denote as usual by $X'
\subset R_{1,t+1}$ the corresponding model of $X_K$, where
$X' \sim 4 C_0 + (\lambda +t)f$.

\medskip
\noindent
{\bf Remark 6.1.} Note that such $X'$ has only double points as singularities (see
3.9).

\medskip
\noindent
{\bf Remark 6.2.} Denote by $H_{X'}$ the hyperplane section of $X' \subset R:= R_{1,t+1}
\subset \Bbb P^{t+3}$. \hb
Since $H_R \sim C_0+(t+1)f$ then 
$$
H_{X'} = H_R \cdot X' \sim \Phi + \Delta, \quad \hbox{where} \quad
\Phi \in g^1_4  \quad \hbox{and} \quad \Delta \in g^{1+t}_{\lambda +t}.
$$
In particular
$$
\deg(H_{X'}) =  \lambda +t+4
$$
and one can easily verify that $X'$ is the embedding of minimum degree of the curve $X_K$.

\medskip
\noindent
{\bf Definition.} A linear system $|D|$ on a curve $X$ is called {\it primitive} if, for
each point $P \in X$, the linear system $|D +P|$ has $P$ as base point. Equivalently,
$\dim |D+P| = \dim |D|$.

\medskip
It is not difficult to see that the following property of $X' \subset \Bbb P^{t+3}$, here
stated for a standard model
$\ss_0$, holds also for any birational model $S_0 \in \grc_L(S)$.

\medskip
\proclaim
Proposition 6.3. Let  $\ss_0 \cong R_{1,t+1} \subset \Bbb P^{t+3}$ be a standard model of $S$. Let $\Phi$ and $\Delta$ be as before and $X' = X_{\Phi + \Delta} \subset R_{1,t+1}$ be as usual.
If $g > 13$ then the following facts hold:
\item {(i)} the divisor $\Phi + \Delta$ is a special divisor on $X$; in particular $K-\Phi-\Delta$ is an effective divisor.
\item {(ii)} The curve  $X'\subset \Bbb P^{t+3}$ is linearly normal.
\par

\goodbreak

\pf 
$(i)$ 
It is enough to show that $h^0({\cal O}(K-\Phi-\Delta)) >0$ or, equivalently by Riemann--Roch Theorem, that $\lambda < g-1$. If $t=0$, it follows immediately from 5.4 $(ii)$. \hb
If $t \ge 1$, still from 5.4 $(iii)$, we have:
$$
\lambda \le t + {{g+3} \over2}
\quad \hbox{and} \quad 
t \le {{g+3} \over4}
\quad \Rightarrow \quad 
\lambda \le {{3g+9} \over4} < g-1
$$
where the last inequality is true since $g > 13$ by assumption. Finally, observe that $\Phi + \Delta$ special implies that $K-\Phi-\Delta$ is an effective divisor. \hb
$(ii)$
Let us recall that (as in 5.1) the surface  $\ss_0$ is naturally embedded, via the
isomorphism $\varphi_{2C_0+(\lambda-2)f}$, in a projective space: namely  $\ss_0 \subset
V_{\ss_0}
\subset
\Bbb P^{g-1+\delta}$, where 
$$
V_{\ss_0} = \Bbb P({\cal O}(\lambda-2-2t) \oplus
{\cal O}(\lambda-2-t) \oplus
{\cal O}(\lambda-2) )
$$
and $t \ge 0$. If $t>0$, denoting by $M:= \langle \varphi_{2C_0+(\lambda-2)f}((\lambda
-3-t) \Phi) \rangle$,  it is clear that 
$$
\pi_M: V_{\ss_0} \longrightarrow 
 \Bbb P({\cal O}(1) \oplus
{\cal O}(t+1)) = R_{1,t+1}.
$$
This map can be factorized as follows: setting $\Sigma$ the divisor of the singular
points of $\xx_0$ and taking into account that $K-\Phi-\Delta$ is an effective divisor on $X$ from $(i)$, put:
$$
L:= \langle \varphi_{2C_0+(\lambda-2)f}(\Sigma) \rangle, \quad
N:= \langle \varphi_{K}(K-\Phi-\Delta) \rangle .
$$
Then we have the following diagram:
$$
\matrix{
&&\xx_0 & \subset & \ss_0 & \subset & V_{\ss_0} & \subset &  \Bbb P^{g-1+\delta} \cr
&{}^{\overline \varphi}\nearrow&\downarrow&&\downarrow&&\downarrow &&\mapdown
{\pi_L} \cr
\hs t \supset X_t & \hfill \mapright \varphi_K,30 & X_K & \subset & S & \subset & V &
\subset & 
\Bbb P^{g-1}
\cr
&{}_{\varphi'}\searrow&\downarrow&&&\searrow&\downarrow &&\mapdown {\pi_N} \cr
&&X' &  & \subset &  & R_{1,t+1} & \subset &  \Bbb P^{t+3} \cr
}
\eqno{(8)}
$$
where $\overline \varphi := \varphi_{2C_0+(\lambda-2)f}$, $\varphi' =
\varphi_{\Phi + \Delta}$ and
$$
\pi_N \circ \pi_L = \pi_M. 
$$
Note that $\xx_0$ is not linearly normal. Namely, $\xx_0$ is not special; if it was
linearly normal, then $\dim \langle \Phi \rangle = 3$ in $\Bbb P^{g-1+\delta}$, while
$\xx_0$ is contained in the scroll $V_{\ss_0}$  which is ruled by planes. \hb 
Hence we
have to consider its  normalization $\wxx \subset  \Bbb P^{g-1+2\delta}$, and the
corresponding scroll 
$$
W:= \bigcup_{\Phi \in g^1_4} \; \langle \Phi \rangle \subset  \Bbb P^{g-1+2\delta}.
$$
It is easy to see that $W$ is ruled by planes. Setting $\wll:= \langle \Sigma  \rangle
\subset  \Bbb P^{g-1+2\delta}$, the projection $\pi_{\wll}$ factorizes through the
normalization map, say $\Pi$, as follows:
$$
\matrix{
\wxx &  \subset & W & \subset &  \Bbb P^{g-1+2\delta} \cr
\downarrow&&\downarrow &&\mapdown {\Pi} \cr
\xx_0 &  \subset & V_{\ss_0} &\subset & \Bbb P^{g-1+\delta}\cr
\downarrow&&\downarrow &&\mapdown {\pi_L} \cr
X_K  & \subset  & V & \subset &  \Bbb P^{g-1} \cr
}
\eqno{(9)}
$$
and
$$
\pi_L \circ \Pi = \pi_{\wll}.
$$
Setting
$$
\wmm := \langle (\lambda -3-t) \Phi \rangle \subset  \Bbb P^{g-1+2\delta}
$$
and keeping into account $(8)$ and $(9)$ we finally obtain:
$$
\matrix{
\wxx &  \subset & W & \subset &  \Bbb P^{g-1+2\delta} \cr
\downarrow&&\downarrow &&\mapdown {\pi_{\wll}} \cr
X_K &  \subset & V &\subset & \Bbb P^{g-1}\cr
\downarrow&&\downarrow &&\mapdown {\pi_N} \cr
X'  & \subset  & R_{1,t+1} & \subset &  \Bbb P^{t+3} \cr
}
$$
where
$$
\pi_N \circ \pi_{\wll} = \pi_{\wmm}. 
$$
Since $ \pi_{\wmm}: \; \wxx \longrightarrow X'$ and $\wxx$ is linearly normal, than also
$X'$ is linearly normal. \hb
If $t=0$,  the  proof runs in a similar way.
\cvd

 \bigskip
\proclaim
Proposition 6.4. Let  $\ss_0 \cong R_{1,t+1} \subset \Bbb P^{t+3}$, $\Phi$, $\Delta$  and 
$X' = X_{\Phi + \Delta}$ be as usual.
If $g > 13$ then the following facts hold:
\item {$i)$} The linear system $|\Delta|$ defined before is
primitive;
\item {$ii)$} if $B \subset \Delta$ is a divisor on $X'$ such that $B \in g^1_\beta \ne
g^1_4$, then $B \sim \Delta - A_1 - \cdots -A_t$,  for suitable $A_i \in X' \setminus C_0$ for
all $i$. In particular, $\beta = \lambda$.
\par
\pf $i)$ Assume that there exists $P \in X'$ such that $\Delta +P \in
g^{2+t}_{\lambda +t+1}$ and  consider the model of $X_K$ given by
$X_{\Delta +P} \subset \Bbb P^{t+2}$. Keeping into account 6.3, we have that 
$X' = X_{\Phi + \Delta}$ is linearly normal in $\Bbb P^{t+3}$. Hence we can consider the
following diagram:
$$
\matrix{
&  & X_{\Phi + \Delta} & \subset R_{1,t+1} \subset &  \Bbb P^{t+3} \cr
&\nearrow &\mapdown {}&& \mapdown {\pi_{\langle \Phi - P \rangle}}\cr
\cr
X & \longrightarrow & X_{\Delta +P}  & \subset  &  \Bbb P^{t+2} \cr
&\searrow &\mapdown {}&& \mapdown {\pi_{P}}\cr
\cr
&  & X_{\Delta}  &  \subset & \Bbb P^{t+1} \cr
}
$$
 therefore $\Phi - P$ is a triple point of $X'= X_{\Phi + \Delta}$,  in contrast  with  6.1. \hb
$ii)$  The result is obvious for $t=0$, so we can assume that $t>0$. \hb
Since $\langle \Phi \rangle$ is a fibre of $R_{1,t+1}$, then  the projection centered
in the line $\langle \Phi \rangle$ maps $R_{1,t+1}$ onto a cone:
$$
\eqalign{
\pi_{\langle \Phi \rangle} : & \; \Bbb P^{t+3} \; \;  \longrightarrow \Bbb P^{t+1}\cr
& R_{1,t+1} \; \mapsto \; R_{0,t}   . \cr
}
$$
Moreover, recalling that $H_{X'} \sim \Phi + \Delta$, we have $\pi_{\langle \Phi \rangle} (X') =
X_\Delta = \varphi_\Delta(X) \subset R_{0,t}$. Since all the
singularities of $X'$ belong to $C_0$ (see 4.2), then necessarily $X_\Delta$ has only
one singular point in $C:=\pi_{\langle \Phi \rangle}(C_0)$, which is the vertex of the
cone
$R_{0,t}$. \hb
In order to obtain a linear series of dimension 1 on $X_\Delta \subset \Bbb P^{t+1}$, it
is necessary to project it from $t$ points, say $A_1, \dots, A_t$, of $X_\Delta$. If each
of these points if different from $C$, then we get the required $B \in g^1_\beta$, where
$\beta =
\deg(\Delta) - t = \lambda$. If, for some $i$, it occurs that $A_i = C$, then
$\pi_C(R_{0,t}) = {\cal C} \subset \Bbb P^t$, where ${\cal C}$ is a rational normal curve
of degree $t$: in this case $B \in g^1_4$,  in contrast with the  assumption $g^1_\beta \ne g^1_4$.
\cvd

\medskip
\noindent
{\bf Definition.} A linear system $|\Delta|$ on the curve $X$ is called {\it minimal}
if it satisfies the conditions $i)$ and $ii)$ of 6.4.

\medskip
\noindent
{\bf Remark 6.5.} Note that, if we perform the previous construction with respect to a
birational model $S_0 \in \grc_L(S)$ which is not a standard model, then the
corresponding series $|\Delta|$ is primitive but not minimal. 

\medskip
\noindent
{\bf Remark 6.6.} If $t=0$, i.e.  $|\Delta| = g^1_\lambda$, then $|\Delta|$ is minimal
if and only if is primitive.

\medskip
We have seen in 6.4 that, if $R_{1,t+1}$ is isomorphic to a standard model, then the
associated series $|\Delta|$ on $X'$ is minimal. The converse is also true, as the
following result shows.

\bigskip
\proclaim
Proposition 6.7.  Let $X$ be as usual and consider two divisors $\Phi \in g^1_4$ and $\Delta \in g^{1+t}_{\lambda +t}$. If  the  linear series $|\Delta|$  is minimal on $X$, then
$ X_{\Phi + \Delta} \subset R_{1,t+1}$ is isomorphic to a standard model of $X_K \subset S$.
\par
\pf We have to consider two cases: either 
$\dim \langle \varphi_{\Phi + \Delta}(\Phi) \rangle = 1$ or 
$\dim \langle \varphi_{\Phi + \Delta}(\Phi) \rangle = 2$. \hb
$(1)$ In this case, since $\deg(\Phi) = 4$, then   $X_{\Phi + \Delta}$ is contained in a geometrically ruled surface as a four--secant curve. Moreover, since $\dim |\Delta| = t+1$, then the invariant of such ruled surface is $t$. Therefore 
$X_{\Phi + \Delta} \subset R_{h,t+h}$ for a suitable $h \ge 1$. \hb
Assume first that  $h \ge 2$. With a construction as in the proof of 6.4 $(ii)$,
consider the projection
$$
\pi_{\langle \Phi \rangle} : \;  R_{h,t+h} \longrightarrow  R_{h-1,t+h-1}
$$
where $\pi_{\langle \Phi \rangle}(X_{\Phi + \Delta}) = X_{\Delta}$. \hb
Note that $H_R \sim U+ hf$, where $U$ is a unisecant curve of degree $t+h$. Therefore,
as noted in 6.2,
$$
\Phi + \Delta = H_R \cdot X_{\Phi + \Delta} \sim
h \Phi + U \cdot X_{\Phi + \Delta}.
$$ 
Since $h \ge 2$, it follows that $\Delta \sim 
(h-1) \Phi + U \cdot X_{\Phi + \Delta}$, so
$\Phi \subset \Delta$.  Hence $\Delta - \Phi \in g^{t-1}_{\lambda+t-4}$.
Therefore there exist $t-2$ points, say $A_1, \dots, A_{t-2}$, such that 
$\Delta - \Phi - A_1- \cdots - A_{t-2} \in g^1_{\lambda-2}$. But this is impossible
since $|\Delta|$ is minimal, hence it satisfies $(ii)$ of 6.4. 
This proves that $h=1$, so  $X_{\Phi + \Delta} \subset R_{1,t+1}$. \hb
If $X_{\Phi + \Delta}$ has a multiple point $P$ not belonging to $C_0$, then we can
project it from $P$ and $t-1$ general points of the curve, obtaining a divisor $B \subset
\Delta$ such that $B \in g^1_{\overline \lambda}$ and $\overline \lambda < \lambda$.
Therefore all the singular points of $X_{\Phi + \Delta} \subset R_{1,t+1}$ belong to
$C_0$ and this implies (from 4.2) that $R_{1,t+1}$ is a standard model. 
\smallskip
\noindent
$(2)$ In this case the curve is contained in the scroll $V$, ruled by planes, whose
fibers are $\langle \varphi_{\Phi + \Delta}(\Phi) \rangle$, $\Phi \in g^1_4$. So we set,
for suitable $a \le b \le c$:
$$
X_{\Phi + \Delta} \subset V = \Bbb P({\cal O}(a) \oplus {\cal O}(b) \oplus {\cal
O}(c)).
$$
Clearly, among the unisecant curves $U^b$ of degree $b$ such that $U^b \subset R_{a,b}
\subset V$, we can choose one of them, say $U$, which does not meet 
$X_{\Phi + \Delta}$ (otherwise $X_{\Phi + \Delta}$ would be contained in the ruled
surface $R_{a,b} \subset V$, against the assumption). Therefore, if we consider 
the projection
$$
\pi_{\langle U \rangle}: \quad V \longrightarrow R_{a,c}
$$
it is clear that $\pi_{\langle U \rangle}(X_{\Phi + \Delta})$ is again a curve, say
$\xx_{\Phi + \Delta}$, whose hyperplane divisor is still $\Phi + \Delta$, but
 $\xx_{\Phi + \Delta} \subset R_{a,c}$,  contrary to  the assumption as well.
\cvd

\bigskip
The remaining part of this Section is devoted to the case $t=0$. Here the linear series $|\Delta|$ will be denoted by $|\Lambda|$, since its degree is $\lambda$, as noted in 6.6. \hb
We will show that this linear series  is, in general,  not unique. In order to determine all such series $g^1_\lambda$, let us describe the situation and notation.

\medskip
Let $X_K \subset S \subset V$ be as usual and assume that $t(S)= 0$.
Let $\Phi \in g^1_4$, $\Lambda' \in g^1_{\lambda'}$ (where $\lambda' >4$) and 
$X_{\Phi+\Lambda'} := \varphi_{\Phi+\Lambda'}(X) \subset R_{1,1}$. Denote by $|l|$ and $|l'|$ the two
rulings of $R_{1,1}$. 

\medskip
\noindent
{\bf Notation.} If $P \in R_{1,1}$, denote by $l_P$ and $l'_P$ the lines of the two rulings passing through $P$. Moreover, if $A$ is a double point of  $X_{\Phi+\Lambda'}$,  denote by $A_1$ and $A_2$ the corresponding points on the canonical model of the curve, i.e. $A_1, A_2 \in X_K$ are such that  $\varphi_{\Phi+\Lambda'}(A_1) = \varphi_{\Phi+\Lambda'}(A_2) = A$.

\medskip
\proclaim
Proposition 6.8.  In the above situation,   each   pair of double
points,   $A$ and $B$ say,  of $X_{\Phi+\Lambda'}$  such that $l_A \ne l_B$ and $l'_A \ne l'_B$, determines a linear series $|\overline \Lambda'| \ne |\Lambda'|$ of degree
$\lambda'$.
\par
\pf Consider the four--gonal divisors and the $\lambda'$-gonal divisors of $|\Lambda'|$ containing, respectively, the
two double points, i.e.
$$
\eqalign{
A_1 + A_2 + A'_1 + A'_2\in g^1_4, \quad
& A_1 +A_2+ P_1 + \cdots + P_{\lambda'-2} \in |\Lambda'| \cr
B_1 + B_2 + B'_1 + B'_2 \in g^1_4, \quad
& B_1 + B_2 + Q_1 + \cdots + Q_{\lambda'-2} \in |\Lambda'| \cr}.
$$
Consider the divisor $\overline \Lambda' = \Phi + \Lambda' - (A_1 + A_2+B_1 + B_2)$; it is clear
that $|\overline \Lambda'|$ is a linear series of degree $\lambda'$ which is distinct
from $|\Lambda'|$. 
\cvd

\medskip
\noindent
{\bf Remark 6.9.}  Let $X_K \subset S$ be as usual and assume that $t=0$ and $\lambda$ are the invariants of $S$. Let  $\Phi \in g^1_4$, $\Lambda \in g^1_{\lambda}$ be two divisors on $X$.
 In the general case, the
$\delta$ double points of $X' = X_{\Phi + \Lambda} \subset R_{1,1}$  belong to different lines of the two rulings $|l|$ and 
$|l'|$. Therefore from the above result it is clear that there are $\delta \choose 2$ linear series $|\Lambda|$ of degree
$\lambda$; to each of them we can associate a model of $X$ lying on $R_{1,1}$. 
In particular, if $|\overline \Lambda|$ is one of these series,  the  corresponding model 
$X_{\Phi+\overline \Lambda}$  still has $\delta$ double points since the   pair  $(A,B)$ has been replaced by $(A',B')$, where 
$A' := \varphi_{\Phi+ \overline\Lambda}(A'_1) = \varphi_{\Phi+ \overline\Lambda}(A'_2) $ and 
$B' := \varphi_{\Phi+ \overline\Lambda}(B'_1) = \varphi_{\Phi+ \overline\Lambda}(B'_2) $, following the notation in 6.8.

\medskip
\proclaim
Theorem 6.10. Let $X_K \subset S \subset V$ and let $S$ be a surface ruled by conics of
minimum degree. Let $t$ and $\lambda$ be the invariants of $S$ defined before. If
$t=0$ then the invariant
$\lambda$ is the minimum degree of a linear series distinct from the
 $g^1_4$, i.e.
$$
\lambda = \min \{r \; | \; X \; 
\hbox{has a complete and base--point--free linear series $g^1_r$ and}
\; r > 4\}.
$$
Moreover, assume that $|\Lambda|$ and $|\Lambda'|$ are two distinct linear series of degree $\lambda$
and let $S$ and $S'$ be the associated surfaces. Then the following facts hold:
\item{$(i)$} if  $\lambda \ne {{g+3} \over 2}$, then $S= S'$;
\item{$(ii)$}  if $\lambda = {{g+3} \over 2}$,  then $S$ and
$S'$ are not necessarely coincident but belong to a pencil of surfaces, ruled by conics, each of them   associated  to a linear series of degree $\lambda$ and has degree ${{3g-7} \over 2}$.
\par
\pf
Recall that $\lambda$ is defined at the beginning of this Section as the invariant of $X$ such that a standard model of $X$ is a divisor of type $(4,\lambda)$ on $R_{1,1}$.
Consider a linear series
$g^1_{\lambda'}
\ne g^1_{\lambda}$;
we  need  to show that $\lambda' \ge \lambda$. Suppose that $\lambda' < \lambda$. \hb
If $g^1_{\lambda'}$ is minimal, consider $\Lambda' \in g^1_{\lambda'}$. Clearly,
$X_{\Phi+\Lambda'} \subset R_{1,1}$ is a standard model. \hb
If $g^1_{\lambda'}$ is not minimal, then it is not primitive (from 6.6); so there
exist $t'$ points, say $A_1, \dots, A_{t'}$ such that $\Delta :=\Lambda' + A_1 + \cdots +
A_{t'}$ is both primitive and minimal. Therefore $X_{\Phi+\Delta} \subset R_{1,t'+1}$ is a
standard model. Hence the corresponding surface $S'$ ruled by conics is such that
$X_K \subset S' \subset V$ and $\deg(S')=g+\lambda'-t'-5$. Assume that $S' \ne S$;
since $X_K \subseteq S \cap S'$, by $(IF)$ we have:
$$
\deg(X_K) \le \int_V S\cdot S'=2\ \deg(S)+2\  \deg(S')-4 \ \deg(V)
$$
hence
$$
2g-2\le 2(2g+\lambda+\lambda'-t-t'-10)-4(g-3)
\quad \Rightarrow \quad
\lambda+\lambda' \ge t+t'+g+3.
$$
Since  $\lambda' < \lambda$ then the above relation gives:
$$
\lambda > {{g+3}\over 2} + {{t+t'}\over 2} ={{g+3}\over 2} + {{t'}\over 2}
$$
where the last equality comes from the assumption $t=0$. \hb
On the other hand, $\lambda \le  {{g+3} \over 2}$ from 5.4. Hence $t' < 0$ and this is
impossible. \hb
 Therefore we have proved that, if $S' \ne S$ then $\lambda' \ge \lambda$.\hb
Assume now that $S' = S$. Clearly, $t'=t=0$ and $\deg(S) = \deg(S')$. Hence, from 4.4,
it follows  that $\lambda = \lambda'$.  \hb
In this way, we have proved the first part of the statement. 
\smallskip
\noindent
\item{$(i)$}
Assume now that $\lambda \ne {{g+3} \over 2}$ and
$S \ne S'$. Then we can use the $(IF)$ as before and, from the assumption $\lambda =
\lambda'$, we obtain 
$$
\lambda \ge {{g+3}\over 2} + {{t'}\over 2}. 
$$
Again we apply 5.4 to $S$, so:
$$
\lambda \le {{g+3} \over 2}.
$$
Comparing these inequalities, we obtain:
$$
t'=0
\quad \hbox{hence} \quad 
\lambda =  {{g+3} \over 2} 
$$
 contrary to the assumption.
\item{$(ii)$} Suppose now that $\lambda =  (g+3)/2$. In this case, from 4.4,
$$
\deg(S) = g+\lambda  -5= {{3g-7} \over 2}.
$$
Therefore 
$$
\deg(S') = g-\lambda-t'-5 \le \deg(S)
$$ 
and this implies $t'=0$ and 
$$
\deg(S') = \deg(S) ={{3g-7} \over 2}.
$$
So, by 1.1,  the result follows.
\cvd

\bigskip
\noindent
{\bf 7. Bounds for the invariants  $a$ and $b$}
\medskip
In this section we determine the range of the invariants
$a$ and $b$ of the four--gonal curve $X$. \hb
Let us keep the notation of Section 5, where $\xx_0 \subset \ss_0\subset \vv$ are standard models of $X_K \subset S \subset V$ and $\pi: \Bbb P^{g-1+\delta} \longrightarrow \Bbb P^{g-1}$ is the projection centered  on  the singular locus of $\xx_0$. \hb
Recall also that $V= \Bbb P({\cal O}(a) \oplus {\cal O}(b) \oplus {\cal O}(c))$ and
$\vv=V_{\ss}= \Bbb P({\cal O}(\lambda-2-2t) \oplus
{\cal O}(\lambda-2-t) \oplus
{\cal O}(\lambda-2) )$. Moreover, from 2.2 $(iii)$, we have
$\delta= 3(\lambda -t-1) -g$
and,  from 5.4, we obtain the following range of the  invariant $\lambda$:
$$
{{g+3}\over 3}  \le \lambda -t \le 
{{g+3}\over 2} .
\eqno(10)
$$
 
\medskip
\noindent
{\bf Remark 7.1.} Note that, from the above expression of $\vv$, it follows that $a \le
\lambda -2-2t$, $b \le \lambda-2-t$, $c \le\lambda -2$. Moreover, since $a+b+c
= g-3$, there are only two independent invariants, $a$ and $b$   say.

\medskip
\noindent
{\bf Notation.}  Clearly, if $a<b$, there exists a unique
directrix on $V$ having degree $a$. In this case, 
let us denote by $A$ such directrix of $V$, 
 by $\aa \subset \vv$ the preimage of $A$ via $\pi$, by
$\delta_A$  the number of the double points (possibly  infinitely near) of $\xx_0$
lying on $\aa$ and by $\overline a$ the degree of $\aa$.   Then
$$
a = \overline a - \delta_A.
\eqno(11)
$$

\medskip
\proclaim
Proposition 7.2.  Let $t>0$ and $U$ be a directrix on $\ss_0$. If
$\deg(U) <
\lambda -2$, then
$U = C_0$.
\par
\pf  It is enough to consider the isomorphism
$$
\varphi_{2C_0+(\lambda - 2)f}: \quad \hs{t}   \longrightarrow \ss_0
$$
and the unisecant irreducible curves $C_0$ and $U=C_0 + \al f$ on $\hs t$. \hb
If $U \ne C_0$, then $\al \ge t$ from 0.1. So
$$
\deg_{\ss_0}(U) = \int_{\ss_0} \; (C_0 + \al f) \cdot (2C_0+(\lambda - 2)f)
= \lambda -2+ 2 \al -2t \ge  \lambda -2
$$
 and the result follows. 
\cvd

\goodbreak

\medskip
\proclaim
Proposition 7.3. Let $t \ge 0$. Then the directrix $\aa$ of $\vv$ is contained
in $\ss_0$
\par
\pf Assume that $\aa \not \subset \ss_0$.  Then,  taking into account that $\deg(\ss_0) = 4(\lambda -t-2)$ as computed in 4.4 and $\deg(\vv) = 3(\lambda -t-2)$, using the Intersection Formula we have:
$$
\int_{\vv} \; \xx_0\cdot \aa \le \int_{\vv} \; \ss_0\cdot \aa =
\deg(\ss_0) + 2 \deg(\aa) - 2 \deg(\vv)
=2 \overline a  - 2 \lambda + 2t + 4.
$$
Therefore, if the $\delta_A$ singular points are distinct, it  follows that:
$$
\delta_A \; \le \; {1\over2} \int_{\vv} \; \xx_0\cdot \aa =
\overline a  - \lambda + t + 2.
$$
In the case of infinitely  near points, it is not so difficult to show that the same
relation holds. \hb 
In this way, from $(11)$, we have the following bound  of $a$:
$$
a = \overline a - \delta_A \ge \lambda -t -2 ,
$$
which is the minimum degree of a directrix of $V$. \hb
Consider the directrix $\pi(C_0) \subset V$. Since $\deg_{\vv}(C_0) =\lambda -2t
- 2$ and the center of $\pi$ contains at least one point of $C_0$, then
$\deg_V(\pi(C_0)) \le \lambda -2t  -3 < \lambda -t -2$; this  concludes the proof.
\cvd

\bigskip
 Next we determine bounds for the invariant $a$.

\medskip
\noindent
{\bf Remark 7.4.} 
Consider the unisecant $\aa \subset \ss_0 \cong \hs{t}$. Clearly, from 0.1, we have:
$$
\aa \sim C_0 + \al f, \quad \hbox{for some $\al \ge t$ or $\al =0$.}
$$
Therefore, as computed in the proof of 7.2,  we have:
$$
\overline a  = \deg_{\ss_0} (\aa) =  \lambda - 2t +2
\al -2
\eqno{(12)}
$$

$$
\aa \cdot \xx_0   = \int_{\ss_0} (C_0 + \al f)(4C_0 + (\lambda +t) f) =  
\lambda - 3 t + 4\al
$$

$$
\delta_A \le {{\aa \cdot \xx_0} \over 2} = {{\lambda - 3 t + 4\al} \over 2}. 
\eqno{(13)}
$$
It is immediate to see that, from $(11)$, $(12)$ and $(13)$:
$$
a = \overline a - \delta_A \ge 
{\lambda-t-4 \over 2}.
\eqno{(14)}
$$
 Note that this bound of $a$ does not depend on $\al$.

\medskip
\noindent
{\bf Remark 7.5.} Note that, since $\delta_A \le \delta$, from $(11)$ we have:
$$
a = \overline a - \delta_A \ge \overline a - \delta
$$
so, taking into account that $\delta = 3(\lambda - t -1) -g$, from $(12)$  we immediately obtain
$$
a  \ge 
\lambda - 2t +2\al -2 - 3(\lambda -1-t) +g = g-2\lambda +t + 2\al +1 \ge
g-2\lambda +t  +1.
\eqno{(15)}
$$

\medskip
\noindent
{\bf Remark 7.6.} In order to compare the two  bounds  of $a$ given by
$(14)$ and $(15)$, just note that
$$
{\lambda-t-4 \over 2} < g-2\lambda +t  +1
\quad \Leftrightarrow \quad
\lambda < {2g+3t+6\over5}.
$$
This leads us to consider the best lower bound of $a$ in each of the two ranges of
$\lambda$.

\bigskip
Keeping into account the previous remarks, we have immediately:

\goodbreak

\medskip
\proclaim
Proposition 7.7. The invariant $a$ has the following lower bound:
$$
a_{\min} := a_{\min}(g,\lambda,t) = 
\left\{
\matrix{
\displaystyle{\Big\lceil {\lambda -t-4 \over 2} \Big\rceil} & \hbox{if} &
\lambda
\ge {2g+3t+6\over5}
\cr 
\cr
\displaystyle{g-2\lambda +t  +1} & \hbox{if} & \lambda \le {2g+3t+6\over5} \cr 
} \right. 
$$
and these bounds are attained if and only if $\aa = C_0$.
\cvd

\medskip
\noindent
{\bf Remark 7.8.}
We can also obtain an ``absolute'' lower bound of $a$, just observing that
$a_{\min}$  can be realized when $\delta_A = \delta$ hence when 
${\lambda-t-4 \over 2} = g-2\lambda +t  +1$ or, equivalently (from 7.6) when 
$\lambda = {2g+3t+6\over5}$. \hb
It is immediate to see that, on this line of the plane $(t,\lambda)$ the two functions giving $a_{\min}(g,\lambda,t)$ coincide and are equal to
$$
a_{\min}(g,t) = {g-t-7\over5}.
\eqno{(16)}
$$
Clearly, the minimum value of $a$ is obtained for the maximum value of $t$ (if $t>0$). Therefore, keeping into account that  $\lambda  \ge 3t$ (by 2.2), it is clear that the minimum value of $a$ corresponds to the common point of the lines $\lambda = {2g+3t+6\over5}$ and $\lambda =3t$. 
We  finish  the argument by observing that
$$
{2g+3t+6\over5} =3t 
\quad \Leftrightarrow \quad
t={g+3 \over 6}
$$
and substituting this value in $(16)$ we obtain:
$$
a_{\min}(g)= {g-9 \over 6}.
$$
Note that, in this case, $\lambda = 3t = {g+3 \over 2}$.  Summing up  we have proved that:
$$
\hbox{if} \quad t >0 \quad \hbox{then} \quad
a_{\min}(g) = {g-9 \over 6}, \quad
\hbox{for}
\quad t = {g+3
\over 6} 
\quad \hbox{and} \quad \lambda = {g+3 \over 2}.
$$
Note also that, if $t=0$,  the value of $a_{\min}$  of $(16)$ can be  realized for
$\lambda = {2g+6\over5}$ and we immediately have:
$$
\hbox{if} \quad t =0 \quad \hbox{then} \quad
a_{\min}(g) = {g-7 \over 5}, \quad
\hbox{for}
 \quad \lambda = {2g+6 \over 5}.
$$

\medskip
Therefore, from  7.8, we  obtain:

\medskip
\proclaim
Corollary 7.9.  With the notation above we have: 
$$
\hbox{for all}  \;
t \ge 0, 
\quad
a \ge {g-9 \over 6} \qquad
\hbox{while, if} \; \;
t=0,
\quad
a \ge {g-7 \over 5}.
$$
In particular, $V_S$ is not a cone for $t \ge 0$ and $g \ge 10$ or  $t=0$ and  $g \ge 8$.

\cvd

\bigskip
\proclaim
Proposition 7.10.  Keeping the notation above, the invariants $a$ and $b$
can vary in the following two ranges:
$$
a_{\min} \le a \le {{g -3} \over 3 }
\eqno{(R_2)}
$$

$$
g-\lambda -1 \le a+b \le {{2(g -3)} \over 3 } .
\eqno{(R_3)} 
$$
\par
\pf The two inequalities on the right in
$(R_2)$ and $(R_3)$ follow from $a \le b \le c$ and $a+b+c = g-3$. 
 For  the left inequality
of $(R_3)$, note that  $c \le \lambda -2$ by 7.1, hence
$a+b=g-3-c \ge g-3-(\lambda -2)$, as requested.
\cvd

\medskip
\noindent
{\bf Remark 7.11.} If $a < {{g-\lambda -1} \over 2}$ then  $a<b$,
hence $A$ is unique. 

\goodbreak

\bigskip
\noindent
{\bf 8. Existence of curves of given invariants $\lambda, a,b$  when  $t=0$. }
\medskip
\noindent
{\bf Remark 8.1. }Let us examine  the situation corresponding to
$t=0$. Here a standard model $\ss_0$ of $S$ is isomorphic to the quadric $\hs0$ via
$$
\varphi_{2l+(\lambda - 2)l'}: \quad \hs0   \longrightarrow \ss_0 \subset \Bbb P^{3\lambda -4}
$$
and $\xx_0 \sim 4l + \lambda l'$ on $\ss_0$.
Moreover, the projection from $\vv$ to $V$ is 
$\pi: \Bbb P^{3\lambda -4} \longrightarrow \Bbb P^{g-1}$, 
  $\vv = \Bbb P({\cal O}(\lambda-2)^{\oplus 3})$
and the previous 2.2 $(iii)$, $(10)$, $(R_2)$, $(R_3)$ become, respectively:
$$
\delta= 3(\lambda -1) -g
\eqno{(17)}
$$

$$
{{g+3}\over 3} \le \lambda \le 
{{g+3}\over 2} .
\eqno(R_1)
$$

$$
a_{\min} \le a \le {{g -3} \over 3 }
\eqno{(R_2)}
$$

$$
g-\lambda -1 \le a+b \le {{2(g -3)} \over 3 }
\eqno{(R_3)}
$$
where
$$
a_{\min} = 
\left\{
\matrix{
\displaystyle{\Big\lceil {\lambda -4 \over 2} \Big\rceil} & \hbox{if} &
\lambda
\ge {2g+6\over5}
\cr 
\cr
\displaystyle{g-2\lambda   +1} & \hbox{if} & \lambda \le {2g+6\over5} \cr 
} \right. .
$$

\medskip
\noindent
Note that ${2g+6 \over 5}$ belongs to the range of $\lambda$ given in $(R_1)$. 
Moreover, $\lambda = {2g+6 \over 5}$ if and only if  $\delta = {\lambda \over 2 }$.

\bigskip

At this point, beside the map $\varphi:=\varphi_{2l+(\lambda - 2)l'}$ defined before, it is useful to introduce a further model of $S$ given by the following  isomorphism
$$
\psi:=\varphi_{4l+\lambda l'}: \; \hs0 \longrightarrow S' \subset \Bbb
P^{5\lambda +4}.
$$

\medskip
\noindent
{\bf Notation.} {}From now on, we  denote a geometrically ruled surface
$\varphi_{nl+ml'}(\hs0)\subset \Bbb P^{(n+1)(m+1)-1}$  by $S_{n,m}$.

\medskip
\noindent
In this way, $S' = S_{4,\lambda}$ and we set $f : \;S' \longrightarrow \ss_0$ the
isomorphism   being   given by  $\varphi = f \circ \psi$.

\medskip
\noindent
{\bf Remark 8.2.} A hyperplane section $H \cdot S'$ of $S'
\subset \Bbb P^{5\lambda +4}$ corresponds, via the morphism
$\psi$, to a curve 
$X_H \subset \hs0$ of type $(4,\lambda)$.
It is not difficult to show, using 3.9,   that  $P \in \hs0$ is a double point of $X_H$
if and only if $H$ contains the tangent plane $T_P(S')$ (here $P$
means $\psi(P) \in S'$).

\medskip
\noindent
{\bf Remark 8.3.} Let $S:= S_{n,m} \subset \Bbb P^{(n+1)(m+1)-1}$ and $Y
\subset S$ be a divisor whose decomposition into irreducible and reduced 
components is
$Y = Y_1 \cup \dots \cup Y_s$. Let $P_1, \dots, P_\delta$ be points
of
$Y$ and denote by $\delta_i$ the number of these points  belonging
to the component $Y_i$.  Let  
$$
L:= \Bigg\langle T_{P_1}(S),\dots,T_{P_\delta}(S) 
\Bigg\rangle
$$
be  the linear space spanned by the $\delta$ tangent planes.
Clearly, if $H$ is any hyperplane
containing $L$, then $H$ intersects $Y_i$ in at least $2\delta_i$
points. Therefore, if $2\delta_i > \deg(Y_i)$, then $H $ contains  $Y_i$. 

\medskip
The above observation leads to the following:

\medskip
\noindent
{\bf Definition.} 
We say that $P_1, \dots, P_\delta$  {\it trivially degenerate} the
component $Y_i$ if $2\delta_i > \deg(Y_i)$. Moreover, we say that 
$P_1, \dots, P_\delta$  {\it trivially degenerate} the curve $Y$ if
this occurs for at least one component of $Y$.

\medskip
\noindent
{\bf Remark 8.4.} Let $S' = S_{4,\lambda}$ be as before. Assume that $a \le b \le c$
fulfil the relations
$(R_1),(R_2),(R_3)$.
\item{$(a)$} Let $M \sim l$ be a divisor of $S'$. Clearly $\deg(M) = H \cdot M =
\lambda$.  Let us consider $\lambda - 2-a$ distinct points of $M$, say $P_1,
\dots, P_{\lambda-2-a}$. Clearly $P_1, \dots, P_{\lambda-2-a}$ do not trivially degenerate $M$ if and only
if
$$
2(\lambda - 2 -a) \le \deg(M) = \lambda
\quad \Leftrightarrow \quad a \ge {\lambda -4 \over 2}
$$
and this is true by $(R_2)$.
\item{$(b)$} In the same way, if $N \sim l$ is a divisor of $S'$ and $P_1, \dots,
P_{\lambda-2-b}$ are distinct points of $N$, then 
$$
2(\lambda - 2 -b) \le 2(\lambda - 2 -a) \le \deg(N) = \lambda
$$
again  by $(R_2)$. So $P_1, \dots, P_{\lambda-2-b}$ do not trivially degenerate
$N$.
\item{$(c)$} Consider now a divisor $Q \sim (\lambda -2-c) l'$ consisting of
$\lambda -2-c$ distinct components and  a set of distinct points  $P_1, \dots,
P_{\lambda-2-c}$, one on each component of $Q$. Obviously
 $P_1, \dots, P_{\lambda-2-c}$ do not trivially degenerate $Q$.

\bigskip
\proclaim
Theorem 8.5. Let $g,a,b,\lambda$ be positive integers,  with  $g\ge 10$, and
consider the following inequalities:
$$
{{g+3}\over 3} \le \lambda \le 
{{g+3}\over 2} 
\eqno(R_1)
$$
$$
a_{\min}
 \le a \le {{g -3} \over 3 }
\eqno (R_2)
$$
$$
g-\lambda -1 \le a+b \le {{2(g -3)} \over 3 }
\eqno (R_3)
$$
where
$$
a_{\min} = 
\left\{
\matrix{
\displaystyle{\Big\lceil {\lambda -4 \over 2} \Big\rceil} & \hbox{if} &
\lambda
\ge {2g+6\over5}
\cr 
\cr
\displaystyle{g-2\lambda   +1} & \hbox{if} & \lambda < {2g+6\over5} \cr 
} \right. .
$$ 
Then there exists a $4$--gonal curve of genus $g$ and invariants
$a,b,\lambda$ if and only if $(R_1),(R_2),(R_3)$ are verified.
\par
\pf
If there exists a $4$--gonal curve of genus $g$ and invariants
$a,b,\lambda$ then $(R_1),(R_2),(R_3)$ come from 8.1.  \hb
Conversely, let us choose $g,\lambda,a,b$ satisfying the inequalities
$(R_1)$, 
$(R_2)$, $(R_3)$. Using 8.2, it is enough to show that there exists an
irreducibile hyperplane section $H \cdot S'$ of $S' = S_{4,\lambda}$, i.e. a curve
$X_H \sim 4l +\lambda l'$ on $\hs0$, of genus $g$ and invariants $a,b$.  \hb
Take the following three divisors of $S'$: $M$, $N$, $Q$,
where $M \sim l \sim N$ ($M\ne N$) and $Q \sim (\lambda-2-c) l'$
consists of distinct lines; moreover consider 
$\lambda-2-a$ distinct points of $M$, 
$\lambda-2-b$ distinct points of $N$ and
$\lambda-2-c$ distinct points of $Q$, one on each line and none
belonging to $M$ or $N$. \hb
Note that  $M+N+Q \in |2l +(\lambda-2-c)l'|$ and the equality
$(\lambda-2-a)+(\lambda-2-b)+(\lambda-2-c)=\delta$ holds from $(17)$. \hb
Therefore, taking into account also  8.4, it is immediate to see that 
the hypotesis of the  forthcoming lemma 9.4 are verified; 
then we can deduce  that the linear space $L$ spanned by the tangent
planes to
$S'$ at the above $\delta$ points does not contain any further point of
$S'$. In particular, a general hyperplane $H \supset L$
corresponds to an irreducible curve $X_H \sim 4l+\lambda l'$ having exactly $\delta$
nodes; so its genus is
$g(X_H)=3(\lambda-1) - \delta =g$. \hb 
Consider the isomorphism  $f: \; S'  \longrightarrow  \ss_0$ defined before
and set $\aa:=f (M)$, $\bb:=f(N)$. Clearly
$$
\deg(\aa) =\deg(\bb) = \lambda -2.
$$
Set $\xx_0:=\varphi(X_H)
\subset \ss_0$ and denote by $\delta_A$ and $\delta_B$  the number of the double
points of $\xx_0$ lying on $\aa$ and on $\bb$,
respectively. {}From the construction, it is clear that:
$$
\delta_A = \lambda -2-a
\quad \hbox{and} \quad
\delta_B = \lambda -2-b.
$$
Setting $A,B \subset S \subset V$ the projections of
$\aa$ and $\bb$, respectively, via 
$\pi_{\langle \del \rangle}: \ss_0 \rightarrow S$, from $(11)$ we have that
$\deg(A)  = \deg(\aa) - \delta_A=
\lambda -2- \delta_A = a $ and
$\deg(B) = \deg(\bb) - \delta_B=
\lambda -2-\delta_B = b$. \hb
In this way one can easily deduce that
$V = V_S= \Bbb P(
{\cal O}(a) \oplus 
{\cal O}(b) \oplus 
{\cal O}(c))$,
so $a$ and $b$ are the other two invariants of $X$. 
\cvd

\medskip
In order to complete the proof of the Theorem above, we need to prove the ``Key--lemma'' stated in 9.4. Next section will be devoted to this purpose.

\goodbreak

\bigskip
\bigskip
\noindent
{\bf 9. Proof of the Key--lemma}
\medskip
In order to prove  the Key--lemma 9.4, we need some preliminary technical
results.

\medskip
\proclaim 
Lemma 9.1. Let $S := S_{n,m}$ and $D \sim
hl+kl' \subset S$  be a divisor, where
 $h \le n+1$ and $k \le m+1$.   Then the following facts hold:
\item{$i)$} 
$$
\dim \langle D \rangle= h(m+1)+k(n+1)-hk-1.
$$
\smallskip
\noindent
Moreover, if $D$ is  irreducible: 
\item{$ii)$}  $D$ is a non--special curve;
\item{$iii)$}  $D$ is a linearly normal curve in $\langle D \rangle$. 
\par
\pf 
\item{$i)$} Assume first that $h \le n$ and $k \le m$. It is clear that, setting
$S':=S_{n-h,m-k}$, we have
$\dim \langle D \rangle = h^0({\cal O}_S(1)) - h^0({\cal O}_{S'}(1)) -1$ and
this proves the above relation.\hb 
The remaining cases are: $h=n+1$ and $k \le m+1$ or $h\le n+1$ and $k = m+1$. In both
of them,  $D \sim hl + kl'$ cannot be contained in any
hyperplane section $H \cdot S \sim nl + ml'$ of $S$. Hence 
$\langle D \rangle = \langle S \rangle$, so $\dim \langle D \rangle = \dim \langle
S \rangle = (n+1)(m+1) -1$ and this gives the formula in the statement
when $h= n+1$ or $k=m+1$.
\item{$ii)$} It is enough to show that $\deg(D) > 2p_a(D) -2$. Taking into account
that $\deg(D) = hm+kn$ and $p_a(D) = hk - h - k+1$, and using the assumption
$n \ge h-1$ and $m \ge k-1$, we obtain:
$$
\deg(D) = hm+kn \ge h(k-1) + (h-1)k >2hk -2h -2k = 2p_a(D) -2.
$$
\item{$iii)$} It is enough to prove that $h^0(D, {\cal O}_D(1))= \dim \langle D \rangle +1$. \hb
Since $D$ is non--special, as proved before, applying
the Riemann--Roch Theorem, we obtain
$$
h^0({\cal O}_D(1)) = \deg(D) - p_a(D) +1
$$ 

and this coincides with $\dim\langle
D \rangle +1$, as one can easily verify. Hence $D$ is linearly normal in $\langle D \rangle$.
\cvd

\medskip
\proclaim
Lemma 9.2. Let $S:=S_{2,k}$, where $k \ge 2$, and consider 
$d$ distinct points: $P_1, \dots, P_d \in S$, where $d
\le 2k+1$. Setting $J:=\langle P_1, \dots, P_d \rangle$, if $\dim(J) < d-1$, then
there exists a unisecant curve $U$ on $S$ such that $\# (U \cap \{P_1, \dots, P_d
\})\ge \deg(U) +1$. In particular, $U \subset S \cap J$. 
\par
\pf Assume for simplicity that the considered points belong to distinct fibres of $S'$.
\hb
Since $\dim|l+kl'| = 2k+1 \ge d$, there exists a unisecant curve
linearly equivalent to $l+kl'$ containing $P_1, \dots, P_d$. Therefore we can
find a unisecant,  $U'$  say, of minimum degree containing $P_1, \dots, P_d$.
Clearly, $U' \sim l+ \epsilon l'$, where $\epsilon \le k$; moreover $U' = U+
l'_1 + \cdots + l'_\al$, where $U$ is irreducible, $P_1, \dots, P_{d-\al} \in
U$ and $P_{d-\al+i} \in l'_i \setminus U$, for $i=1, \dots,\al$. Let us show
that $U$ is the required unisecant curve.  Were this not the case, setting
$$
\be:= \deg(U)+1-(d-\al)
$$
it follows  that  $\be > 0$. Consider the linear space
$T:= \langle J, A_1, \dots, A_\be \rangle$,
where $A_j \in U$.  Clearly $U \subset T$, hence $T$ meets each fiber $l'_i$ in
two points: $P_{d-\al+i}$ and $U \cap l'_i$. Since the fibers are conics then,
choosing $B_i \in l'_i$, the linear space
$$
\Sigma: =\langle J, A_1, \dots, A_\be, B_1, \dots, B_\al \rangle
$$
contains $\langle U'\rangle$. Therefore $\dim \langle U'\rangle \le
\dim(\Sigma) \le \dim(J) + \al + \be=\dim(J)  +\deg(U)+1-d+2\al$. On the other
hand, using 9.1,  $\dim \langle U'\rangle = \deg(U') = \deg(U) + 2\al$, so  $\dim(J) \ge
d-1$, against the assumption. \hb
It is not difficult to generalize this proof to the case  where  at most two of the $d$
points belong to the same fibre.
\cvd

\medskip
\proclaim
Lemma 9.3. Let $S:=S_{4, \lambda}$, where $\lambda \ge 4$,  and
$\wdd \in |2l+ \epsilon l'|$ be a bisecant curve on $S$ such that
 $\wdd$ does not contain any fiber of $S$.  
 Consider $d +1$ points $P,P_1,\dots, P_d$ as follows: $P \in S$, 
$P_1,\dots, P_d \in \wdd$ such that they do not trivially
degenerate $\wdd$ and at most two of them belong to the same fibre. Assume that 
$P_1, \dots, P_m$ are double points of $\wdd$ (for $0 \le m \le d$) and 
$P_{m+1}, \dots, P_d$ are simple points of $\wdd$. Let
$$
T:= \langle
P, T_{P_1}(S), \dots, T_{P_m}(S), t_{P_{m+1}}(\wdd), \dots,
t_{P_{d}}(\wdd)
\rangle
$$
where $T_{P_i}(S)$ and $t_{P_i}(\wdd)$ denote the tangent plane to $S$ and the tangent line to $\wdd$, respectively, at $P_i$. \hb
If $\epsilon \le \lambda$ and $d \le \lambda$, then $\dim(T) = 2 d +m$.
\par
\goodbreak

\pf For simplicity, assume that $P \in \wdd$ and $P_1,\dots, P_d$ belong to distinct
fibres of $S$. In this situation, $T \subseteq \langle \wdd \rangle$ and $m \le d \le \epsilon$. \hb
Claim: $T$ is a proper subspace of $\langle \wdd \rangle$. \hb
In order to prove this, observe that, by 9.1 and the assumption $d \le \lambda$, we have
$$
\dim \langle \wdd \rangle =  2\lambda+3\epsilon +1
\ge 2d + 3 \epsilon +1.
$$
As noted at the beginning, $m \le \epsilon$ hence $\dim \langle \wdd \rangle \ge
2d + 3m+1 > 2d +m \ge \dim(T)$ and this proves the claim.\hb 
Let $N:= \dim \langle \wdd \rangle$ and consider the projection $\pi_T: \, \Bbb P^N \rightarrow \Bbb P^n$ with
center $T$, for a suitable $n$. Clearly, by the claim above, $n>0$. \hb
Let $R:=R(\wdd)$ be the ruled surface generated by $\wdd$ via the ruling on $S$.
Since $T$ is a multisecant space of this ruled surface and $P_1, \dots, P_d$
belong to distinct fibers, then $T\cap R$ contains a unisecant curve (see [{\bf 4}], 1.5),    $Y$ say. Therefore $\pi_T(R) = \pi_T(\wdd)$ is a rational normal curve of degree $n$ in $\Bbb P^n$. In particular:
$$
N-n= \dim \langle \wdd \rangle - \dim \langle \pi_T(\wdd)\rangle
= \dim(T) +1.
\eqno(18)
$$
In order to prove the statement, observe that  it holds that $\dim(T)  \le  2 d +m$. \hb
First case: $\wdd$ is irreducible. \hb
 Since $\pi_{T|\wdd}$ is a map of degree two, then
$$
n = \deg(\pi_T(\wdd))=
{{\deg(\wdd) - \int T\cdot \wdd} \over 2}.
\eqno(19)
$$ 
Moreover, from 9.1 $(iii)$ we have that:
$$
N =\dim \langle \wdd \rangle =h^0({\cal O}_{\wdd}(1)) -1 = \deg(\wdd) - p_a(\wdd)
$$
so, using $(18)$ we finally obtain:
$$
\dim(T) = N-n-1=
\deg(\wdd) - p_a(\wdd) - {{\deg(\wdd) - \int T\cdot \wdd} \over 2} -1=
{{\deg(\wdd) + \int T\cdot \wdd} \over 2} - p_a(\wdd) -1.
$$
Note that $\deg(\wdd)= 4\epsilon + 2 \lambda$ and $p_a(\wdd)=
\epsilon-1$; moreover, by
the definition of $T$,   $\int T\cdot \wdd \ge 2d + 2m+1$. Hence we obtain
$$
\dim(T) \ge \epsilon +  \lambda+d + m  +1/2.
$$
 Thus, if we assume $\dim(T) < 2d+m$,  we get
$$
\epsilon +  \lambda+d + m  +1/2 < 2d + m
\quad \Rightarrow \quad
d > \lambda + \epsilon +1/2
$$
 contrary to  the assumption $d \le \lambda$. \hb
Second case: $\wdd$ is reducible. \hb
Let
$\wdd = U_1 + U_2$, where $U_i$ are irreducible unisecant curves. Let  $d_i$ be  the number of points among $P_1, \dots, P_d$ belonging
to $U_i$. Clearly, $P_1, \dots, P_m$ belong to $U_1 \cap U_2$, so $d=d_1 + d_2 - m$.  
 Moreover, we have
$$
\dim \langle \wdd \rangle = \dim \langle U_1 \rangle +
\dim \langle U_2 \rangle - \int U_1 \cdot U_2 +1.
\eqno(20)
$$
Since $T$ is a proper subspace of $\langle \wdd \rangle$ as proved in the previous claim, then $\wdd \not \subset T$; therefore only two cases can occur: either  $U_i \not \subset T$ for
$i=1,2$ or (for instance) $U_1 \subset T$ and $U_2 \not \subset T$. \hb
If $U_i \not \subset T$ for $i=1,2$, then $  \pi_T(\wdd) =
\pi_T(U_1)= \pi_T(U_2)$ so
$$
n= \dim \langle \pi_T(\wdd)\rangle = \dim \langle \pi_T(U_i)\rangle = \deg (\pi_T(U_i))=
\deg( U_i) - \int T \cdot U_i
\quad \hbox{for} \quad i=1,2 .
\eqno(21)
$$
Adding the previous relations $(21)$  for $i=1$ and $i=2$, we obtain that $2n =
\deg(U_1+U_2) - \int T \cdot (U_1+U_2)$, so this equality coincides with
$(19)$ and we conclude the proof as in the first case.\hb
We are left to study the case $U_1 \subset T$, i.e. $U_1 = Y$. Since $T$ contains the tangent lines to $U_2$ at all the $d_2$ points defined before  and since  $U_1 \subset T$ and the $m$ double points of  $\wdd$ belong to $U_1 \cap U_2$, then
$$
\int  T \cdot U_2 = 2d_2 + \int U_1 \cdot U_2 - m.
$$
In this case
$(21)$ holds only for $U_2$, so it becomes:
$$
\dim \langle \pi_T(\wdd)\rangle = 
\deg(U_2) - \left(2d_2 + \int U_1 \cdot U_2 - m \right).
$$
Therefore,  using the relation above and $(20)$, and taking into account that $\dim \langle U_i \rangle = \deg(U_i)$, we obtain:
$$
\dim \langle \wdd \rangle - \dim \langle \pi_T(\wdd)\rangle = \deg(U_1) + 2d_2 -m +1.
$$
Now we substitute $d_2 = d +m -d_1$ and use $(18)$,  obtaining
$$
\dim(T)+1 = \deg(U_1) + 2d+2m-2d_1 -m +1.
$$
Finally recall that the $P_i$'s do not trivially degenerate $\wdd$, hence
 $2d_1 \le \deg(U_1)$; so we obtain
$$
\dim(T)+1  \ge 2d + m +1
$$
as required. In the general case, the proof runs in a similar way.
\cvd

\bigskip
\noindent
{\bf Notation.} Since we will consider, in the following result, both
$S':=S_{4,\lambda}$and $S_{2,c+2}$, we denote the divisors on these surfaces by: 
$D_4, \wdd_4, \dots$ and $D_2, \wdd_2, \dots$, respectively.

\medskip
\proclaim
Key--Lemma 9.4.  Let $g,a,b,c,\lambda$ be positive integers satisfying
$(2)$, $(R_1),(R_2),(R_3)$. \hb
Let $S':= S_{4,\lambda} \subset \Bbb P^{5\lambda +4}$ and 
$D_4 \in |2l+(\lambda-2-c)l'|$ be a curve on $S'$ of type 
$$
D_4 = \wdd_4 + \sum_{i=1}^{\al} l_i'
$$ 
where $\al$ is an integer such that $0 \le \al \le \lambda -2-c$ and $\wdd_4$ is a
suitable bisecant divisor not containing any irreducible component linearly equivalent
to $l'$. \hb 
Let us take $\delta=3(\lambda-1)-g$ distinct points on $D_4$ which do not trivially degenerate $D_4$ and set
$$
P_1, \dots, P_{\delta-\al} \in \wdd_4 
\quad \hbox{and} \quad P'_1, \dots, P'_{\al} \in \sum_{i=1}^{\al}l_i'
$$ 
such that $P'_i \in l_i' \setminus \wdd_4$ for $i=1,
\dots,\al$.  Consider the linear space
$$
L:= \big\langle T_{P_1}(S'),\dots,T_{P_{\delta - \al}}(S'),
T_{P'_1}(S'),\dots,T_{P'_\al}(S') 
\big\rangle
$$ 
spanned by the tangent planes to $S'$ at these $\delta$ points.  \hb
If $P \in S'$ is any further point such that $P \not \in L$ and 
$L':=\langle P, L \rangle$, then:
$$
\dim(L') = 3 \delta.
$$
In particular, $\dim (L) =3\delta -1$, i.e. $L$ is of maximum dimension and the intersection of 
$L$ and $S'$ consists only of the points  $P_1, \dots, P_{\delta-\al}, P'_1, \dots, P'_{\al} $. 
\par
\noindent
\pf Note first that $\dim(L') \le 3 \delta$ and  $\dim(L) \le 3 \delta -1$. So it is enough
to show that $\dim(L') \ge 3 \delta$. \hb
Assume first that $P \not \in \wdd_4$. 
\smallskip
\noindent
{\bf Step 1.} Computation of the dimension of 
 $\Sigma:= \langle L',D_4\rangle$. \hb
Among the choosen points $P_1, \dots, P_{\delta - \al} \in \wdd_4$, consider those which
are singular points of $\wdd_4$, say $P_1, \dots, P_m$, for some
$0 \le m \le \delta- \al$.

 \bigskip
 \centerline{
 \epsfxsize=10cm 
 \epsfbox{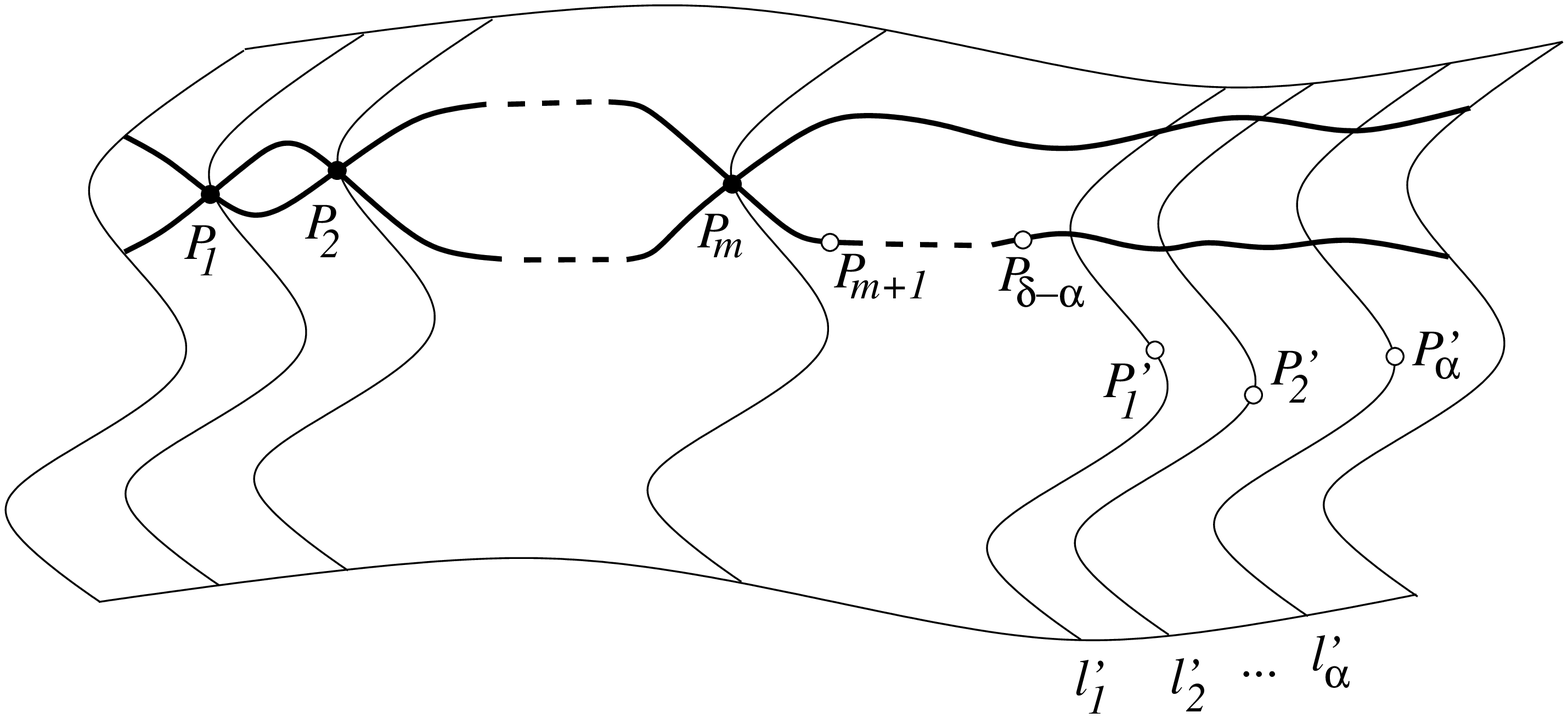}
 }
 \medskip
 \centerline {Figure 5}
  \medskip

\noindent
Clearly, since they are double points of $\wdd_4$, the tangent plane at each of them is
contained in $\langle \wdd_4 \rangle$. On the other hand,  the tangent plane at the
remaining $\delta -m$ points intersects $\langle D_4\rangle$ in a line (either tangent to $\wdd_4$ for $P_{m+1}, \dots, P_{\delta -\al}$, or tangent to $l'_i$ for the points of type $P'_i$). Briefly:
$$
\eqalign{
T_{P_i}(S')  \subset  \langle \wdd_4 \rangle, \; \,  \qquad \qquad \qquad \qquad \qquad
&
\hbox{for} \;  i=1, \dots,m \cr
T_{P_i}(S') \cap \langle D_4\rangle  = t_{P_i}(D_4) = t_{P_i}(\wdd_4), \qquad &
\hbox{for}
\; i=m+1,\dots,\delta-\al \cr 
T_{P'_j}(S') \cap \langle D_4\rangle  = t_{P'_j}(D_4) = t_{P'_j}(l_j'), \qquad &
\hbox{for} \; j=1, \dots, \al . \cr
} \eqno {(22)}
$$
Consider now the projection
$$
\pi:=\pi_{\langle D_4\rangle} : \; S'=S_{4,\lambda} \longrightarrow
S_{2,c+2}
$$
and set
$$
J:= \pi(\Sigma) = \langle \pp, \pp_{m+1}, \dots, \pp_{\delta - \al},
\pp'_{1}, \dots, \pp'_{\al}\rangle \quad 
$$
where
$$
\pp:= \pi(P), \;  \pp_i:= \pi (T_{P_i}(S')),
\; \hbox{for $i= m+1, \dots, \delta - \al$, and} \quad
\pp'_j:= \pi (T_{P'_j}(S') ), \;  \hbox{for  $j=1, \dots, \al$.} 
$$
By the definition of $J$, we clearly have:
$$
\dim(\Sigma) = \dim(J) +\dim \langle D_4 \rangle +1 .
\eqno(23)
$$
{\bf Step 2.} Computation of the dimension of $J$. 
\smallskip
\noindent
Observe that the isomorphisms $\varphi_{4l+\lambda l'}$ and $\varphi_{2l+(c+2) l'}$
induce a canonical isomorphism, say $\chi$, as follows
$$
\matrix{
 \hs 0 \; \cr
{}^{\varphi_{4l+\lambda l'}}\swarrow \qquad  \searrow^{\varphi_{2l+(c+2) l'}} \cr
\cr
S_{4,\lambda} \quad \mapright \chi ,30 \quad S_{2,c+2} \cr
}
$$
and $\chi$ coincides with $\pi$ on $S_{4,\lambda}  \setminus D_4$. \hb
Therefore, setting $D_2 := \chi (D_4) \subset S_{2,c+2}$, the points 
$\pp_{m+1}, \dots, \pp_{\delta - \al}, \pp'_{1}, \dots, \pp'_{\al}$ belong to
$D_2$. \hb
Clearly, $\dim(J) \le \delta -m$. We want to show that $\dim(J) = \delta -m$. \hb
Assume  that $\dim(J) < \delta -m$. In order to apply 9.2, we need to compare the
number of points spanning $J$ with the integer $c$.
\hb
On one hand, from $(17)$ and $(R_1)$ we have:
$$
\delta =3(\lambda -1) -g \le {g+3 \over 2}.
$$
On the other hand, from $(R_3)$, we get $c \ge {g-3 \over 3}$, i.e. $g \le 3c +3$.
Therefore we obtain:
$$
\delta -m \le \delta \le {g+3 \over 2} \le {3c+6 \over 2} < 2c +5 \quad \Rightarrow
\delta -m +1 \le  2(c+2) +1.
$$
So, we can apply Lemma 9.2  to $J$ (which is spanned  by
$\delta - m+1$ points and has dimension smaller than $\delta - m$) and
$S_{2,c+2}$. In this way we  obtain that there exists a unisecant curve  
$\uu \subset J \cap S_{2,c+2}$ such that, setting $r$ the number of the points among 
$\pp, \pp_{m+1}, \dots, \pp_{\delta - \al}, \pp'_{1}, \dots, \pp'_{\al}$ belonging to
$\uu$, then 
$$
\deg(\uu) \le r-1.
$$
Let $\uu \sim
l+\epsilon l'$; then $\deg(\uu) = c+2+2\epsilon$. \hb
{\bf Claim.} The unisecant $\uu$ is not contained in $D_2$. \hb
If not, let $U:=
\chi^{-1}(\uu)$ and $h$ be the number of the points among $P$, the $P_i$'s and
the $P'_j$'s belonging to $U$. On one hand, since these points do not trivially
degenerate $D_4$ (by assumption) and $U \subset D_4$ (since $\uu \subset D_2$ by the assumption of the Claim), then $2h \le \deg(U)$. \hb
On the other hand,  $h \ge r$ by the definitions of $h$ and $r$ and from $\chi(U) = \uu$. {}From all these observations,  it follows
$$
\deg(U) \ge 2h \ge 2r \ge 2(\deg(\uu)+1) = 2(c+3+2\epsilon).
$$
Since  $\deg(U) = \lambda +4 \epsilon$, we obtain $2c + 6 \le \lambda$.  Using
the bound $c \ge (g-3)/3$,  we finally get $\lambda \ge (2/3)g + 4$, against $(R_1)$.
In this way the claim is proved.
\smallskip
\noindent
Since $\uu$ is not contained in $D_2$, we can consider their intersection, which
surely contains the $r$ points introduced before. So
$$
r \le \int_{S_{2,c+2}} \uu \cdot D_2= (l + \epsilon l') \cdot (2l + (\lambda -2-c)l')
= 
\lambda -2-c+2\epsilon.
$$
 The above relation and  $\deg(\uu) \le r-1$ give:
$$
c+2 + 2 \epsilon = \deg(\uu) \le r-1 \le \lambda -3-c+2\epsilon
$$
so $\lambda \ge 2c+5$ and this leads to a
contraddiction, as in the proof of the claim above. \hb
Hence such unisecant curve $\uu$ does not exist and this implies
$$
\dim(J) = \delta -m. \eqno{(24)}
$$

\medskip
\noindent
{\bf Step 3.} Computation of the dimension of $L'$.
\smallskip
\noindent
Putting together $(23)$ and $(24)$ we finally obtain:
$$
\dim(\Sigma) = \dim \langle D_4 \rangle + \delta - m+1.
\eqno (25)
$$
Now let us compare $\dim(\Sigma)$ with $\dim(L')$. Consider the linear space
$$
T:= \langle
P, T_{P_1}(S'), \dots, T_{P_m}(S'), t_{P_{m+1}}(\wdd_4), \dots,
t_{P_{\delta - \al}}(\wdd_4)
\rangle \subseteq L'.
$$
Note that, from $(R_1)$, we have $g \ge 2\lambda -3$; hence
$$
\delta - \al \le \delta = 3(\lambda-1) - g \le \lambda .
$$
Therefore the assumption in  9.3 are satisfied  by
$S_{4,\lambda}$,
$\wdd_4$ and $T$ with respect to the points $P, P_1, \dots, P_{\delta -\al}$: we then obtain
$$
\dim(T)=2(\delta - \al) +m. \eqno(26)
$$
Since $T \subseteq \langle \wdd_4, P \rangle$ by $(22)$,
 there exist $\be$ points, say $R_1, \dots, R_\be \in \wdd_4$ such that
$\langle T, R_1, \dots, R_\be \rangle$ coincides with $\langle \wdd_4, P
\rangle$, where 
$$
\be = \dim \langle \wdd_4, P \rangle - \dim(T) \le \dim \langle \wdd_4
\rangle - \dim(T) +1.   \eqno(27)
$$
Therefore the linear space $\langle L', R_1, \dots, R_\be
\rangle$ contains $\langle \wdd_4 , P \rangle$, so it meets each fibre $l'_{P'_j}$ (for
$j=1, \dots, \al$) in four points:  two of them are
$l'_{P'_j} \cap
\wdd_4$ and the remaining  ones are $l'_{P'_j} \cap T_{P'_j}(S')$.
Hence, if we add to this space a further point, say $A_j$,  on  each fiber, the obtained linear space contains also the quartic curves $l'_{P'_1}, \dots,
l'_{P'_\al}$, hence the whole divisor $D_4$. In this way we have proved
that
$$
\langle L', R_1, \dots, R_\be, A_1, \dots, A_\al \rangle 
\supset \langle L', D_4 \rangle = \Sigma
$$
so 
$$
\dim(\Sigma) \le \dim(L') + \al +\be.
\eqno(28)
$$
Using $(25)$ and $(28)$ we obtain:
$$
\dim \langle D_4 \rangle + \delta - m +1 = \dim(\Sigma)
\le \dim(L') + \al +\be 
$$
and from this, using $(27)$ we get:
$$
\dim \langle D_4 \rangle + \delta - m +1 
\le \dim(L') + \al + \dim \langle \wdd_4 \rangle  - \dim(T) +1.
$$
Finally, using $(26)$ we obtain:
$$
\eqalign{
\dim(L') & \ge
\delta - m + \dim \langle D_4 \rangle - \dim \langle \wdd_4 \rangle - \al + 2(\delta -\al) + m
= \cr
& = 3\delta - 3\al  + \dim \langle D_4 \rangle - \dim \langle \wdd_4 \rangle =\cr
& = 3 \delta}
$$
where the last equality easily comes from 9.1. \hb
Note that the statement has been proved in the case $P \not \in
\wdd_4$, but the case $P \in \wdd_4$ runs in a similar way, with some cautions. Namely, in Step 1, the main difference concernes the linear space $J:= \pi(\Sigma) = \langle  \pp_{m+1}, \dots, \pp_{\delta - \al},
\pp'_{1}, \dots, \pp'_{\al}\rangle $ obtained from $\Sigma$ by projecting from ${\langle D_4\rangle}$ and the relation $(23)$ still holds. In Step 2, since $\delta -m +1 \le 2(c+2) +1$ then, {\it a fortiori}, it holds $\delta -m  \le 2(c+2) +1$. So also in this case Lemma 9.2 can be applied to $J$, which is spanned by $\delta -m$ points and it is assumed to have dimension smaller then $\delta -m -1$. With the same argument can be proved the analogous of $(24)$ i.e. $\dim(J) = \delta -m -1$. Finally, in Step 3 we obtain the analogous of $(25)$ and precisely $\dim(\Sigma) = \dim \langle D_4 \rangle + \delta - m$. In the following argument the result 9.3 is used; since it holds for any $P$,  also in this case $(26)$ is verified. Now it is immediate to see that $(27)$ becomes $\be = \dim \langle \wdd_4 \rangle - \dim(T)$ and we obtain again that 
$$
\dim \langle D_4 \rangle + \delta - m  = \dim(\Sigma)
\le \dim(L') + \al +\be. 
$$
Using the new form of $(27)$ we finally obtain:
$$
\dim \langle D_4 \rangle + \delta - m \le  \dim(L') + \al + \dim \langle \wdd_4 \rangle - \dim(T)
$$
which leads to the end of the proof as in the general case.
\cvd

\medskip
\noindent
{\bf Remark 9.5.} The result stated in 9.4 holds   also if
at most two of the points $P_{1}, \dots, P_{d}$ belong  to the same fibre.

\medskip
 The following immediately follows from 9.4:

\medskip
\proclaim
Corollary 9.6. For every  curve 
$\dd \sim 2l+(\lambda -2-c)l' \subset \ss_0 \cong  \hs 0$ 
and for every choice of $P_1, \dots, P_\delta \in \dd$ which do not
trivially degenerate $\dd$, there exists a curve
$\xx_0 \subset  \ss_0$ whose double points are exactly 
$P_1, \dots, P_\delta$ and whose characters are
$a,b,\lambda$, where $a+b=g-3-c$.
\cvd
\par

\medskip
We conclude this section with some remark about the construction of the bisecant curves
$D_4$ and $\wdd_4$.
\goodbreak

\medskip
\noindent
Let us consider  a geometrically ruled surface contained in $V$ and having minimum
degree; each of such surfaces corresponds to a quotient of type 
$$
{\cal F}:=
{\cal O}(a) \oplus {\cal O}(b) \oplus {\cal O}(c)
\longrightarrow
{\cal O}(a) \oplus {\cal O}(b)
\longrightarrow 0
\eqno(29)
$$
i.e. it is of the type $R:=R_{a,b}= \Bbb P({\cal O}(a) \oplus {\cal O}(b))$.

\medskip
\noindent
{\bf Remark 9.7.} 
Since the above quotients correspond to the sections of 
${\cal F}(-c)$, tensorizing $(29)$ by ${\cal O}(-c)$ we obtain:
$$
0
\longrightarrow
{\cal O}
\longrightarrow
{\cal O}(a-c) \oplus {\cal O}(b-c) \oplus {\cal O}
\longrightarrow
{\cal O}(a-c) \oplus {\cal O}(b-c)
\longrightarrow 0
$$
so
$$
h^0({\cal F}(-c))= \;
\left\{ \matrix{
3 & \hbox{if} & a=b=c \cr
2 & \hbox{if} & a<b=c\cr
1 & \hbox{if} & \hfill b<c\cr
} \right. 
\quad \hbox{or, equivalently:}\quad
\dim|R_{a,b}|= \;
\left\{ \matrix{
2 & \hbox{if} & a=b=c \cr
1 & \hbox{if} & a<b=c\cr
0 & \hbox{if} & \hfill b<c\cr
} \right. .
$$

\medskip
\noindent
{\bf Remark 9.8.}
Set $\vv:=V_{\ss_0}$ and let as usual $\Sigma$ be the set of the double points of $\xx_0$.
We have the diagram
$$
\matrix{
\ss_0 & \subset & \vv & \supset &\rr \cr
\mapdown {} & & \mapdown {\pi_{\Sigma}} & & \mapdown {}\cr
S & \subset & V & \supset & R \cr
}
$$
where $\rr:=\pi_{\Sigma}^{-1}(R)$.
Setting $\delta_R:= \sharp(\Sigma \cap \rr)$, i.e. the number of
the double points (possibly   infinitely near) of $\xx_0$ lying on $\rr$, it is clear that
$\deg(\rr)=\deg(R)+\delta_R=a+b+\delta_R$.
\medskip

\medskip
\proclaim 
Lemma 9.9. Let $R \in |R_{a,b}|$ be a fixed ruled surface on
$V={\cal O}(a) \oplus {\cal O}(b) \oplus {\cal O}(c)$ and $\ss_0 = S_{2,\lambda-2}$ be
as usual. Then
$$
\wdd :=\rr \cdot \ss_0 \sim 2l +(\lambda-2 -c-\delta + \delta_R)l' 
$$
and there exists a unique bisecant curve 
$\dd \sim 2l +(\lambda-2 -c)l' \subset \ss_0$ such that 
$\Sigma \subset \dd$ and $\dd \supseteq \wdd$.
Moreover, as soon as $R$ varies in
$|R_{a,b}|$,  $\dd$  varies in a linear system of dimension
$0,1,2$ if $b<c$, $a<b=c$, $a=b=c$, respectively.
\par
\pf Let $H_{\vv}$ be a hyperplane section of $\vv$ containing
$\rr$. Since each hyperplane section cannot contain any other
unisecant component out of $\rr$, then  $H_{\vv} \sim \rr + \tau
F_{\vv}$, where $F_{\vv}$ is the generic fibre of $\vv$ and $\tau$ is a non negative
integer. \hb 
Clearly, since
$\deg(H_{\vv})=\deg(\vv)=\deg(V)+\delta=a+b+c+\delta$ and
$\deg(\rr)=a+b+\delta_R$, we obtain that
$$
\rr \sim H_{\vv}-(c+\delta-\delta_R)F_{\vv}.
$$
Taking into account that $H_{\vv} \cdot \ss_0= 2l + (\lambda -2)l'$
and $F_{\vv} \cdot \ss_0= l'$, we obtain:
$$
\rr \cdot \ss_0 \sim
2l + (\lambda -2)l'  -(c+\delta-\delta_R)l'
=2l +(\lambda -2-c-\delta+ \delta_R)l'
$$ 
as required. Note that only $\delta_R$
points of $\Sigma$ lie on $\widetilde D$ and the remaining
$\delta -\delta_R$ lie on $\delta -\delta_R$ fibres (possibly
coincident) of $\ss_0$, say $l'_1, \dots, l'_{\delta -\delta_R}$.
Hence
$$
\Sigma \subset\;
\widetilde D \cup l'_1 \cup \dots \cup l'_{\delta -\delta_R} \sim
2l +(\lambda -2-c)l'
$$
so, setting $\dd:=\widetilde D \cup l'_1 \cup
\dots \cup l'_{\delta -\delta_R}$, we obtain that $\dd$ is linearly equivalent
to $2l+(\lambda-2-c)l'$ and  contains both  $\Sigma$ and
$\wdd$, as required. 
Finally, from the above construction, the  divisor 
$\dd$  is unique, for each $\rr$. 
The last statement follows from 9.7.
\cvd

\medskip
\noindent
Keeping the notation above,  one can immediately compute the
degree of $\dd$:
$$
\deg(\dd)=
\int(2l+(\lambda-2-c)l') \cdot (2l+(\lambda-2)l')=4(\lambda-2)-2c.
\eqno(30)
$$

\medskip
\noindent
Observe that  $\rr$ is the ruled
surface generated by the ruling of $\vv$ on $\widetilde D$, i.e.
$$
\rr = {\bigcup}_{P,Q \in \widetilde D \cap F_{\vv}} \; l_{P,Q}
$$
where $l_{P,Q}$ denotes the line passing through the points $P$
and
$Q$. In particular, $\rr$ is determined by $\widetilde D$;   to
stress this fact, we will write $\rr=\rr(\widetilde D)$.

\goodbreak

\bigskip
\noindent
{\bf 10. Moduli spaces of $4$--gonal curves with $t=0$}

\medskip
In this section we study the moduli spaces of $4$--gonal curves with 
given invariants; in particular we determine whether they are
irreducible and find their dimension. Moreover we give a
stratification of these spaces using the invariants introduced in
the previous sections.

\medskip

Let $X$ be a
$4$--gonal curve of genus $g$ and consider its canonical model
$X_K \subset S \subset V \subset \Bbb P^{g-1}$, where (from 1.1) $S$ is a surface ruled by conics, of minimum degree and unique,  unless $g$ is odd and
$\deg(S)= {3g-7 \over 2}$. In this case, there is  a pencil of such surfaces.

Assume that $S$ has invariant $t=0$, i.e. its (embedded) standard model is the quadric surface $R_{1,1} \subset \Bbb P^3$, on which $X$ can be realized  as a curve $X' \sim 4l +\lambda l'$
having only double points as singularities: we will write $X=X(g,
\lambda)$. Moreover, if $V= \Bbb P({\cal O}(a) \oplus {\cal O}(b) \oplus {\cal O}(c))$, then
$a$ and $b$ are further invariants of $X$ and we will write
$X=X(g, \lambda,a,b)$.

\medskip
\noindent
{\bf Remark 10.1.}  If $X$ is as before,  then by  6.9  it is clear that it has a finite
number  of models $X'$, at most  $\delta \choose 2$, on 
$R_{1,1}$ unless $g$ is odd and $\deg(S)= {3g-7 \over 2}$. In this case, there is  a
one--dimensional family of such models of $X$. 
More precisely, one model comes from another via an elementary transformation of type $elm_{A,B}$, where $A$ and $B$ are two double points of $X'$ as in 6.9. In this way, denoting by
$X''$  another model of $X$ on $R_{1,1}$ and by $\xi$  an elementary transformations as before,
the set
$$
\Xi_{X'} := \{\xi: X' \longrightarrow X'' \}
$$   
consists of at most $\delta \choose 2$
elements if $\deg(S) \le \lceil {{3g-8} \over 2} \rceil$, while $\dim(\Xi_{X'}) = 1$
if $\deg(S)= {3g-7 \over 2}$. \hb
Note that $\Xi_{X'} $ has exactly $\delta \choose 2$ elements in the general case.

\bigskip
Let us denote by ${\cal A}_\lambda$ the open subset of the linear
system $|4l+\lambda l'|$ on $R_{1,1}$ parametrizing the
irreducible curves of such linear system and set
$$
\displaylines{
{\cal W}^{\lambda}_g  := \{X'  \in {\cal A}_\lambda \; | \;
X =X(g,\lambda) \; \hbox{and $X'$ has $\delta$ double points on distinct fibres} \} \cr
{\cal W}^{\lambda}_g(a,b)  := \{X'  \in {\cal W}^{\lambda}_g \; 
| \; X=X(g,\lambda,a,b)\}  . \cr}
$$
Let us denote by ${\cal M}_{g,4}$ the moduli space of $4$--gonal curves of genus $g$ and let
$$
\theta: {\cal W}^{\lambda}_g \longrightarrow
{\cal M}_{g,4}
$$
be the usual projection defined by $\theta(X') = [X]$, where $[X]$ is the isomorphism class of the four--gonal curve $X$ in  ${\cal M}_{g,4}$. Finally set
$$
{\cal M}^{\lambda}_g := \theta({\cal W}^{\lambda}_g), 
\quad
{\cal M}^{\lambda}_g(a,b):=\theta({\cal W}^{\lambda}_g(a,b)). 
$$
It is clear that, in order to compute the dimension of these moduli spaces, we need to find both the dimensions of ${\cal W}^{\lambda}_g$ (resp. ${\cal W}^{\lambda}_g(a,b)$) and of the general fibre of $\theta$.

\medskip
\noindent
{\bf Remark 10.2.} {}From 8.5, the locally closed subsets 
${\cal W}^{\lambda}_g(a,b)$ and hence
${\cal W}^{\lambda}_g$ are not empty, as soon as $a,b,\lambda$
fulfil $(R_1),(R_2),(R_3)$.

\medskip
\proclaim
Lemma 10.3. Let $X',Y' \in {\cal W}^{\lambda}_g$  be two curves
on $R_{1,1}$. If
$[X] =[Y]$ in 
${\cal M}^{\lambda}_g$,  there exists an automorphism $\beta$ of the quadric
surface $R_{1,1}$ and a morphism $\xi \in \Xi_{Y'}$ such that
$$
Y' = \xi(\beta(X')).
$$ 
Therefore the dimension of the general fibre of $\theta$ is:
$$
\dim(\theta^{-1}([X]))=
\cases{
7 & if  $g$ is odd and $\lambda=\left\lceil {{g+2} \over 2}
\right\rceil$ \cr
6 & otherwise \cr}.
$$
\par
\pf Since $X \cong Y$, then  $X_K \cong Y_K$ and there exists a linear automorphism, 
  $\al$ say, of 
$\Bbb P^{g-1}$ such that $\al(X_K) = Y_K$. \hb
Let $S_X$ and $S_Y$ be the surfaces, ruled by conics and of minimum degree such that
 $X_K \subset S_X \subset \Bbb P^{g-1}$ and  $Y_K \subset S_Y \subset \Bbb P^{g-1}$. 
 Assume that these surfaces are unique: therefore $\al(S_X) = S_Y$. \hb
 Let us consider the diagram $(8)$ for both $X$ and $Y$: defining with obvious notation
 $N_X:= \langle \varphi_X(K_X - \Phi_X - \Lambda_X) \rangle$ and $N_Y$ analogously, we  have
 $$
 \matrix{
  \Bbb P^{g-1} & \supset & S_X & \supset & X_K & \mapright \pi_{N_X},40 & 
  X_{\Phi_X +\Lambda_X}= X' & \subset & R_{1,1}(X)\cr
  \lmapdown {\al} &&  \lmapdown {\al} &&  \lmapdown {\al} &&&&  \lmapdown {\beta} \cr
  \cr
   \Bbb P^{g-1} & \supset & S_Y & \supset & Y_K &  \mapright \pi_{N_Y},40 & 
   Y_{\Phi_Y +\Lambda_Y}= Y' & \subset & R_{1,1}(Y)\cr
  }
  $$
where $\beta$ is the isomorphism between the quadrics $R_{1,1}(X)$ and $R_{1,1}(Y)$ induced by $\al$. Up to a linear change of coordinates in $\Bbb P^3$, we can assume that $R_{1,1}(X) =R_{1,1}(Y)$ so $\beta\in {\rm Aut}(R_{1,1})$. \hb
Consider then the curves $Y'$ and $\beta(X')$ lying on $R_{1,1}$: from the construction above, we obtain that they are both models of $Y$ on a quadric. Therefore, applying 10.1, we get that there exists $\xi \in \Xi_{Y'}$ such that
$Y' = \xi(\beta(X'))$, as requested. \hb
  When $S_X$ and $S_Y$ are not unique they vary in a pencil (see 1.1) and the proof runs in a similar way. \hb
The second part of the statement follows from the first part; namely, it is clear that
$$
\dim(\theta^{-1}([X]))= \dim ({\rm Aut}(R_{1,1})) + \dim (\Xi_X).
$$
On one hand, observe that
${\rm Aut}(R_{1,1}) \cong 
{\rm Aut}(\Bbb P^1 \times \Bbb P^1) \cong PGL(2) \times PGL(2)$
has dimension 6. \hb
On the other hand, by 10.1,
$$
 \dim (\Xi_X) = 
 \cases{
 1 & if $g$ is odd and $\deg(S) =  {{3g-7} \over 2}$ \cr
 0 & otherwise\cr}.
 $$
Finally note that  (using 4.4):
$$
g+\lambda-5 = \deg(S)= {{3g-7} \over 2}
$$
or, equivalently
$$
\lambda= {{g+3} \over 2} = 
\left\lceil {{g+2} \over 2} \right\rceil
$$
where the last equality holds since $g$ is odd.
\cvd 

\medskip
\noindent
Let us recall (see Section 8) that, if  $X' \in {\cal W}^\lambda_g$ then $X' \subset R_{1,1} \cong \hs0$ and $\varphi_{4l+\lambda l'} : \hs0 \longrightarrow S' \subset \Bbb P^{5\lambda+4}$; in particular, we can associate  to $X'$ a hyperplane $H_X$ of $\Bbb P^{5\lambda+4}$. 
By 8.2 we have  that $X' $ has
$P_1,
\dots,P_\delta$ as double points if and only if $H_X$ contains the linear space
$$
L_{P_1,\dots,P_\delta}:= 
\langle 
T_{P_1}(S'), \dots, T_{P_\delta}(S') 
\rangle.
$$
In this way we can identify ${\cal W}^\lambda_g$ with its
image via the injective morphism
$$
\eqalign{
i: \quad {\cal W}^\lambda_g  &\longrightarrow
\check{\Bbb P}^{5\lambda+4} \cr
X' & \mapsto H_X \cr}.
$$

\medskip
In order to compute the dimension of ${\cal W}^{\lambda}_g$ and of ${\cal W}^{\lambda}_g(a,b)$ and to prove their irreducibility, we need further preliminary observations.

 \medskip
\noindent
{\bf Remark 10.4.} The Key--Lemma 9.4 has been proved under the assumption 
that $(P_1, \dots, P_\delta)$  are distinct points. 
For instance, if $\delta =2$, this result says that
$$
\dim L_{P_1,P_2} = \dim \langle T_{P_1}(S'), T_{P_2}(S') \rangle = 5.
$$
 If
$P_2$ is infinitely near  to
$P_1$, given a local system of coordinates of $S'$ in a neighbourhood of
$P_1$,  the  tangent plane to
$S'$ at
$P_1$ is generated by $P_1$ and the first derived vectors both along the bisecant $\wdd$ and
along the fibre $l'_1$. Hence  it is easy to see that the linear space 
$L_{P_1,P_2}$ is generated by the above generators of $T_{P_1}(S')$ and by two
further  second derived vectors and a third derived vector. 
One can show that all of them are
linearly independent so, also in this case, $\dim L_{P_1,P_2} = 5$. \hb
It is not difficult to prove that, if $k$ is any integer ($1 \le k \le \delta -1$) and
the  considered points are
$P_1, P_2, \dots, P_{k+1}, \dots, P_\delta$ where  
$P_2, \dots, P_{k+1}$ are infinitely near to $P_1$, 
then
$$
\dim L_{P_1, \dots P_\delta} \ge 3 \delta - k.
$$

\goodbreak

\medskip
\proclaim
Lemma 10.5.  Let us consider the morphism 
$$
\eqalign{
\Psi: \; {\cal W}^{\lambda}_g  & \longrightarrow
\sym^{\delta}(R_{1,1}) \cr
X' & \mapsto (P_1, \dots,P_{\delta}) \cr}
$$
where $\Sigma = P_1 + \cdots + P_{\delta}$ is the singular locus of $X' \subset R_{1,1}$. Then the general fibre of $\Psi$ has dimension
\item{$i)$} $\dim (\Psi^{-1}(P_1, \dots,P_{\delta})) = 5\lambda +4 - 3 \delta$ if $P_1, \dots,P_{\delta}$ are distinct points;
\item{$ii)$} $\dim (\Psi^{-1}(P_1, \dots,P_{\delta})) \le  5\lambda +3 - 3 \delta +k$ if $P_2, \dots,P_{k+1}$ are infinitely near to $P_1$, for some $k \ge 1$.
\par
\pf
By definition, ${\cal W}^{\lambda}_g$ consists of the irreducible curves of type $(4, \lambda)$ on $R_{1,1}$ having $\delta$ double points on distinct fibres. So, taking into account the above injective morphism 
$i: \; {\cal W}^\lambda_g  \longrightarrow \check{\Bbb P}^{5\lambda+4}$ and the fact that $X'  \in  {\cal W}^\lambda_g $ has $P_1, \dots,P_\delta$ as double points if and only if the hyperplane $H_X:= i(X')$ contains the linear space
$L_{P_1,\dots,P_\delta}$, it is clear that the general fibre $\Psi^{-1}(P_1, \dots,P_{\delta})$ is isomorphic to
 an open subset  of
$\{H \in \check {\Bbb P}^{5\lambda +4} \;| \; 
H \supset L_{P_1, \dots,P_{\delta}}\}$, since the general hyperplane containing $ L_{P_1, \dots,P_{\delta}}$ contains the tangent planes to $S'$ only at the choosen points. This means exactly that
$$
\dim (\Psi^{-1}(P_1, \dots,P_{\delta})) = 5\lambda +4- (\dim L_{P_1, \dots P_\delta} +1).
$$
\item{$i)$} If $P_1, \dots,P_{\delta}$ are distinct, then in the Key--Lemma 9.4 we have shown that the
dimension of $L_{P_1, \dots,P_{\delta}}$ is $3 \delta -1$ independently on the position
of the considered points. So, in this case, $\Psi^{-1}(P_1, \dots,P_{\delta})$ is irreducible of dimension $5\lambda +4-3\delta$.
\item{$ii)$} If $P_1, \dots,P_{\delta}$ are not distinct -- as in the assumption -- then the  fibre  of $\Psi$ could have bigger dimension. Nevertheless, we can get an upper bound  on  this dimension by taking into account 10.4, obtaining that $\dim (\Psi^{-1}(P_1, \dots,P_{\delta}))$ is at most $5\lambda +4-(3\delta-k +1)$ and this proves the second part of the statement. 

\cvd

\medskip
\proclaim
Proposition 10.6. For each $\lambda$
satisfying 
$$
{{g+3} \over 3} \le \lambda \le 
\left\lceil{{g+2} \over 2}\right\rceil
\eqno{(R_1)}
$$ 
the locally closed subset 
${\cal W}^{\lambda}_{g}$ is irreducible of dimension $g+2\lambda+7$.
\par
\pf 
Setting $\sym:= \sym^{\delta}(R_{1,1})$, consider the map $\Psi: \; {\cal W}^{\lambda}_g   \rightarrow \sym$ defined in 10.5. Note that $\Psi$ is dominant and $\dim(\sym) = 2 \delta$. \hb
Recall also that the $\delta$ singular points of the general curve $X' \in  {\cal W}^{\lambda}_g$ are in general position  on $R_{1,1}$ by 9.4. \hb
If $P_1, \dots,P_{\delta}$ are distinct points, by 10.5 we get that $\dim (\Psi^{-1}(P_1, \dots,P_{\delta})) = 5\lambda +4 - 3 \delta$. Therefore
$$
\eqalign{
\dim({\cal W}^{\lambda}_g) &=
\dim (\Psi^{-1}(P_1, \dots,P_{\delta})) + \dim(\sym)= \cr
{} &= 5\lambda+4- \delta =  \cr
{} & = g + 2\lambda + 7\cr}
$$
where the last equality follows from $\delta = 3(\lambda -1) -g$. \hb
Assume now that $P_2, \dots, P_{k+1}$ are infinitely near to $P_1$ for some $k \ge 1$. Then 
the fibre of $\Psi$ at the point $(P_1, \dots,P_{\delta}) \in \sym$ has dimension at most $5 \lambda +3 - 3\delta +k$ by 10.5. The difference between such integer and $5\lambda +4-3\delta$ is at most 
$$
k-1 < 2k = \codim_{\sym}(\Delta)
$$
where  $\Delta: = \{ (Q_1, \dots, Q_{\delta}) \in \sym \; | \; Q_1 = \cdots = Q_{k+1} \}$. Clearly $\Delta$  is a closed subset of
$\sym$ and contains  the considered element $(P_1, \dots,P_{\delta})$.  
Therefore, the variety consisting of the fibres on the points of $\Delta$ is a proper closed subset of $ {\cal W}^{\lambda}_g$. 
\cvd

\medskip
\noindent
{\bf Remark 10.7.} Recall that ${\cal M}_{g,4}$ is a closed
irreducible subset of the moduli space ${\cal M}_{g}$ and has
dimension $2g+3$. Let us set the maximum value of
$\lambda$ (see
$(R_1)$):
$$
\lambda_{\max}:=
\left\lceil{{g+2} \over 2}\right\rceil.
$$
Then, from 10.6
$$
\dim({\cal W}_{g}^{\lambda_{\max}}) = g+2\lambda_{\max}+7.
$$
Let us recall that the fibre of 
$\theta: {\cal W}_{g}^{\lambda_{\max}} \rightarrow
{\cal M}_{g}^{\lambda_{\max}}$ has dimension either 6 or 7,
accordingly   to wheter  $g$ is even or odd, respectively (from 10.3). Hence
$$
\dim({\cal M}_{g}^{\lambda_{\max}}) =
\cases{
g+2{{g+2} \over 2}+1=2g+3, & if $g$ is even ;\cr
g+2{{g+3} \over 2}=2g+3, & if $g$ is odd. \cr}
$$
Therefore, in both cases, we have that 
$\dim({\cal M}_{g}^{\lambda_{\max}}) =
\dim({\cal M}_{g,4})$; in other words, the general $4$--gonal curve
has invariant $\lambda_{\max}$.

\medskip
\noindent
{\bf Remark 10.8.} We know that, if $t>0$, then $X$ admits a standard model $X'
\subset R_{1,t+1}$. Nevertheless, also in this case, it is possible to define another model of $X$,  
$X''$ say,  on a quadric surface $R_{1,1}$. Clearly, in this situation, $X''$ will have not only
double points as singularities, but also triple points. \hb
Namely, let $Q_1, \dots, Q_t$ be simple points of $X'$, belonging to $t$ distinct
fibres of $R_{1,t+1}$ and consider the projection from these points:
$$
\matrix{
X' & \subset & R_{1,t+1} \cr
\mapdown {} && \mapdown {\pi_{Q_1, \dots, Q_t}} \cr
\cr
X'' & \subset & R_{1,1} \cr
}
$$
Since $X'$ meets each fibre of $R_{1,t+1}$ in the four points of the gonal
divisor,  the singularities of $X''$ are the $\delta$ double points of $X'$ and, in
addition, $t$ triple points, all of them belonging to the same line $l$. \hb
It is clear that the closure $\overline {{\cal W}^\lambda _g}$ of 
${\cal W}^\lambda _g$ in ${\cal A}_\lambda$ contains also the curves of invariants $g,
\lambda$ and $t>0$ and it is not difficult to see that the closed subset consisting of
such curves has dimension smaller then $\dim({\cal W}^\lambda _g)$.

\medskip
Using 10.2, 10.3, 10.6, 10.7 and 10.8, we immediately obtain
the following result,  which is the first part of the Main Theorem stated in the Introduction
(here $\overline{\cal M}^{\lambda}_{g}$ denotes the closure of
${\cal M}^{\lambda}_{g}$ in the moduli space
${\cal M}_{g,4}$ of $4$--gonal curves):

\medskip
\proclaim
Theorem 10.9. There exists a stratification of the moduli space
${\cal M}_{g,4}$ of $4$--gonal curves given by:
$$
{\cal M}_{g,4}=
\overline{\cal M}^{\left\lceil{{g+2} \over 2}\right\rceil}_g
\supset 
\overline{\cal M}^{\left\lceil{{g} \over 2}\right\rceil}_g
\supset \cdots \supset
\overline{\cal M}^{\lambda}_{g}
\supset \cdots \supset
\overline{\cal M}^{\left\lceil{{g+3} \over 3}\right\rceil}_g
$$
and $\overline{\cal M}^{\lambda}_g$ are irreducible locally closed
subsets of dimension $g+2\lambda+1$, for each $\lambda$
satisfying ${{g+3} \over 3} \le \lambda < 
\left\lceil{{g+2} \over 2}\right\rceil$.
\cvd
\par

\medskip
In order to show the second part of the Main Theorem, let us
start with some preliminary fact. 

\medskip
\noindent
We keep the notation of 9.9,
where $\widetilde D$ denotes a divisor of $\ss_0 = S_{2,\lambda-2} \subset \Bbb
P^{g-1+\delta}$ linearly equivalent to
$2l +(\lambda-2-c-\delta+\delta_R)l'$ and containing $\delta_R$
points among $P_1, \dots, P_\delta$.  \hb
Recall also that, referring to Section 7, the unisecant $\aa \subset \vv$ is the preimage, via $\pi$, of the (unique if $a<b$) unisecant of degree $a$ of $V$. Moreover, $\rr := \pi^{-1}(R)$, where $R:=R_{a,b}$, so $\aa \subset \rr = \rr(\wdd)$ as described in 9.9. \hb
In the forthcoming computations we will use a few times the following relations (coming from $a+b+c = g-3$ and from $(17)$):
$$
c =g-3-a-b, \qquad 3\lambda = \delta +g+3  . \eqno{(31)} 
$$

\medskip
\proclaim
Lemma 10.10. Let $\wdd \subset \ss_0$ and
$\rr:=\rr(\widetilde D)$ be as before.  
Let $\aa \in Un^{a+\delta_R}(\rr)$ and
$\Gamma:= \widetilde D \cdot \aa$. Assume that $a \ge (g-\lambda
-1)/2$. Then:
\item {$i)$} $\deg(\Gamma)  = 4(\lambda-2)
-2b-2c-2(\delta-\delta_R)$;
\item {$ii)$} $h^0({\cal O}_{\rr}(\aa)) =
h^0({\cal O}_{\widetilde D}(\Gamma))$;
\item {$iii)$} assume also that $\delta_R=\delta$ and  
either $a > (g-\lambda -1)/2$ or $a = (g-\lambda -1)/2$ and $a<b$;
then:
$$
H^0({\cal O}_{\rr}(\aa)) \cong
H^0({\cal O}_{\widetilde D}(\Gamma)).
$$
\par
\pf $i)$ Recall that, keeping the notation in 9.9,  $\dd = \wdd + (\delta - \delta_R)
l'$. So
$\deg(\wdd )=
\deg(\dd)-2(\delta-\delta_R)$ since $\ss_0$ is ruled by conics. Hence, using $(30)$,
we obtain that $\deg(\wdd )=
 4(\lambda-2)-2(c+\delta-\delta_R)$.
Therefore, applying $(IF)$ and 9.8, we have that    
$$
\eqalign{
\deg(\Gamma) & =
2\deg(\aa)+\deg(\widetilde D)-2\deg(\rr)= \cr
{} & = 2(a+\delta_R) + 4(\lambda-2)-2(c+\delta-\delta_R) -
2(a+b+\delta_R)= \cr {} & = 4(\lambda-2) -2b-2c-2(\delta-\delta_R).
\cr}
$$
$ii)$ 
Let us show first that $\Gamma$ is a non special divisor on $\wdd$. Since $\widetilde D$
is of type $(2,\lambda -2-c-(\delta-\delta_R))$ on the quadric, then
$p_a(\widetilde D)=\lambda-3-c-(\delta-\delta_R)$. A sufficient condition in
order to have $\Gamma$ non special is
$\deg(\Gamma) > 2p_a(\widetilde D)-2$, or, equivalently:
$$
4(\lambda -2)-2b-2c-2(\delta-\delta_R) > 2(\lambda
-3-c-(\delta-\delta_R)) -2
\quad \Longleftrightarrow \quad
\lambda -b> 0
$$
and this is true since $b \le \lambda -2$.
Therefore
$h^1({\cal O}_{\wdd}(\Gamma))=0$ and, by Riemann--Roch
theorem, using also $(31)$, we obtain that
$$
h^0({\cal O}_{\widetilde D}(\Gamma))-1= 
\deg(\Gamma)  -p_a(\widetilde D)= a-b+\delta_R+1.
$$
Moreover $h^0({\cal O}_{\rr}(\aa))-1=
\dim_{\rr}(|\aa|)=\dim({\rm Un}^{a+ \delta_R}(\rr))=
a-b + \delta_R +1$ by $(UF)$. Hence we obtain that
$$
h^0({\cal O}_{\widetilde D}(\Gamma))= a-b+\delta_R+2=
h^0({\cal O}_{\rr}(\aa)).
$$
$iii)$ In order to prove the claim, consider the exact sequence
$$
0 \longrightarrow
{\cal I}_{\widetilde D/\rr}(\aa) \longrightarrow
{\cal O}_{\rr}(\aa) \longrightarrow
{\cal O}_{\widetilde D}(\Gamma) \longrightarrow 0.
\eqno{(32)}
$$
By $ii)$, it suffices to show that  the map
$f: \; H^0({\cal O}_{\rr}(\aa)) \rightarrow
H^0({\cal O}_{\widetilde D}(\Gamma))$ 
induced by $(32)$ is injective. \hb
Clearly this holds if and
only if there exists a unique $\aa \in Un^{a+\delta_R}(\rr)$
passing through
$\Gamma$ and this holds if $\int {\aa\ }^2 < \deg(\Gamma)$. {}From
$(IF)$ and  9.8 we obtain that
$$
\int {\aa\ }^2 = 2\deg(\aa)-\deg(\rr) =
2(a+\delta_R)-(a+b+\delta_R)=a-b+\delta_R.
$$
Therefore the  condition  $\int {\aa\ }^2 < \deg(\Gamma)$ becomes
$$
a-b+\delta_R < 4(\lambda -2) -2b-2c-2(\delta - \delta_R).
$$
Using again $(31)$, the above inequality is equivalent to:
$$
\lambda -g +a+b+1 -(\delta -\delta_R) >0.
$$
By assumption $\delta -\delta_R= 0$, so 
$$
a+b >g-\lambda -1
$$
and using the further assumptions on $a$ and $b$, the claim is proved.
\cvd

\bigskip
Before stating the second part of the Main Theorem, let us set
$$
\epsilon := \cases{
0, & if $b<c$ \cr
1 , & if $a<b=c$ \cr
2 , & if $a=b=c$ \cr}
\quad , \quad
\tau := \cases{
0, & if $a<b$ \cr
1 , &  if $a=b$ \cr}
\quad \hbox{and} \quad
\xi := \cases{
1, & if $\lambda={{g+3} \over 2}$ \cr
0 , & otherwise \cr}.
$$

\medskip
\proclaim
Theorem 10.11. Let $g,\lambda,a,b$ be positive integers satisfying
$(R_1)$, $(R_2)$, $(R_3)$ and $c=g-3-a-b$.
If $a \ge (g-\lambda -1)/2$ then 
${\cal M}^{\lambda}_g(a,b)$ is an irreducible variety of dimension
$2(2a+b+\lambda)+10-g-\epsilon-\tau-\xi$.
\par
\pf {}From 10.2 and 10.3, it is enough to show that  ${\cal W}^{\lambda}_g(a,b)$ is
irreducible of the right dimension. \hb
Keeping the notation in 10.5, set $Y^{\lambda}_g(a,b):=\Psi({\cal W}^{\lambda}_g(a,b))$. \hb
Claim:   $\Psi^{-1}(Y^{\lambda}_g(a,b)) \subset {\cal W}^{\lambda}_g(a,b)$. \hb
 This is  equivalent to the following property: let $X'' \in
{\cal W}_g^\lambda$ be such that 
$\Psi(X'')=(P_1,\dots,P_\delta)=\Psi(X')$, where 
$X' \in {\cal W}_g^\lambda(a,b)$; then 
$X'' \in {\cal W}_g^\lambda(a,b)$. 
This is true, since $\pi_{\langle P_1,\dots,P_\delta \rangle}(\vv)$
is the scroll 
$V=\Bbb{P}({\cal O}(a)\oplus {\cal O}(b)\oplus {\cal O}(c))$
associated both to $X'$ and to $X''$ and this proves the claim.

\smallskip
\noindent
{\bf Step 1.} Irreducibility and dimension of ${\cal W}_g^\lambda(a,b)$. \hb
{}From the claim above we can  consider  the
restriction of $\Psi$
$$
\psi: \; {\cal W}^{\lambda}_g(a,b)   \longrightarrow
Y^{\lambda}_g(a,b).
$$
From 10.5, $\dim (\psi^{-1}(P_1, \dots,P_{\delta})) = 5\lambda +4 - 3 \delta$ if $P_1, \dots,P_{\delta}$ are distinct points. \hb
With the same argument as the one in the proof of 10.6, in the case of infinitely near points one easily shows that  the variety consisting of the fibres on the points of $\Delta$ is a proper closed subset of $ {\cal W}^{\lambda}_g(a,b)$. For this reason, ${\cal W}^{\lambda}_g(a,b)$ is irreducible if
$Y^{\lambda}_g(a,b)$ is irreducible and
$$
\dim({\cal W}^{\lambda}_g(a,b)) = \dim(Y^{\lambda}_g(a,b))
+5\lambda +4-3\delta.
\eqno(33)
$$

\smallskip
\noindent
{\bf Step 2.} Irreducibility and dimension of $Y_g^\lambda(a,b)$ . \hb
Recall that  the singular locus 
$\Sigma = P_1 + \cdots + P_{\delta}$ of $X' \subset R_{1,1}$ is contained
in a suitable bisecant curve $\dd \sim 2l+(\lambda-2-c)l' \subset R_{1,1}$ by 9.9 (there the result concerns $\ss_0$, here $R_{1,1}$). \hb
It is not hard to show that there exists an open subset,  
$Y^0$ say, of $Y^{\lambda}_g(a,b)$ whose elements $(P_1,
\dots,P_{\delta})$ fulfil the following property:  there exists
$\dd \in |2l+(\lambda-2-c)l'|$ not containing fibres and such that 
$P_1, \dots,P_{\delta} \in
\dd \cap \aa$, for a suitable $\aa \in 
{\rm Un}^{a+ \delta}(\rr)$, where $\rr := \rr(\dd)$. In
particular, on this subset $\delta_R=\delta$.\hb 
Let us check that the above condition is compatible with the degrees of the involved divisors i.e., setting $\Gamma := \dd \cap \aa$, we must have that $\delta \le \deg(\Gamma)$. 
{}From 10.10 (ii), taking into account that here $\delta = \delta_R$ and using $(31)$ as usual, it is easy to see that $\deg(\Gamma) = 2a+\lambda-g+1+\delta \ge
\delta$, since $2a+\lambda-g+1\ge0$: namely this is equivalent to $a
\ge (g-\lambda-1)/2$, which holds by assumption. \hb
Consider then the following correspondence: 
$$
Z_{a,b}^{\lambda}
\quad \subset \quad
|2l+(\lambda-2-c)l'| \times \sym^{\delta}(R_{1,1})
$$
defined by:
$$
Z_{a,b}^{\lambda}:=\left\{(\dd,P_1 , \dots , P_{\delta}) 
\;  | \; 
\hbox{there exists}\; \aa \in {\rm Un}^{a+ \delta}(\rr(\dd)) \;
\hbox{such that}\; 
P_1 ,\dots  , P_{\delta} \in \dd \cap \aa 
\right\}.
$$
Consider now the two canonical projections, where $\Omega$ is the
open subset of $|2l+(\lambda-2-c)l'|$ consisting of curves not
containing fibres:
$$
\matrix{
& Z_{a,b}^{\lambda} &
\cr
 &{}_p \swarrow{}\searrow {}_q &\cr
|2l+(\lambda-2-c)l'| \supset \Omega &{}&
Y^0 \subset Y^{\lambda}_g(a,b) \subset \sym^{\delta}(R_{1,1})
\hfil  
\cr }.
$$
By 9.9, every element $(P_1, \dots, P_{\delta})$ of $Y^0$
 determines either a unique
$\dd \sim 2l+(\lambda-2-c)l'$ (if $b<c$) or a pencil 
(if $a<b=c$) or a two--dimensional linear system (if $a=b=c$) of such curves. This
implies that the general fibre of $q$ is irreducible of dimension
$\epsilon$, where $\epsilon=0,1,2$ as soon as  $b<c$, $a<b=c$,
$a=b=c$, respectively.  
Furthermore $p$ is surjective by 9.6. \hb
Denoting by
$Z_{\dd}:=p^{-1}(\dd)$ any fibre of $p$,  we have that: if
$Z_{\dd}$ is irreducible, then 
${Y}^{\lambda}_g(a,b)$ is irreducible and
$$
\eqalign{
\dim(Y^{\lambda}_g(a,b)) & = \dim(Z^{\lambda}_{a,b})-\epsilon = 
\dim(Z_{\dd})+ \dim(|\dd|)   -\epsilon \cr
& = \dim(Z_{\dd})+ 3(\lambda-1-c) -1   -\epsilon. \cr}
\eqno(34)
$$

\smallskip
\noindent
{\bf Step 3.} Irreducibility and dimension of $Z_{\dd}$ \hb
It is clear that
$$
Z_{\dd} \cong \{(P_1 , \dots , P_{\delta})\in \sym^{\delta}(\dd)
 \; | \; 
\hbox{there exists}\; \aa \in Un^{a+\delta}(\rr) \;
\hbox{such that} \; P_1 ,\dots , P_{\delta}
\in\dd \cap \aa \}.
$$
 In order to compute
the dimension and to prove the irreducibility of $Z_{\dd}$, consider
the following correspondence (where $\Gamma = \dd \cap \aa$ is as before):
$$
T_{\dd} := \{(P'_1 , \dots , P'_{\delta},\aa)
 \; | \; P'_1 , \dots , P'_{\delta} \in \Gamma \} \subset \sym^{\delta}(\dd) \times Un^{a+\delta}(\rr)
$$
and the two projections:
$$
\matrix{
  {} & T_{\dd} &  {}  \cr
 {} & {}_{\pi_1} \swarrow{}\searrow {}_{\pi_2}   &  {}  \cr
\sym^{\delta}(\dd)  & {} & Un^{a+\delta}(\rr)  \cr}
$$
Obviously, $Im(\pi_1) =Z_{\dd}$ and $\pi_2$ is a finite
surjective morphism; hence, denoting by $\tau$ the dimension of the
fibres of $\pi_1$, we obtain:
$$
\dim(Z_{\dd})=\dim(T_{\dd})- \tau=
\dim(Un^{a+\delta}(\rr))-\tau= a-b+\delta+1-\tau. \eqno{(35)}
$$
Let us find the possible values of $\tau$.

\goodbreak

\noindent
In the proof of 10.10 (iii) we show that  $\int {\aa\ }^2 = a-b+\delta$; with the same argument used there to prove the uniqueness of the unisecant $\aa$ passing through a certain divisor, it is immediate to see that
$$
\tau=0 \quad \Leftrightarrow \quad
\int {\aa\ }^2 <  \delta
\quad \Leftrightarrow \quad
a-b+ \delta  <  \delta
\quad \Leftrightarrow \quad
a<b.
$$
With the same argument we obtain:
$$
\tau \ge 1 \quad \Leftrightarrow \quad
\int {\aa\ }^2   \ge  \delta
\quad \Leftrightarrow \quad
a-b+ \delta  \ge  \delta
\quad \Leftrightarrow \quad
a=b \; \hbox{and} \; \int {\aa\ }^2   =  \delta.
$$
Hence, necessarily, $\tau = 1$ and $a=b$. \hb
We are left to show that $Z_{\dd}$ is irreducible.  Since $Z_{\dd} = \pi_1(T_{\dd})$, it is enough to show that $T_{\dd}$ itself is irreducible. \hb
Assume first that 
$$
a > {{g-\lambda -1} \over 2}
\quad \hbox{or} \quad
a = {{g-\lambda -1} \over 2} <b.
$$
It follows from 10.10 (iii) that
$H^0({\cal O}_{\rr}(\aa)) \cong
H^0({\cal O}_{\dd}(\Gamma))$, hence
$$
T_{\dd} \cong
\{(P'_1 , \dots , P'_{\delta},\Gamma')
 \; | \; P'_1 , \dots , P'_{\delta} \in \Gamma'\}
\subset \sym^{\delta}(\dd) \times |\Gamma|.
$$
Consider the morphism associated to $|\Gamma|$:
$$
\varphi_{\Gamma}: \quad
\dd \longrightarrow \Bbb P^r
$$
where $r=\dim|\Gamma|=a-b+\delta+1$ (as computed in the proof of 10.10 (ii)); if $\dd'$ denotes the image of $\dd$ in $\Bbb P^r$, it is clear that
$$
T_{\dd} \cong \{(P'_1 , \dots , P'_{\delta},H)\; | \; P'_1 , \dots ,
P'_{\delta} \in H \cap \dd'\}
\subset  
\sym^{\delta}(\dd')  \times \check{\Bbb P}^r.
$$
The irreducibility of $T_{\dd}$ is a consequence of the forthcoming
lemma 10.12.  \hb
Finally, we have to consider the last case: 
$$
a = {{g-\lambda -1} \over 2} =b.
$$
Since $c=g-3-(a+b)=\lambda-2$,  from 10.10 (i) we have
$$
\deg(\Gamma)=4(\lambda-2)-2b-2c=3\lambda-3-g=\delta.
$$
Therefore $\pi_2: T_{\dd} \rightarrow Un^{a+\delta}(\rr)$ is an
isomorphism, hence $T_{\dd}$ is irreducible of dimension $\delta+1$
(since $a=b$). \hb
 Finally observe that, if $\ddd \not \in \Omega$ in Step 2, 
then one can easily prove that 
$\dim(Z_{\ddd})= a-b+\delta_R+1-\tau$. In
particular, $\dim(Z_{\ddd})<\dim(Z_{\dd})$
hence $p^{-1}\left(|2l+(\lambda-2-c)l'| \setminus \Omega\right)$ is a Zariski locally
 closed subset of $Z^{\lambda}_{a,b}$.
 
 \goodbreak
 
 \smallskip
 \noindent
 {\bf Step 4.} Final computation \hb
We can now compute the dimension of the moduli space
using $(33)$, $(34)$,  $(35)$ and $(31)$:
$$
\eqalign{
\dim({\cal W}^{\lambda}_g(a,b)) &= \dim(Y^{\lambda}_g(a,b))
+5\lambda +4-3\delta= \cr
& = \dim(Z_{\dd})+ 3(\lambda-1-c) +3-\epsilon +5\lambda -3\delta=
\cr
& = 2(2a+b+\lambda)+16-g-\epsilon -\tau \cr}
$$
hence, from 10.3, we obtain
$$
\dim({\cal M}^{\lambda}_g(a,b)) =  \dim({\cal W}^{\lambda}_g(a,b))-6-\xi=
2(2a+b+\lambda)+10-g-\epsilon -\tau-\xi
$$
and this proves the claim. 
\cvd

\medskip
We are left to show the following fact:

\medskip
\proclaim
Lemma 10.12. Let $X \subset \Bbb P^r$ be a (smooth) irreducible
curve, $k$ an integer such that $k \le \deg(X)$ and let 
$$
V_X:=\{(P_1, \dots, P_k;H) \;|\; P_1, \dots , P_k \in H \cap X \}
\subset \sym^k(X) \times \check{\Bbb P}^r.
$$
Then the variety $V_X$ is irreducible.
\par
\pf
It is a straightforward generalization of the argument used in the proof of the
Uniform Position Lemma,  [{\bf 9}]. 
\cvd

\medskip
Now we are going to prove the last part of the Main Theorem. We
need first some preliminary results; let us recall that, if 
$a < (g-\lambda-1)/2$, then $\aa \subset \ss_0 \subset \vv$ (from 7.3).

\medskip
\proclaim
Lemma 10.13. 
Let $a < (g-\lambda-1)/2$ and 
$[X] \in {\cal M}^{\lambda}_g(a,b)$. Then in $\theta^{-1}([X])$ there
exists a curve $X'  \subset  R_{1,1}$ such that
$\aa \sim l$. In particular, $\deg(\aa) =\lambda-2$ and
$\delta_A=\lambda-2-a$.
\par
\pf Let $\aa \sim l+ \al l' \subset \ss_0= 
\varphi_{2l+(\lambda -2)l'}(\hs 0) \subset \Bbb P^{3 \lambda -4}$ and assume 
$\al \ge 1$. Since 
$$
\deg_{\ss_0}(\aa) = \int (l+ \al l') \cdot (2l+(\lambda -2)l')=
\lambda - 2 +2 \al
\eqno(36)
$$
and $\deg(A)=a \le \lambda-2$ (from 7.1), 
 the number of double points of $\xx_0$ lying on $\aa$
is, from $(11)$, 
$\delta_A = \deg(\aa) - \deg(A) = \lambda - 2 +2\al -a \ge 2 \al$. Therefore, since $\aa$ meets each line
of the ruling $l$ of $\ss_0$  in $\al$ points, there are at least two double
points of $\xx_0$,  $N_1$ and $N_2$ say,  belonging to $\aa$ and not
belonging to a same line $l$. \hb
Consider now the isomorphism
$$
\varphi_{l+2l'}: \quad
R_{1,1} \cong \ss_0 \longrightarrow \tilde S \cong R_{2,2}
$$
and set $\tilde A:=\varphi(\aa) \sim \tilde l + \al \tilde l'$; for simplicity, we still denote by $N_1$ and $N_2$ the images of these points in $\wss$.
Clearly $\deg(\tilde A)= \al +2$ and the projection
$$
\pi_{\langle N_1,N_2 \rangle}: \quad
\tilde S \longrightarrow R_{1,1}
$$
maps $\tilde A$ to a unisecant curve $\aa^*$ of degree $\al$ (since
$N_1,N_2 \in \tilde A$) lying on $R_{1,1}$; hence 
$\aa^* \sim l + (\al-1)l'$; in particular, from $(36)$,
$\deg_{\ss_0}(\aa^*)=\lambda -2 +2(\al-1)$.  \hb 
Set $X':=
\left(\pi_{\langle N_1,N_2 \rangle} \circ
\varphi_{l+2l'}\right) (X) \subset R_{1,1}$ and $A^* \subset S$ be the curve
corresponding to $\aa^* \subset R_{1,1}$. Since the number of the double points of
$X'$ lying on $\aa^*$ is $\delta_A-2$, we get that
$$
\deg(A^*)= \deg_{\ss_0}(\aa^*) - (\delta_A-2)=
\lambda -2 +2\al -\delta_A = a =
\deg(A)
$$
and this implies that $A^*=A$.
Iterating this procedure we obtain a model of $X$ such that $\al =0$, hence $\aa \sim l$ and the other requirements are fulfilled.
\cvd

\medskip
\proclaim
Corollary 10.14. Let $a < (g-\lambda -1)/2$ and let
$\widetilde{\cal W}^{\lambda}_g(a,b) \subset
{\cal W}^{\lambda}_g(a,b)$ be the following set:
$$
\widetilde{\cal W}^{\lambda}_g(a,b) := \{X' \in {\cal W}^{\lambda}_g(a,b) \;|\; X'
\subset R_{1,1}, \aa
\sim l\}.
$$
Then the restriction
$$
\theta: \;\widetilde{\cal W}^{\lambda}_g(a,b)  \longrightarrow
{\cal M}^{\lambda}_g(a,b) 
$$
is surjective and the fibres have dimension 6 unless 
$g$ is odd and $\lambda = (g+3)/2$: in this case they have dimension 7.
\par
\pf The surjectivity is immediate by 10.13 and the dimension of
the fibres can be computed with the same argument of 10.3.
 \cvd

\medskip
Let us set 
$$
\epsilon := \cases{
0, & if $b<c$ \cr
1 , & if $a<b=c$ \cr}
\quad \hbox{and} \quad
\xi := \cases{
1, & if $\lambda={{g+3} \over 2}$ \cr
0 , & otherwise \cr}.
$$
Note that the case $a=b=c$ (which corresponds to $\epsilon =2$ in 10.11) here does not occur. Namely we now consider the range $a < (g-\lambda -1)/2$: the relation $a=b=c$ would contradict  $(R_1)$.

\medskip
\proclaim
Theorem 10.15. Let $g,\lambda,a,b$ be positive integers satisfying
$(R_1)$, $(R_2)$, $(R_3)$ and $c=g-3-a-b$.
If $a < (g-\lambda -1)/2$ then 
${\cal M}^{\lambda}_g(a,b)$ is an irreducible
variety of dimension $2(a+b)+\lambda +8-\epsilon-\xi$.
\par
\pf
Using 10.14, we can slightly
modify the construction in 10.11; essentially we use 
$\widetilde{\cal W}^{\lambda}_g(a,b)$ instead of
${\cal W}^{\lambda}_g(a,b)$. In particular, we consider models
$X' \subset  R_{1,1}$ of $X$ such that
$\aa \sim l$ and $\aa \subset \dd \sim 2l+(\lambda-2-c)l'$. Namely, if $\aa \not
\subset \dd$, then $\delta_A \le \aa \cdot \dd$; but $\delta_A = \lambda -2-a$
(from 10.13) while 
$ \aa \cdot \dd = \lambda -2 - c$ and this is impossible since $a < c$. \hb
Setting
$\widetilde{Y}^{\lambda}_g(a,b)$ the image of 
$\widetilde{\cal W}^{\lambda}_g(a,b)$ via the map 
$\Psi: {\cal W}^{\lambda}_g \rightarrow  \sym^{\delta} (R_{1,1})$
we have 
$$
\widetilde{Y}^{\lambda}_g(a,b) =
\{(P_1 , \dots , P_{\delta}) \;| \; 
\hbox{there exist}\; \aa \in |l|, \bb \in |l+(\lambda-2-c)l'|:
P_1 ,\dots , P_{\lambda -2-a} \in \aa, 
P_{\lambda -1-a},\dots , P_{\delta} \in \bb \}
$$
and the analogous of $(33)$ holds:
$$
\dim(\widetilde{\cal W}^{\lambda}_g(a,b))=
\dim(\widetilde{Y}^{\lambda}_g(a,b))+5 \lambda+4-3\delta.
\eqno (37)
$$
Consider the following correspondence
$$
Z_{a,b}^{\lambda}
\quad \subset \quad
|l| \times |l+(\lambda-2-c)l'| \times \sym^{\delta}(R_{1,1})
$$
defined by:
$$
Z_{a,b}^{\lambda}:=\left\{(\aa,\bb,(P_1 , \dots , P_{\delta})) 
\;  | \;
P_1 ,\dots , P_{\lambda -2-a} \in \aa, 
P_{\lambda -1-a},\dots , P_{\delta} \in \bb 
\right\}.
$$
Note that $b$ is determined from $a$ and $c$. Consider now the two canonical
projections:
$$
\matrix{
& Z_{a,b}^{\lambda} &
\cr
 &{}_p \swarrow{}\searrow {}_q &\cr
|l| \times |l+(\lambda-2-c)l'| &{}&
\widetilde{Y}^{\lambda}_g(a,b) \subset \sym^{\delta}(R_{1,1}) \hfil  
\cr }.
$$
With the same argument as in 10.11, one can see that the fibres
of $q$ are irreducible of dimension $\epsilon$. Note that, in this
case, $\epsilon$ can assume only the values $0$ and $1$, since
the assumption $a<(g-\lambda-1)/2$ implies $a<b$, 
otherwise $a+b <g-\lambda-1$, against $(R_3)$ (see  8.5). \hb  
Note that $p$ is surjective from 9.6. 
Moreover the general fibre $p^{-1}(\aa, \bb)$ of $p$ is isomorphic
to
$\sym^{\lambda-2-a}(\aa) \times \sym^{\delta-\lambda+2+a}(\bb)$,
hence it is irreducible of dimension $\delta$. Therefore we can
conclude that $Z_{a,b}^{\lambda}$ and hence
$\widetilde{Y}_{g}^{\lambda}(a,b)$ are irreducible and 
$$
\eqalign{
\dim(\widetilde{Y}_{g}^{\lambda}(a,b))
&=\dim(Z_{a,b}^{\lambda}) -\epsilon=
\dim|l| + \dim|l+(\lambda-2-c)l'| + \delta - \epsilon= \cr
&= 2(\lambda-1-c) + \delta - \epsilon
\cr}
$$
so, using $(37)$ we obtain
$$
\dim(\widetilde{\cal W}^{\lambda}_g(a,b)) = 
2(\lambda-1-c) + \delta - \epsilon
+5\lambda +4-3\delta=
2(3\lambda+1-c- \delta) + \lambda -\epsilon.
$$
Using $(31)$, we get  $3\lambda+1-c- \delta=3\lambda+1-(g-3-a-b)-3(\lambda-1) +g=
a+b+7$, so
$$
\dim(\widetilde{\cal W}^{\lambda}_g(a,b)) = 2(a+b)+14+\lambda -
\epsilon.
$$
\goodbreak
\noindent
Appliying 10.14, we obtain that 
$$
\dim({\cal M}^{\lambda}_g(a,b)) = 
\dim(\widetilde{\cal W}^{\lambda}_g(a,b)) -6-\xi= 
2(a+b)+8+\lambda - \epsilon-\xi
$$
as required. 
\cvd

\medskip
\noindent
{\bf Remark 10.16.} If $a < (g-\lambda-1)/2$ then $\delta = 3(\lambda-1)-g  >0$; in
particular,  $\lambda > (g+3)/3$. 
To show this, just remark that $g \le 3\lambda-3$ by
$(R_1)$; hence
$a < {{g-\lambda-1} \over 2} \le
{{3\lambda-3-\lambda-1} \over 2}= \lambda-2
$
so, from 10.13: $\delta \ge \delta_A = \lambda -2 -a >0$.

\medskip
\proclaim
Corollary 10.17. Set, as usual, $a \le b \le c$ and $a+b+c= g-3$. The following facts hold:
\item{1)} The general curve $X(g,\lambda,a,b)$ of  ${\cal M}^{\lambda}_g$ 
satisfies $a+b \ge (2g-8)/3$. 
\item{2)} For the general curve $X(g,\lambda,a,b)$ of  ${\cal M}^{\lambda}_g$, the values of
$a,b,c=g-3-(a+b)$ are determined by the class of $g$ (mod $3$); in
particular: 
$$
\matrix{
(i) \hfill & \hbox{ \sl if} & g=3p \hfill & \hbox{\sl then} \quad &
(a,b,c)=(p-1,p-1,p-1);
\cr (ii) \hfill & \hbox{\sl if} & g=3p+1 & \hbox{\sl then} \quad & (a,b,c)=(p-1,p-1,p);
\hfill
\cr (iii) & \hbox{\sl if} & g=3p+2 & \hbox{\sl then} \quad & (a,b,c)=(p-1,p,p). \hfill
\cr }
$$
\item{3)} Conversely, for the above values of $a$ and $b$ we obtain a stratum of maximal dimension, i.e.
$$
\dim({\cal M}^{\lambda}_g(a,b)) = \dim({\cal M}^{\lambda}_g).
$$
\smallskip
\noindent
Consequentely,
$$
\hbox{a curve} \quad
X(g,\lambda,a,b) \in {\cal M}^{\lambda}_g 
\quad \hbox{is general} 
\quad \Longleftrightarrow \quad
a,b,c \in \left\{ \left[ g-3 \over 3 \right] , \left[ g-1 \over 3 \right]   \right\}.
$$
\par
\pf  $1)$  We have to show that, if $a+b < (2g-8)/3$, then
$\dim({\cal M}^{\lambda}_g(a,b)) < \dim({\cal M}^{\lambda}_g)$. \hb
Let us rewrite the above condition correspondingly to the possible values of 
$g$ (mod $3$):
$$
\eqalign{
\bullet  \quad g=3p \phantom{{}+1} \quad& 
: \; a+b \le 2p-3 \quad \Rightarrow \quad a \le p-2 \quad \Rightarrow \quad  2a+b\le 3p-5; \cr
\bullet  \quad g=3p+1 \quad & 
: \; a+b \le 2p-3 \quad \Rightarrow \quad a \le p-2 \quad \Rightarrow \quad  2a+b\le 3p-5; \cr
\bullet  \quad g=3p+2 \quad & 
: \; a+b \le 2p-2 \quad \Rightarrow \quad a \le p-1 \quad \Rightarrow \quad  2a+b\le 3p-3. \cr
}
$$
Clearly, in all these cases
$$
a+b \le {{2g-9}\over3}
\qquad \hbox{and} \qquad
2a+b \le g-5.
\eqno{(38)}
$$
 {}From 10.11 (resp. 10.15) and using $(38)$ we immediately obtain:
$$
\eqalign{
a \ge {{g-\lambda -1} \over 2}
\; \Rightarrow \; &
\dim({\cal M}^{\lambda}_g(a,b)) \le
2(2a+b+\lambda)+10-g- \xi \le
2(g-5+\lambda) +10-g -\xi= g+2\lambda - \xi \cr
a < {{g-\lambda -1} \over 2}
\; \Rightarrow \; &
\dim({\cal M}^{\lambda}_g(a,b)) \le
2(a+b)+\lambda +8 - \xi \le
{{4g-18} \over 3} +\lambda +8- \xi=
g+\lambda+1 + {{g+3} \over 3} - \xi \cr
}
$$
where, in both cases, $\xi := \cases{
1, & if $\lambda={{g+3} \over 2}$ \cr
0 , & otherwise \cr}.$

\smallskip
\noindent
Note that, in the second case, from 10.16 we have that $(g+3)/3 < \lambda$. 
Therefore, for every value of $a$ it holds
$$
\dim({\cal M}^{\lambda}_g(a,b))< g+2\lambda+1 - \xi.
\eqno{(39)}
$$
Finally recall that, from 10.9, $\dim({\cal M}^{\lambda}_g)= g+2\lambda+1$ for all $(g+3)/3 < \lambda < \lambda_{\max}$, where $\lambda_{\max} = \left\lceil{{g+2} \over 2}\right\rceil$. \hb
On the other hand, from 10.3, 10.6 and 10.7 it turns out that
$$
\dim({\cal M}^{\lambda_{\max}}_g)= \cases{
g+2 \lambda_{\max}, & if $g$ odd \cr
g+2 \lambda_{\max}+1 , & if $g$ even \cr}.
$$
Therefore, if $\lambda < \lambda_{\max}$ or $g$ even, then $\xi =0$ so $(39)$ gives
$$
\dim({\cal M}^{\lambda}_g(a,b))< g+2\lambda+1 = \dim({\cal M}^{\lambda}_g).
$$
Otherwise, $\lambda = \lambda_{\max}$ and $g$ odd; then $\xi = 1$ so $(39)$ gives
$$
\dim({\cal M}^{\lambda}_g(a,b))< g+2\lambda = \dim({\cal M}^{\lambda}_g).
$$
and this proves the first part of the statement. 
\goodbreak

\noindent
$2)$ Let us consider  a general curve $X(g,\lambda,a,b) \in {\cal M}^{\lambda}_g$. We have just proved that $a+b \ge (2g-8)/3$. From the condition $(R_3)$ we get:
$$
{{2g-8} \over 3} \le a+b \le {{2g-6} \over 3}
\quad \Rightarrow \quad
a+b = \left[{{2g-6}\over 3}\right]
$$
hence $a+b$ is uniquely determined. Therefore, since $c=g-3-(a+b)$ and
$a \le b \le c$,  we obtain:
$$
 \matrix{
\bullet & g=3p \hfill & : & a+b= 2p-2 & \Rightarrow & c=p-1
& \Rightarrow & (a,b,c)=(p-1,p-1,p-1)\cr
\bullet & g=3p+1 & : & a+b= 2p-2 & \Rightarrow & c=p \hfill
& \Rightarrow & (a,b,c)=
\cases{(p-1,p-1,p) \cr (p-2,p,p) \cr}\hfill\cr
\bullet & g=3p+2 &: & a+b= 2p-1 & \Rightarrow & c=p \hfill
& \Rightarrow & (a,b,c)=(p-1,p,p) \hfill \cr
} 
$$
Note that the case $g=3p+1$ and $(a,b,c)=(p-2,p,p)$ does not
correspond to a general curve since, in this case,
$X(g,\lambda, a,b)$ belongs to a proper closed subset of ${\cal
M}^{\lambda}_g$. \hb
To show this, let us consider the two ranges of $a$ and the corresponding dimensions of the moduli spaces  found in 10.11 and 10.15, respectively.
\item {(I)} $a \ge {{g-\lambda -1} \over 2}$.
$$
\dim({\cal M}^{\lambda}_g(a,b)) \le
2(2a+b+\lambda)+10-g =
2(3p-4+\lambda) + 10 -(3p+1)=
3p+2\lambda+1
$$
while $\dim({\cal M}^{\lambda}_g) =g+2\lambda+1=3p+2\lambda+2$. 
\item {(II)}  $a< {{g-\lambda -1} \over 2}$. \hb
Substituting $g=3p+1$ in  $(R_1)$  and in the bound of $a$ in the assumption, we obtain respectively:
$$
\eqalign{
\lambda \ge {{g+3}\over 3} = p+{4\over 3}
& \quad \Rightarrow \quad 
\lambda \ge p+2 \cr
p-2= a<{{g-\lambda-1}\over 2 } 
& \quad \Rightarrow \quad 
 \lambda \le p+3. \cr}
$$
Using 10.15, under the assumption $(a,b,c)=(p-2,p,p)$ we obtain that $\epsilon =1$ and $\xi=0$, hence
$$
\dim({\cal M}^{\lambda}_g(a,b))=
2(a+b)+\lambda +8-\epsilon-\xi =
2(2p-2) + \lambda +8 - 1 = 4p + \lambda +3.
$$
On the other hand
$$
\dim({\cal M}^{\lambda}_g) =g+2\lambda+1 = 3p + 2 \lambda +2.
$$
Examining the two possible cases of $\lambda$, we immediately get:
$$
\dim({\cal M}^{\lambda}_g(a,b))= \cases{
5p+5, & if $\lambda= p+2$ \cr
5p+6, & if  $\lambda = p+3$ \cr}
\quad \hbox{while} \quad
\dim({\cal M}^{\lambda}_g)= \cases{
5p+6, & if $\lambda= p+2$ \cr
5p+8, & if  $\lambda = p+3$ \cr}
$$
and this proves the second part.

\smallskip
\noindent
$3)$ We are left to show that the strata corresponding to the values $(i), (ii), (iii)$ of $(a,b,c)$ are maximal. \hb
First note that the inequalities 
$a < {{g-\lambda-1} \over 2}$ and $\lambda \ge {g+3 \over 3}$ (the latter coming from $(R_1)$) become, respectively:
$$
\eqalign{
(i) & \qquad p-1 < {3p-\lambda -1 \over 2} \quad \hbox{and} \quad  \lambda \ge {3p+3 \over 3} \cr
(ii) & \qquad p-1 < {3p-\lambda  \over 2} \quad \hbox{and} \quad  \lambda \ge {3p+4 \over 3} \cr
(iii) & \qquad p-1 < {3p-\lambda+1  \over 2} \quad \hbox{and} \quad  \lambda \ge {3p+5 \over 3} \cr
}
$$
and in cases $(i)$ and  $(ii)$ we get a contraddiction, while in $(iii)$ we get $\lambda =p+2$. So
in cases $(i)$ and  $(ii)$ necessarily $a \ge {{g-\lambda-1} \over 2}$. \hb
Secondly, observe that if $a \ge {{g-\lambda-1} \over 2}$ then 10.11 can be applied and we have
$$
\dim({\cal M}^{\lambda}_g(a,b)) = 2(2a+b+\lambda)+10-g-\epsilon-\tau-\xi \eqno(*)
$$
where  $\xi=1$ if and only if $\lambda = {g+3 \over 2}$. This happens if $g$ is odd, so $\lambda= {g+3 \over 2} =
\left\lceil{{g+2} \over 2}\right\rceil$. Keeping the notation and the result in 10.7, where $\lambda_{\max}:=\left\lceil{{g+2} \over 2}\right\rceil$, we have that
$\dim({\cal M}_{g}^{\lambda_{\max}}) = 2g+3 = \dim({\cal M}_{g,4})$. Otherwise, $\xi = 0$  and $\lambda < \left\lceil{{g+2} \over 2}\right\rceil$; in this case, from 10.9, $\dim({\cal M}_{g}^{\lambda}) = g+2\lambda +1$.  \hb
Consider now each possibility.

\goodbreak
\smallskip
\noindent
Case $(i)$: $g=3p$, $(a,b,c) = (p-1,p-1,p-1)$. \hb
Since $\epsilon =2$ and $\tau =1$, from $(*)$ we obtain:
$$
\eqalign{
\dim({\cal M}^{\lambda}_g(a,b))
& = 2(3p-3+\lambda) + 10 -3p -2 -1 -\xi=  3p + 2\lambda +1 -\xi = \cr
& = g + 2\lambda +1 -\xi. \cr}
$$
Therefore
$$
\eqalign{
\lambda = \left\lceil{{g+2} \over 2}\right\rceil \Rightarrow & \qquad
 \xi = 1 \quad \hbox{and} \quad 
 \dim({\cal M}^{\lambda}_g(a,b)) = g + 2\lambda = g + 2 \; {g+3 \over 2} = 2g+3 = \dim({\cal M}^{\lambda}_g); \cr
\lambda < \left\lceil{{g+2} \over 2}\right\rceil \Rightarrow & \qquad
 \xi = 0 \quad \hbox{and} \quad 
\dim({\cal M}^{\lambda}_g(a,b)) =  g + 2\lambda +1 = \dim({\cal M}^{\lambda}_g). \cr
}
$$

\smallskip
\noindent
Case $(ii)$: $g=3p+1$, $(a,b,c) = (p-1,p-1,p)$. \hb
Since $\epsilon =0$ and $\tau =1$, from $(*)$ we again obtain:
$$
\dim({\cal M}^{\lambda}_g(a,b)) = 
3p+2 \lambda +2 - \xi = g + 2\lambda +1 -\xi.
$$
With the same argument as before we prove the claim.

\smallskip
\noindent
Case $(iii)$: $g=3p+2$, $(a,b,c) = (p-1,p,p)$. 
\item {I)}  If $a \ge {{g-\lambda-1} \over 2}$, the proof runs as above, using $(*)$ where $\epsilon =1$ and $\tau =0$.
\item {II)} If $a < {{g-\lambda-1} \over 2}$, the dimension of the strata is computed in 10.15 where one can find that
$$
\dim \left({\cal M}^{\lambda}_g(a,b)\right) = 2(a+b)+\lambda +8-\epsilon-\xi. \eqno(**)
$$
In our situation, $\epsilon = 1$ and $\xi =0$, since $\lambda \ne {g+3 \over 2}$  being $g=3p+2$ and $\lambda = p+2$, as remarked before. So $(**)$ gives $\dim \left({\cal M}^{\lambda}_g(a,b)\right) = 5p+7$. On the other hand, 
$\dim({\cal M}^{\lambda}_g) = g + 2\lambda +1 = 5p+7$.

\smallskip
\noindent
The final claim comes from  $(2)$ and $(3)$, together with a straightforward computation on the values in  $(i)$, $(ii)$, $(iii)$, taking into account that $a \le b \le c$.
 \cvd

\bigskip
\noindent
{\bf 11.  Moduli spaces of $4$--gonal curves with $t\ge 1$}
\medskip
Let us recall that if $t\ge 1$ and the double points of the standard model $\xx_0$ are distinct, then 
the bounds of the invariants $\lambda$ and $t$ are described in 5.4 $(i)-(iv)$ while
the invariants $a$ and $b$ are determined by
$\lambda$ and $t$ (see 5.1). More precisely, 
$$
{{g+3} \over 3} + t 
\le \lambda \le {{g+3} \over 2} 
\; , \quad
1 \le t \le {{g+3} \over 6}
$$

$$
a=g-2\lambda+t+1, \quad
b=\lambda-t-2, \quad
c=\lambda-2.
$$
  As a consequence,  the subvariety of ${\cal
W}^{\lambda}_g$ parametrizing the curves of invariants $g,\lambda,t,a,b$ can
be simply denoted by 
${\cal W}^{\lambda}_g(t)$. 

In order to describe such variety,  we perform a construction similar to that in
10.11. 

Let us denote by ${\cal A}_\lambda^t$ the open subset of the linear
system $|4C_0+(\lambda +t)f|$ on $R_{1,t+1}$ parametrizing the
irreducible curves of such linear system and set
$$
{\cal W}^{\lambda}_g(t):= \{X' \in {\cal A}_\lambda^t \; | \;
X =X(g,\lambda,t) \; \hbox{and it has $\delta$ distinct double points on $C_0$} \}.
$$
If we consider the morphism
$$
\varphi:= \varphi_{4C_0+(\lambda+t)f}: \quad
R_{1,t+1} \longrightarrow S' \subset \Bbb P^N
$$
  it is clear that $N=h^0(R_{1,t+1},{\cal
O}_{R_{1,t+1}}(4C_0+(\lambda+t)f))-1=5(\lambda -t)+4$ (from [{\bf 4}], 1.8) and we can
identify ${\cal W}^{\lambda}_g(t)$ with the following subset of
$\check{\Bbb P}^N$:
$$
{\cal W}^{\lambda}_g(t) \cong
\{H \in \check{\Bbb P}^N \; | \; H \supset \langle T_{P_1}(S'), \dots
T_{P_\delta}(S') \rangle, P_i \in C_0
\}.
$$
Therefore, consider the following correspondence
$$
\widetilde W=
\{(H; P_1, \dots, P_\delta) \; | \; H \supset \langle
T_{P_1}(S'),
\dots T_{P_\delta}(S') \rangle
\} \subset \check{\Bbb P}^N \times \sym^\delta(\Bbb P^1)
$$
and the projections
$$
\matrix{
{} & {} & \widetilde W & {} & {} \cr
{} & {}^{\pi_1}\swarrow & {} & \searrow^{\pi_2} & {} \cr
\check{\Bbb P}^N & & {} & {} & {} \sym^\delta(\Bbb P^1)
\cr}
$$
Obviously, $\pi_1(\widetilde W)=\overline{{\cal W}^\lambda_g(t)}$  and
$\pi_1$ is an isomorphism on an open subset of ${\cal W}^\lambda_g(t)$.
Moreover, $\pi_2$ is surjective and the fibres have dimension 
$N- \dim \langle T_{P_1}(S'), \dots T_{P_\delta}(S') \rangle$. \hb
One can show (as in 9.4) that also in the case $t \ge 1$ it holds that the
space  $\langle T_{P_1}(S'), \dots T_{P_\delta}(S') \rangle$
has maximum dimension, i.e. $3\delta -1$.  Hence
$\dim({\cal W}^\lambda_g(t))= \dim \widetilde W = 
N - (3\delta -1) +\delta= 5(\lambda -t+1)-2 \delta$, so using 2.2 $(iii)$, 
we obtain:
$$
\dim({\cal W}^\lambda_g(t)) =
2g+t-\lambda+11.
$$ 
As well as in the case $t=0$, one can show that these varieties are not empty.
Furthermore, let us recall that the automorphism group of a rational ruled
surface $R_{1,t+1} \subset \Bbb P^{t+2}$ has dimension $t+5$, if $t \ge 1$,
and 6, if $t=0$ (as we already noted in 10.3).  These two facts, together with the
previous computation of $\dim({\cal W}^\lambda_g(t))$, immediately give the
following result:

\medskip
\proclaim
Theorem 11.1. Let $g,\lambda, t$ be positive integers satisfying: $g \ge 10$, 
$$
{{g+3} \over 3} + t 
\le \lambda \le {{g+3} \over 2} 
\; , \quad
1 \le t \le {{g+3} \over 6}.
$$
Then
${\cal  M}_g^\lambda(t)$ is an irreducible variety of dimension
$2g-\lambda+6$.
\cvd
\par

\goodbreak

\bigskip
\noindent
{\bf References }
\medskip
\item {[{\bf 1}]} E. Arbarello, M. Cornalba, `Footnotes to a Paper of Beniamino Segre',
{\it Math. Ann.} 256 (1981) 341-362

\item {[{\bf 2}]} E. Arbarello, M. Cornalba, P.A. Griffiths, J.
Harris, {\it Geometry of Algebraic Curves} I (Springer--Verlag, New
York, 1985)

\item {[{\bf 3}]} L. Badescu, {\it Algebraic Surfaces}, Universitext (Springer--Verlag, New
York, 2001)

\item {[{\bf 4}]}M. Brundu, G. Sacchiero, `On the varieties
parametrizing trigonal curves with assigned Weierstrass points', 
{\it Comm. in Alg.}  26 N.10 (1998) 3291-3312 

\item {[{\bf 5}]} M. Brundu, G. Sacchiero, `On rational surfaces ruled by conics',
{\it Comm in Alg.} 31 N.8 (2003) 3631-3652

\item {[{\bf 6}]} M. Brundu, G. Sacchiero, `On the singularities of surfaces ruled by conics', {\it to appear}

\item {[{\bf 7}]} R. Friedman, {\it Algebraic surfaces and holomorphic vector bundles},
Universitext (Springer--Verlag, New York, 1998)

\item {[{\bf 8}]} W. Fulton, {\it Intersection Theory} (Springer--Verlag, 
New York, 1998)

\item {[{\bf 9}]} J. Harris, ` The Genus of Space Curves',  {\it Math. Ann.} 249  (1980) 191-204 

\item {[{\bf 10}]} J. Harris, {\it Algebraic Geometry}, Graduate
Texts in Mathematics 133 (Springer--Verlag, New York, 1992)

\item {[{\bf 11}]} R. Hartshorne, {\it Algebraic Geometry}, Graduate
Texts in Mathematics 52 (Springer--Verlag, New York, 1977)

\item {[{\bf 12}]}  F.O. Schreyer,  `Syzygies of Canonical Curves and Special Linear Series',  {\it Math. Ann.}  275 (1986) 105-137

\item {[{\bf 13}]} B. Segre, `Sui moduli delle curve poligonali e sopra un complemento al
teorema di esistenza di Riemann', {\it Math. Ann.} 100 (1928) 537-551

\item {[{\bf 14}]} K.O. St\"orh, P. Viana, `Weierstrass gap sequences and moduli
varieties of trigonal curves', {\it J. Pure Appl. Algebra} 81 (1992) 63-82

\par

\bye